\documentclass[mnsc,nonblindrev]{informs3}

\OneAndAHalfSpacedXI 

\usepackage{algorithm}
\usepackage{algorithmic}
\usepackage{bm}

\newcommand{\R}{\mathbb{R}}
\newcommand{\norm}[1]{\left \Vert #1 \right\Vert}
\newcommand{\abs}[1]{\left \vert #1 \right\vert}

\newcommand{\prth}[1]{\left(#1\right)}

\usepackage{natbib}
 \bibpunct[, ]{(}{)}{,}{a}{}{,}%
 \def\bibfont{\small}%
\TheoremsNumberedThrough     
\ECRepeatTheorems

\EquationsNumberedThrough    

\begin{document}


\RUNAUTHOR{Lam and Qian}

\RUNTITLE{Combating Conservativeness}

\TITLE{Combating Conservativeness in Data-Driven Optimization under Uncertainty: A Solution Path Approach}

\ARTICLEAUTHORS{%
\AUTHOR{Henry Lam}
\AFF{Department of Industrial Engineering and Operations Research, Columbia University, New York, NY 10027, \EMAIL{henry.lam@columbia.edu}} 
\AUTHOR{Huajie Qian}
\AFF{Department of Industrial Engineering and Operations Research, Columbia University, New York, NY 10027, \EMAIL{h.qian@columbia.edu}}
} 

\ABSTRACT{%
In data-driven optimization, solution feasibility is often ensured through a ``safe" reformulation of the uncertain constraints, such that an obtained data-driven solution is guaranteed to be feasible for the oracle formulation with high statistical confidence. Such approaches generally involve an implicit estimation of the whole feasible set that can scale rapidly with the problem dimension, in turn leading to over-conservative solutions. In this paper, we investigate a validation-based strategy to avoid set estimation by exploiting the intrinsic low dimensionality among all  possible solutions output from a given reformulation. We demonstrate how our obtained solutions satisfy statistical feasibility guarantees with light dimension dependence, and how they are asymptotically optimal and thus regarded as the least conservative with respect to the considered reformulation classes. We apply this strategy to several data-driven optimization paradigms including (distributionally) robust optimization, sample average approximation and scenario optimization. Numerical experiments show encouraging performances of our strategy compared to established benchmarks.


}%


\KEYWORDS{optimization under uncertainty, data-driven optimization, stochastic constraint, chance constraint, validation, dimension reduction}

\maketitle

%

\section{Introduction}\label{sec:intro}

We focus on optimization with stochastic or probabilistic constraints that, on a high level, can be written in the form
\begin{equation}
\min_x f(x)\text{\ \ subject to\ \ }H(x)\geq\gamma\label{stoc_opt} 
\end{equation}
where $H(x)$ is an expectation of a random function of the decision variable $x$. Formulation \eqref{stoc_opt} is ubiquitous in decision-making problems under multiple tradeoffs, where the constraint $H(x)\geq\gamma$ signifies a restriction on the risk level or resource capacity (e.g., \cite{atlason2004call,krokhmal2002portfolio}). Moreover, when the random function is an indicator of an event, formulation \eqref{stoc_opt} is a so-called probabilistically constrained or chance-constrained problem (CCP) (\cite{prekopa2003probabilistic}). This important formulation posits the decision to tolerate a small probability on catastrophic events such as system failures and big losses, and comprises a primary approach for safe decision-making when facing uncertainty.




We are interested in the situation where the probability distribution governing $H$ is unknown but only observed through data. Finding good solutions under this setting has been studied prominently in the data-driven optimization literature, harnessing various tools from (distributionally) robust optimization (e.g., \cite{bertsimas2011theory,ben2009robust,wiesemann2014distributionally}) to sample average approximation (e.g., \cite{shapiro2009lectures}) and scenario optimization (e.g., \cite{campi2008exact}). From a statistical viewpoint, the problem challenge and the focus of these studies can be cast as a balancing between feasibility and optimality. Due to data noise, feasibility is at best guaranteed with a high statistical confidence, and accounting for this uncertainty incurs a price on the achieved objective value -- resulting in conservativeness. This impact on optimality from ensuring feasibility depends heavily on the efficiency in assimilating statistical information into the data-driven formulation. In the following, we first explain how the established estimation frameworks can face severe ``looseness" in this regard and lead to over-conservative solutions. This motivates our study that, on a high level, aims to investigate a strategy to substantially tighten the feasibility-optimality tradeoff compared to the previous methods.

\subsection{Existing Frameworks and Motivation of Our Approach}\label{sec:existing}

To facilitate discussion, suppose for concreteness that the decision variable $x$ lies in a $d$-dimensional deterministic space $\mathcal X\subset\mathbb R^d$. Denote $H(x):=\mathbb E_{F}[h(x,\xi)]$ where $\mathbb E_F[\cdot]$ is the expectation under  $\xi\sim F$, and $h(\cdot,\cdot):\R^d\times \R^m\to\R$ is a function of $x\in\mathcal X$ controlled by the randomness $\xi\in\mathbb R^m$. Also, since our focus is on handling uncertain constraints, we assume that the objective function $f$ is deterministic (this can be relaxed with proper modifications of our subsequent discussion). Suppose we have i.i.d. observations $\xi_1,\ldots,\xi_n$.

Let us first consider a natural idea to replace the unknown $H(\cdot)$ with some point estimate, say the sample average $\hat H(\cdot)=(1/n)\sum_{i=1}^nh(x,\xi_i)$, in the constraint. Though simple, this approach is typically inadequate to ensure feasibility in any statistical sense. To explain, suppose the ``true" optimal solution $x^*$ is at the boundary of the feasible region, i.e., $H(x^*)=\gamma$. If we use $(1/n)\sum_{i=1}^nh(x,\xi_i)\geq\gamma$ as the constraint, then, with significant probability an obtained solution $\hat x^*$ (conceivably also at the boundary of the data-driven constraint) can have $H(\hat x^*)$ below $\gamma$  (when $(1/n)\sum_{i=1}^nh(\hat x^*,\xi_i)>\mathbb E_F[h(\hat x^*,\xi)]$), which is infeasible for the original problem. This issue may not arise if $x^*$ or $\hat x^*$ is in the interior of the feasible region, but a priori we do not know our decision. In other words, the nature of constrained optimization enforces us to put some ``safety" margin in addition to the point estimate, in order to achieve any reasonable confidence in feasibility. Here, we can plausibly use a data-driven constraint $\hat H(x)-\epsilon(x)\geq\gamma$, where $\epsilon(x)$ is a properly chosen positive function such that $H(x)\geq\hat H(x)-\epsilon(x)$ for any $x$ with high confidence (such as the scheme in \cite{wang2008sample}, among others).





We place the above discussion in a more general framework. Let $\mathcal F$ be the (unknown) feasible region of \eqref{stoc_opt}. Given the data $\xi_1,\ldots,\xi_n$, by a valid procedure we mean one that is able to output a solution $\hat x^*$ that is truly feasible with a given high confidence level, say $1-\beta$ (e.g., $95\%$). That is, \begin{equation}
P_{data}(\hat x^*\in\mathcal F)\geq1-\beta\label{conf guarantee}
\end{equation}
where $P_{data}$ refers to the probability with respect to the data. By a data-driven reformulation, we mean replacing $\mathcal F$ with $\hat{\mathcal F}$ that is constructed solely from the data $\xi_1,\ldots,\xi_n$. This gives
\begin{equation}
\min_{x\in\mathcal X}\ f(x)\text{\ \ subject to\ \ }x\in\hat{\mathcal F}\label{data_opt}
\end{equation}which outputs solution $\hat x^*$.
If we can choose $\hat{\mathcal F}$ such that
\begin{equation}
P_{data}(\hat{\mathcal F}\subset\mathcal F)\geq1-\beta\label{feasibility guarantee}
\end{equation}then we clearly have \eqref{conf guarantee} since $P_{data}(\hat x^*\in\mathcal F)\geq P_{data}(\hat{\mathcal  F}\subset\mathcal F)$. In the example above, we have used $\hat{\mathcal F}=\{x\in\mathcal X:\hat H(x)-\epsilon(x)\geq\gamma\}$, in the hope that \eqref{feasibility guarantee} holds in order to achieve \eqref{conf guarantee}.

We contend that most approaches in data-driven optimization rely on the above reasoning and are based on \eqref{feasibility guarantee}. In particular, \eqref{feasibility guarantee} provides a convenient way to certify feasibility, by requiring that \emph{all} solutions feasible for \eqref{data_opt} are also feasible for \eqref{stoc_opt} with high confidence. This set-level guarantee generally hinges on a simultaneous estimation task across all $x$ in the decision space $\mathcal X$, for which a proper control of the statistical error can lead to a substantial shrinkage of the size of $\hat{\mathcal F}$ that exacerbates with problem dimension (either of the decision space or the probability space).




We provide several examples to illustrate the phenomenon above. Some of these examples apply most relevantly to CCP, where $H(x)$ is in the form $\mathbb P_F(G(x,\xi)\leq b)$ with $G(x,\xi):\mathbb R^d\times\mathbb R^m\to\mathbb R$.




\begin{example}[Sample average approximation (SAA)]
In the case of CCP, the SAA approach sets $\hat{\mathcal F}=\{x\in\mathcal X:
    \frac{1}{n}\sum_{i=1}^n\mathbf{1}(G(x,\xi_i)+\epsilon\leq b)\geq\gamma+\delta\}$,
where $\epsilon$ and $\delta$ are suitably tuned parameters. For example, when $G$ is Lipschitz continuous in $x$, selecting $   \delta= \Omega(\sqrt{(d/n)\log (1/\epsilon)})$ can guarantee \eqref{feasibility guarantee} (\cite{luedtke2008sample}), and similar relations also hold in discrete decision space (\cite{luedtke2008sample}) and expected value constraints (\cite{wang2008sample}). These estimates come from concentration inequalities in which union bounds are needed and give rise to the dependence on the dimension $d$. Note that the resulting margin $\delta$ scales in order $\sqrt d$, and to get any reasonably small $\delta$, $n$ must be of higher order than $d$.
\hfill\halmos
\end{example}

\begin{example}[Robust optimization (RO) and safe convex approximation (SCA)]\label{example:RO and SCA}
Focusing on CCP,  RO sets \begin{equation}
\hat{\mathcal F}=\{x\in\mathcal X:
   G(x,\xi)\leq b,\text{ for all }\xi\in\mathcal U\}
      \label{RO}
\end{equation}
where $\mathcal U$ is known as the uncertainty set, and $\xi$ in \eqref{RO} is viewed as a deterministic unknown (\cite{bertsimas2011theory,ben2009robust}). A common example of $\mathcal U$ is an ellipsoidal set $\{\xi:(\xi-\hat\mu)'\hat\Sigma^{-1}(\xi-\hat\mu)\leq \rho\}$ where $\hat\mu\in\mathbb R^d$, $\hat\Sigma\in\mathbb R^{d\times d}$ a positive semidefinite matrix, and $\rho\in\mathbb R$. Here the center $\hat\mu$ and shape $\hat\Sigma$ typically correspond to the mean and covariance of the data, and $\rho$ controls the set size. A duality argument shows that, in the case of linear chance constraint in the form $G(x,\xi)=x'\xi$, \eqref{RO} is equivalent to the quadratic constraint $\hat\mu'x+\sqrt \rho\|\hat\Sigma^{1/2}x\|_2\leq b$. Using such type of convex constraints as inner approximations for intractable chance constraints is also known as SCA (e.g., \cite{ben2000robust,nemirovski2003tractable,nemirovski2006convex}).

It is known that, if for instance the random variable $\xi$ has a known bounded support, the above approach guarantees an obtained solution has a satisfaction probability of order $1-e^{-\rho/2}$ via Hoeffding's inequality, and $\rho$ is chosen by matching this expression with the tolerance level $\gamma$. Although $\rho$ calibrated this way may not explicitly depend on the problem dimension, its tightness varies heavily based on problem instance (due to the worst-case nature of concentration bounds), and its validity relies on a priori distributional information (e.g., support) rather than an efficient utilization of data.
Another viewpoint that has been taken recently in data-driven RO (\cite{bertsimas2018data}; \cite{tulabandhula2014robust}; \cite{goldfarb2003robust}; \cite{hong2017learning}) is to select $\mathcal U$ to be a set that contains $\gamma$-content of the distribution of $\xi$, i.e., $\mathbb P_F(\xi\in\mathcal U)\geq\gamma$, with a confidence level $1-\beta$. In this case, any solution $\hat x$ feasible for \eqref{RO} would satisfy $\mathbb P_F(G(\hat x,\xi)\leq b)\geq\mathbb P_F(\xi\in\mathcal U)\geq\gamma$ with at least $1-\beta$ confidence, thus achieving \eqref{feasibility guarantee} as well. Such generated uncertainty set however typically has a size that scales with the dimension of the probability space. For example, consider $G(x,\xi)=x'\xi$ with $\xi\in \R^m$ being standard multivariate Gaussian and the uncertainty set $\mathcal U$ is an ellipsoid with $\hat{\mu}$ and $\hat{\Sigma}$ being the true mean and covariance, i.e., $\mathcal U=\{\xi\in\R^m:\norm{\xi}_2^2\leq \rho\}$. Then, in order to make $\mathcal U$ a $\gamma$-content set the radius $\rho$ has to be at least of order $m$ since $\norm{\xi}_2^2$ has a mean $m$, resulting in the robust counterpart $\sqrt{\rho}\norm{x}_2=\Theta(\sqrt{m})\norm{x}_2 \leq b$. However, the exact chance constraint in this case can be rewritten as $z_{\gamma}\norm{x}_2 \leq b$, where $z_{\gamma}$ is the $\gamma$-quantile of the univariate standard normal, which is independent of the dimension.\hfill\halmos

%

%

\end{example}

\begin{example}[Distributionally robust optimization (DRO)]\label{example:DRO}
DRO sets
\begin{equation}\hat{\mathcal F}= \{x\in\mathcal X:   \inf_{Q\in\mathcal U}\mathbb E_Q[h(x,\xi)]\geq\gamma\}\label{DRO general constraint}
\end{equation}
where $\mathcal U$ is a set in the space of probability measures that is constructed from data, and is often known as the ambiguity set or uncertainty set. The rationale here is similar to RO, but views the uncertainty in terms of the distribution. If $\mathcal U$ is constructed such that it contains the true distribution $F$ with high confidence, i.e., $P_{data}(F\in\mathcal U)\geq1-\beta$, then any solution $\hat x$  feasible for the DRO constraint \eqref{DRO general constraint} would satisfy $\mathbb P_F(G(\hat x,\xi)\leq b)\geq\gamma$ with at least $1-\beta$ confidence so that \eqref{feasibility guarantee} holds.

Popular choices of $\mathcal U$ include moment sets, i.e., specifying the moments of $Q$ (to be within a range for instance) (\cite{ghaoui2003worst}; \cite{delage2010distributionally}; \cite{xu2012optimization}; \cite{wiesemann2014distributionally}; \cite{goh2010distributionally}; \cite{natarajan2008incorporating}; \cite{van2016generalized}; \cite{doan2015robustness}; \cite{hanasusanto2015distributionally}), and distance-based sets, i.e., specifying $Q$ in the neighborhood ball surrounding a baseline distribution, where the ball size is measured by a statistical distance such as $\phi$-divergence (\cite{petersen2000minimax}; \cite{ben2013robust}; \cite{glasserman2014robust}; \cite{lam2016robust}; \cite{lam2017sensitivity}; \cite{hu2013kullback}; \cite{jiang2016data}; \cite{gotoh2018robust}; \cite{dupuis2016path}; \cite{bayraksan2015data}) or Wasserstein distance (\cite{esfahani2018data}; \cite{blanchet2019quantifying}; \cite{gao2016distributionally}; \cite{xie2018distributionally}).

Ensuring $P_{data}(F\in\mathcal U)\geq1-\beta$ means that $\mathcal U$ is a confidence region for $F$. In the moment set case, this boils down to finding confidence regions for the moments whose sizes in general scale with the probability space dimension. To explain, when only the mean $\mathbb E_{F}[\xi]$ is estimated, the confidence region constructed from, say the delta method (\cite{marandi2019extending}), takes the form $\{\hat{\mu}+\hat{\Sigma}^{\frac{1}{2}}v:v\in\R^m,\norm{v}_2^2\leq \chi_{m,1-\beta}^2\}$, where $\hat{\mu}$ and $\hat{\Sigma}$ are the sample mean and covariance and $\chi_{m,1-\beta}^2$ (which is of order $m$) is the $1-\beta$ quantile of the $\chi^2$ distribution with degree of freedom $m$, therefore the diameter of the confidence region scales as $\sqrt{m}$. When the mean and covariance are jointly estimated, the dimension dependence scales up further. In the distance-based set case, one needs to estimate statistical distances. If the Wasserstein distance is used to construct the ball surrounding the empirical distribution, results from measure concentration (\cite{fournier2015rate}) indicate that the ball size needs to be of order $n^{-\frac{1}{m}}$ to ensure $P_{data}(F\in\mathcal U)\geq1-\beta$. Alternatively, if $\mathcal U$ is constructed as a $\phi$-divergence ball surrounding some nonparametric kernel-type density estimate, results from kernel density estimation (see Section 4.3 in \cite{wand1994kernel}) suggest that the estimation error is of order $n^{-\frac{4}{m+4}}$. In either case, the required size of the uncertainty set exhibits exponential dependence on the dimension. Recently, the empirical or the profile likelihood method has also been proposed to calibrate the ball size such that $\mathcal U$ can be (much) smaller than what is needed in being a confidence region for $F$, while at the same time \eqref{feasibility guarantee} still holds (\cite{lam2017empirical}; \cite{duchi2016statistics}; \cite{lam2016recovering}; \cite{blanchet2016sample}). However, the ball size in this approach scales as the supremum of a so-called $\chi^2$-process over the decision space (e.g., \cite{lam2016recovering}). An analysis using metric entropy (e.g., Example 2 in Section 14 in \cite{lifshits2013gaussian}) shows that the $\chi^2$-process supremum can scale linearly in the decision space dimension $d$, a much better but still considerable dependence on the dimension.\hfill\halmos
\end{example}

Finally, we discuss the only two exceptional paradigms, to our best knowledge, in providing guarantee \eqref{conf guarantee} using \eqref{feasibility guarantee}. First, \cite{gupta2019near} studies a Bayesian framework to define feasibility guarantees for (stochastic) constraints with unknown parameters, focusing on DRO formulations. The idea is to ensure the obtained data-driven solution satisfies the constraints with a high posterior probability on the unknown parameters. This definition of feasibility does not utilize the concept of experimental repetitions in the frequentist sense as we have considered, but views the unknown parameters as random and considers the frequency of feasibility from the posterior belief, thus bypassing the set-level guarantee in \eqref{feasibility guarantee}. Indeed, \cite{gupta2019near} shows that under suitable convexity assumption on the constraints (with respect to the unknown parameter) and discreteness of the underlying distribution, the size of the uncertainty set in DRO can be chosen lightly dependent on the problem dimension.

The second exceptional paradigm that we are aware of is scenario optimization (SO) (e.g., \cite{calafiore2005uncertain,campi2008exact}), which applies to the case of CCP. In its basic form, this approach sets $$\hat{\mathcal F}=\{x\in\mathcal X:
   G(x,\xi_i)\leq b\text{\ \ for all\ \ }i=1,\ldots,n\}$$
i.e., using sampled constraints formed from the data. As the number of constraints increases, $\hat{\mathcal F}$ is postulated to populate the decision space in some sense and ensure the obtained solution $\hat x^*$ lies in $\mathcal F$. While the sample size required in the basic SO is linear in the decision dimension $d$, recent works reduce this dependence by an array of generalizations, including using regularization (\cite{campi2013random}), tighter support rank estimates (\cite{schildbach2013randomized,campi2018wait}) and sequential and validation-based schemes (\cite{care2014fast,calafiore2017repetitive}).

The approach that we propose in this paper aims to avoid using the set-level guarantee in \eqref{feasibility guarantee} and the need to control its simultaneous estimation errors, which can cause over-conservativeness as discussed. Our approach operates under a frequentist framework, nonparametric assumptions on the underlying distributions, and applies to all the exemplified methods mentioned above (SAA, RO, DRO and SO). It is thus different from the Bayesian parametric framework in \cite{gupta2019near}. Our idea is closest to some of the validation-type schemes suggested for SO, but more general as it applies to stochastic constraints beyond CCP and to data-driven reformulations beyond SO. Akin to these SO studies, our main results concern the power of our validation procedures in guaranteeing feasibility, which informs the required sample size in relation to the problem dimension. Our results also introduce a notion of optimality with respect to the chosen reformulation class, and deduce joint optimality-feasibility guarantees. In these regards, one main contribution of our work can be viewed as a rigorous construction of the first general-purpose validation framework for data-driven constrained optimization to systematically reduce conservativeness.

\section{Overview of Our Framework and Rationale}\label{sec:overview}
Our framework, as discussed, aims to bypass the set-level guarantee in \eqref{feasibility guarantee} and the need to control its simultaneous estimation errors. Our starting observation is the following. In all the described approaches above, the data-driven reformulation involves a key parameter that controls the level of conservativeness:
\begin{enumerate}
\item SAA: safety margin $\delta$
\item RO and SCA: uncertainty set size $\rho$
\item DRO: divergence ball size or moment set size
\item SO: number of constraints
\end{enumerate}
These parameters have the properties that setting it to one extreme (e.g., 0) would signal no uncertainty in the formulation, leading to a solution very likely infeasible, while setting it to another extreme (e.g., $\infty$) would cover the entire decision space, leading to a solution that is very conservative. In the established approaches, the parameter value is chosen to ensure \eqref{feasibility guarantee}, which tend to locate towards the latter extreme.


On the other hand, given a specific data-driven reformulation, it is easy to see that no matter how we choose this ``conservativeness" parameter, the solution must lie in a low-dimensional manifold. More precisely, denote a given data-driven reformulation as
\begin{equation}\label{para_opt}
\min_{x\in\mathcal X}\ f(x)\text{\ \ subject to\ \ }x\in\hat{\mathcal F}(s)
\end{equation}
where $s\in S$ denotes the conservativeness parameter, and we highlight the dependence of the data-driven feasible region $\hat{\mathcal F}(s)$ on $s$. We denote the obtained solution from \eqref{para_opt} as $x^*(s)$. The \emph{solution path} $\{x^*(s):s\in S\}$ contains all possible obtainable solutions from the data-driven reformulation \eqref{para_opt}. Intuitively, any statement on feasibility suffices to focus on this solution path, instead of the whole decision space.

Nonetheless, besides the conservativeness parameter, a data-driven reformulation could have other parameters playing various roles (e.g., center and shape of an ellipsoidal uncertainty set in RO, baseline distribution in distance-based DRO etc.). The flexibility of these parameter values can enlarge the obtainable solution space and elevate its dimensionality. Suppose we want to contain this enlargement, and at the same time be able to select the optimal candidate within the low-dimensional manifold $\{x^*(s):s\in S\}$. We propose the following two-phase framework to achieve this rigorously.

\begin{algorithm}
\caption{The Two-Phase Framework}
\label{two-phase framework}
\begin{algorithmic}
\STATE {\textbf{Input:} data $\bm\xi_{1:n}=\{\xi_1,\ldots,\xi_n\}$; numbers of data $n_1,n_2$ allocated to each phase ($n_1+n_2=n$); a confidence level $1-\beta$; a given method to construct data-driven reformulation with a (possibly multi-dimensional) parameter $s\in S$; a discrete mesh $\{s_1,s_2,\ldots,s_p\}\subseteq S$.}
\vspace{1ex}
\STATE\textbf{Phase one:}
\STATE \textbf{1.}~Use $n_1$ observations, which we index as $\{\xi_{n_2+1},\ldots,\xi_{n}\}$ for convenience, to construct the data-driven reformulation $OPT(s)$ in the form \eqref{para_opt} parameterized by $s\in S$.
\STATE \textbf{2.}~For each $j=1,\ldots,p$, compute the optimal solution $x^*(s_j)$ of $OPT(s_j)$.
\vspace{1ex}
\STATE \textbf{Phase two:}

Use a validator $V$ to select $(\hat s^*,x^*(\hat s^*))= V(\{\xi_1,\ldots,\xi_{n_2}\},\{x^*(s_1),\ldots,x^*(s_p)\},1-\beta)$, where $x^*(\hat s^*)$ is a solution and $\hat s^*$ is the associated parameter value.

\vspace{1ex}
\STATE {\textbf{Output:} $x^*(\hat s^*)$.}
\end{algorithmic}
\end{algorithm}

Our procedure (Algorithm \ref{two-phase framework}) splits the data into two groups. With the first group of data, we construct a given data-driven reformulation parametrized by a conservativeness parameter $s$ that varies over a space $S$, which we call $OPT(s)$. We obtain the optimal solution $x^*(s)$ for a range of values $s=s_j,j=1,\ldots,p$. This step assumes the availability of an efficient solver for $OPT(s)$. Next, the second group of data is fed into a \emph{validator} $V$ that aims to identify the best feasible solution $x^*(\hat s^*)$ among $\{x^*(s_j):j=1,\ldots,p\}$. The number of points $p$ required to validate depends on the size of $S$, which is constructed to be low-dimensional. There are multiple ways to set up the validator $V$, each with its own benefits and requirements. In the next two sections, we will introduce two classes of validators, one we call \emph{Gaussian supremum validator} (Section \ref{sec:gaussian supremum margin}), and another one we call \emph{univariate Gaussian validator} (Section \ref{sec:standard gaussian margin}). We will present their rationales, theoretical statistical guarantees, and implications on the feasibility and optimality of the obtained solution. Section \ref{sec:formulation showcase} will then tie back the applicability of these validators to the exemplified approaches in Section \ref{sec:existing}.

\section{Validation via Multivariate Gaussian Supremum}\label{sec:gaussian supremum margin}


Our first validator uses a simultaneous estimation of $H(x)$ in the constraint in \eqref{stoc_opt} to assess feasibility over the discretized solution path of $x^*(s)$. More precisely, given the solution set $\{x^*(s_j):j=1,\ldots,p\}$, we use a sample average with an appropriately calibrated safety margin, i.e., $\frac{1}{n_2}\sum_{i=1}^{n_2}h(x,\xi_i)-\epsilon$, to replace the unknown $H(\cdot)$ in \eqref{stoc_opt} and output the best solution among the set. The margin $\epsilon$ is calibrated via the limiting distribution of $(\frac{1}{n_2}\sum_{i=1}^{n_2}h(x^*(s_j),\xi_i))_{j=1,\ldots,p}$ which captures the estimation error of $H(\cdot)$ and is multivariate Gaussian. It contains a critical value $q_{1-\beta}$ that is the quantile of a Gaussian supremum. Algorithms \ref{algo:supremum_unnormalized} and \ref{algo:supremum_normalized} describe two variants of this validator, one unnormalized while another one normalized by the standard deviation at each $s_j$. In the following, we denote $N_p(0,\Sigma)$ as a $p$-dimensional Gaussian vector with mean zero and covariance $\Sigma$.

\begin{algorithm}
\caption{$V$: Unnormalized Gaussian Supremum Validator}
\label{algo:supremum_unnormalized}
\begin{algorithmic}
\STATE {\textbf{Input:} $\{\xi_1,\ldots,\xi_{n_2}\},\{x^*(s_1),\ldots,x^*(s_p)\},1-\beta$}
\vspace{1ex}
\STATE \textbf{1.}~For each $j=1,\ldots,p$ compute the sample mean $\hat H_j=(1/n_2)\sum_{i=1}^{n_2}h(x^*(s_j),\xi_i)$ and sample covariance matrix $\hat \Sigma$ with $\hat \Sigma(j_1,j_2) = (1/n_2)\sum_{i=1}^{n_2}(h(x^*(s_{j_1}),\xi_i)-\hat H_{j_1})(h(x^*(s_{j_2}),\xi_i)-\hat H_{j_2})$.
\STATE \textbf{2.}~Compute $q_{1-\beta}$, the $(1-\beta)$-quantile of $\max\{Z_1,\ldots,Z_p\}$ where $(Z_1,\ldots,Z_p)\sim N_p(0,\hat \Sigma)$, and let
\begin{equation}\label{validate:gs unnormalized}
    \hat s^* = \text{argmin}\Bigg\{f(x^*(s_j)):\hat H_j\geq \gamma + \frac{ q_{1-\beta}}{\sqrt{n_2}}, 1\leq j\leq p\Bigg\}.
\end{equation}
\STATE {\textbf{Output:} $\hat s^*,x^*(\hat s^*)$.}
\end{algorithmic}
\end{algorithm}



\begin{algorithm}
\caption{$V$: Normalized Gaussian Supremum Validator}
\label{algo:supremum_normalized}
\begin{algorithmic}
\STATE {\textbf{Input:} $\{\xi_1,\ldots,\xi_{n_2}\},\{x^*(s_1),\ldots,x^*(s_p)\},1-\beta$}
\vspace{1ex}

\STATE \textbf{1.}~Same as in Algorithm \ref{algo:supremum_unnormalized}.
\STATE \textbf{2.}~Denote $\hat\sigma_j^2=\hat \Sigma(j,j)$. Compute $q_{1-\beta}$, the $(1-\beta)$-quantile of $\max\{Z_j/\hat\sigma_j:\hat\sigma_j^2>0,1\leq j\leq p\}$ where $(Z_1,\ldots,Z_p)\sim N_p(0,\hat \Sigma)$, and let
\begin{equation}\label{validate:gs normalized}
    \hat s^* = \text{argmin}\Bigg\{f(x^*(s_j)):\hat H_j\geq \gamma + \frac{ q_{1-\beta}\hat\sigma_j}{\sqrt{n_2}}, 1\leq j\leq p\Bigg\}.
\end{equation}
\STATE {\textbf{Output:} $\hat s^*,x^*(\hat s^*)$.}
\end{algorithmic}
\end{algorithm}

The first Gaussian supremum validator (Algorithm \ref{algo:supremum_unnormalized}) is reasoned from a joint central limit theorem (CLT) that governs the convergence of $\sqrt{n_2}(\hat H_{1}-H(x^*(s_1)),\ldots,\hat H_{p}-H(x^*(s_p)))$ to $N_p(0,\Sigma)$, where $\Sigma(j_1,j_2)=\mathrm{Cov}_F(h(x^*(s_{j_1}),\xi),h(x^*(s_{j_2}),\xi))$. Using the sample covariance $\hat{\Sigma}$ from Step 1 of Algorithm \ref{algo:supremum_unnormalized} as an approximation of $\Sigma$, we have, by the continuous mapping theorem,
\begin{equation*}
    \max_{1\leq j\leq p}\sqrt{n_2}(\hat H_{j}-H(x^*(s_j)))\approx \max_{1\leq j\leq p}Z_j\text{\ \ in distribution}
\end{equation*}
where $(Z_1,\ldots,Z_p)\sim N_p(0,\hat\Sigma)$. Therefore using the $1-\beta$ quantile $q_{1-\beta}$ of the Gaussian supremum in the margin leads to
\begin{equation*}
    H(x^*(s_j))\geq \hat H_{j}-\frac{q_{1-\beta}}{\sqrt{n_2}}\text{ for all }j=1,\ldots,p,\text{ with probability}\approx 1-\beta.
\end{equation*}
The second validator (Algorithm \ref{algo:supremum_normalized}) uses an alternate version of the CLT that is normalized by the componentwise standard deviation $\sigma_j$, i.e., $\sqrt{n_2}((\hat H_{1}-H(x^*(s_1)))/\sigma_1,\ldots,(\hat H_{p}-H(x^*(s_p)))/\sigma_p)$ converges to $N_p(0,D\Sigma D)$, where $D$ is a diagonal matrix of $1/\sigma_j,j=1,\ldots,p$. Note that the quantile $q_{1-\beta}$ in both validators can be computed to high accuracy via Monte Carlo.

Let us make the above reasoning precise. We present our results for two cases that need separate treatments: When $H(x)\geq\gamma$  is a ``light-tailed" stochastic constraint, and when it is a chance constraint.

\subsection{Performance Guarantees for General Stochastic Constraints}
Recall that $H(x)=\mathbb E_F[h(x,\xi)]$. Denote $$\sigma^2(x):=\mathrm{Var}_F(h(x,\xi))$$
as the variance of $h$ for each decision $x\in\mathcal X$. We assume the following on optimization problem \eqref{stoc_opt}:
\begin{assumption}[Light-tailedness]\label{sub-gaussian bound}
There exists a constant $D_1\geq 1$ such that for all $x\in\mathcal X$ with $\sigma^2(x)>0$, we have
\begin{equation*}
\mathbb E_F\Big[\exp\Big(\frac{\lvert h(x,\xi)-H(x)\rvert^2}{D_1^2\sigma^2(x)}\Big)\Big]\leq 2\text{\ \ and\ \ }\mathbb E_F\Big[\Big(\frac{\lvert h(x,\xi)-H(x)\rvert}{\sigma(x)}\Big)^{2+k}\Big]\leq D_1^k\text{ for }k=1,2.
\end{equation*}
\end{assumption}
This assumption stipulates that the distribution of $h(x,\xi)$ after being centered and normalized by its standard deviation is sufficiently light-tailed at each $x$. Note that no other regularity property, e.g., convexity or continuity, is assumed for the function $h$ itself. We have the following finite-sample feasibility guarantees for the solution output by Algorithm \ref{algo:supremum_unnormalized} or \ref{algo:supremum_normalized}:
\begin{theorem}[Finite-sample feasibility guarantee for unnormalized validator]\label{feasibility_supremum_unnormalized:general constraint}
Suppose Assumption \ref{sub-gaussian bound} holds. Let $\overline{H}=\max_{1\leq j\leq p}H(x^*(s_j))$ and $\bar{\sigma}^2=\max_{1\leq j\leq p}\sigma^2(x^*(s_j))$. For every solution set $\{x^*(s_j):1\leq j\leq p\}$, every $n_2$, and $\beta\in(0,\frac{1}{2})$, the solution output by Algorithm \ref{algo:supremum_unnormalized} satisfies
\begin{align*}
P_{\bm{\xi}_{1:n_2}}(x^*(\hat s^*)\text{ is feasible for \eqref{stoc_opt}})\geq 1-\beta-C\prth{\Big(\frac{D_1^2\log^7(pn_2)}{n_2}\Big)^{\frac{1}{6}}+\exp\big(-\frac{cn_2\epsilon^2}{D_1^2\bar{\sigma}^2}\big)+p\exp\big(-\frac{cn_2}{D_1^4}\big)}
\end{align*}
with
\begin{equation}\label{validation accuracy:supremum unnormalized}
    \epsilon = \prth{\overline{H}-\gamma-C\bar{\sigma}\sqrt{\frac{\log(p/\beta)}{n_2}}}_+
\end{equation}
where $C$ and $c$ are universal constants, and $P_{\bm{\xi}_{1:n_2}}$ denotes the probability with respect to Phase two data $\{\xi_1,\ldots,\xi_{n_2}\}$ and conditional on Phase one data $\{\xi_{n_2+1},\ldots,\xi_n\}$.
\end{theorem}
\begin{theorem}[Finite-sample feasibility guarantee for normalized validator]\label{feasibility_supremum_normalized:general constraint}
Let $\overline{s}\in\argmax \{H(x^*(s_j)):j=1,\ldots,p\}$, i.e., $H(x^*(\overline{s}))=\overline{H}$. Under the same conditions of Theorem \ref{feasibility_supremum_unnormalized:general constraint}, the solution output by Algorithm \ref{algo:supremum_normalized} satisfies
\begin{eqnarray*}
&&P_{\bm{\xi}_{1:n_2}}(x^*(\hat s^*)\text{ is feasible for \eqref{stoc_opt}})\\
&\geq& 1-\beta-C\prth{\Big(\frac{D_1^2\log^7(pn_2)}{n_2}\Big)^{\frac{1}{6}}+\frac{D_1^2\log^2(pn_2)}{\sqrt n_2}+\exp\big(-\frac{cn_2\epsilon^2}{D_1^2\sigma^2(x^*(\overline{s}))}\big)+p\exp\big(-\frac{cn_2^{2/3}}{D_1^{10/3}}\big)}
\end{eqnarray*}
with
\begin{equation}\label{validation accuracy:supremum normalized}
    \epsilon = \prth{\overline{H}-\gamma-C\sigma(x^*(\overline{s}))\sqrt{\frac{\log(p/\beta)}{n_2}}}_+
\end{equation}
where $C$ and $c$ are universal constants.
\end{theorem}
In both Theorems \ref{feasibility_supremum_unnormalized:general constraint} and \ref{feasibility_supremum_normalized:general constraint}, the finite-sample coverage probability consists of two sources of errors. The first source comes from the CLT approximation that decays polynomially in the Phase 2 sample size $n_2$. The second error arises from the possibility that none of the solutions $\{x^*(s_1),\ldots,x^*(s_p)\}$ satisfies the criterion in \eqref{validate:gs unnormalized} or \eqref{validate:gs normalized}, which vanishes exponentially fast. When $\epsilon$ in \eqref{validation accuracy:supremum unnormalized} or \eqref{validation accuracy:supremum normalized} is of constant order, the CLT error dominates. In this case the finite-sample error depends logarithmically on $p$, the number of candidate parameter values, and the bounds dictate a coverage tending to $1-\beta$ when $p$ is as large as $\exp(o(n_2^{1/7}))$.

The derivation of the logarithmic dependence on $p$ in Theorem \ref{feasibility_supremum_unnormalized:general constraint} builds on a high-dimensional CLT and an associated multiplier bootstrap approximation recently developed in \cite{chernozhukov2017central} (Appendix \ref{clt:adaptation}). The proof of Theorem \ref{feasibility_supremum_normalized:general constraint} further requires a Hoeffding-type inequality for U-statistics to control the errors of the sample variance estimates, as well as the so-called Nazarov's inequality, an anti-concentration inequality for multivariate Gaussian, to control the coverage errors when using estimated standard deviations in the margin (Appendix \ref{clt:extension to normalized case}). Appendices \ref{multiplier bootstrap:unnormalized and normalized} and \ref{main proof: gaussian supremum validators} detail the proofs of Theorems \ref{feasibility_supremum_unnormalized:general constraint} and \ref{feasibility_supremum_normalized:general constraint} that put together the above mathematical developments.

We explain the implication on the dimensionality of the problem. Note that to sufficiently cover the whole solution path, $p$ is typically exponential in the dimension of $S$, denoted $\text{dim}(S)$ (this happens when we uniformly discretize the parameter space $S$). The discussion above thus implies a requirement that $n_2$ is of higher order than $\text{dim}(S)^7$. Here the low dimensionality of $S$ is crucial; for instance, a one-dimensional conservativeness parameter $s$ would mean $\text{dim}(S)=1$, so that a reasonably small $n_2$ can already ensure adequate feasibility coverage. Moreover, the margin adjustments in Algorithms \ref{algo:supremum_unnormalized} and \ref{algo:supremum_normalized} both depend only on $\text{dim}(S)$. Thus, choosing $\hat s^*$ relies only on $\text{dim}(S)$, but not the dimension of the whole decision space.
Note that Theorems \ref{feasibility_supremum_unnormalized:general constraint} and \ref{feasibility_supremum_normalized:general constraint} provide guarantee conditional on Phase one data. However, the universality of the involved constants means that analogous unconditional feasibility guarantees also hold if Assumption \ref{sub-gaussian bound} can be verified uniformly or with high probability with respect to Phase one data, an observation that persists for other subsequent results.

Comparing between the two validators, we also see that the normalized one (Algorithm \ref{algo:supremum_normalized}) is statistically more efficient than the unnormalized one (Algorithm \ref{algo:supremum_unnormalized}) when the variance $\sigma^2(x)$ exhibits high variability across solutions. More specifically, in order to make the exponential error non-dominant, one needs at least $\epsilon>0$. In the case of Algorithm \ref{algo:supremum_unnormalized}, expression \eqref{validation accuracy:supremum unnormalized} suggests that, after ignoring the logarithmic factor $\log(p/\beta)$, this requires an $n_2$ to be of order $\overline{\sigma}^2/(\overline{H}-\gamma)^2$. In contrast, for Algorithm \ref{algo:supremum_normalized} this becomes $\sigma^2(x^*(\overline{s}))/(\overline{H}-\gamma)^2$, where the maximal variance is replaced with the variance at the solution that optimizes the $H$-value, which in general does not have the maximal variance.


Theorems \ref{feasibility:gs unnormalized} and \ref{feasibility:gs normalized} also give immediately the following asymptotic feasibility guarantee (proof in Appendix \ref{main proof: gaussian supremum validators}):
\begin{corollary}[Asymptotic feasibility guarantee]\label{asymptotic feasibility:supremum + general constraint}
Suppose Assumption \ref{sub-gaussian bound} holds. Let $\overline{H}=\max_{1\leq j\leq p}H(x^*(s_j))$. For every solution set $\{x^*(s_j):1\leq j\leq p\}$ such that $\overline{H}>\gamma$ and every $\beta\in(0,\frac{1}{2})$, the solution output by Algorithm \ref{algo:supremum_unnormalized} or \ref{algo:supremum_normalized} satisfies
\begin{equation*}
\liminf_{n_2\to\infty\text{ and }p\exp(-n_2^{1/7})\to 0}P_{\bm{\xi}_{1:n_2}}(x^*(\hat s^*)\text{ is feasible for \eqref{stoc_opt}})\geq 1-\beta.
\end{equation*}
\end{corollary}

\subsection{Performance Guarantees for Chance Constraints}
Underlying the finite-sample bounds in Theorems \ref{feasibility_supremum_unnormalized:general constraint} and \ref{feasibility_supremum_normalized:general constraint} is the light-tailedness condition in  Assumption \ref{sub-gaussian bound}. However, in a CCP that takes the form
\begin{equation}\label{chance constraint}
\min_{x\in\mathcal X}\ f(x)\text{\ \ subject to\ \ }P(x):=\mathbb P_{F}((x,\xi)\in A)\geq 1-\alpha
\end{equation}
where $A\subseteq \R^d\times \R^m$ is a deterministic set and $1-\alpha$ is a tolerance level for the satisfaction probability, the tail of the normalized indicator function $\mathbf{1}((x,\xi)\in A)$ can be arbitrarily heavy as the satisfaction probability approaches $0$ or $1$ and hence violates Assumption \ref{sub-gaussian bound}. Thus, instead, we present different finite-sample error bounds for \eqref{chance constraint} than Theorems \ref{feasibility_supremum_unnormalized:general constraint} and \ref{feasibility_supremum_normalized:general constraint} whose derivations rely on the Bernoulli nature of the underlying function:
\begin{theorem}[Finite-sample chance constraint feasibility guarantee for unnormalized validator]\label{feasibility:gs unnormalized}
Let $\bar{\alpha}=1-\max_{1\leq j\leq p}P(x^*(s_j))$. For every solution set $\{x^*(s_j):1\leq j\leq p\}$, every $n_2$, and $\beta\in(0,\frac{1}{2})$, the solution output by Algorithm \ref{algo:supremum_unnormalized} satisfies
\begin{equation*}
P_{\bm{\xi}_{1:n_2}}(x^*(\hat s^*)\text{ is feasible for \eqref{chance constraint}})\geq 1-\beta-C\prth{\Big(\frac{\log^7(pn_2)}{n_2\alpha}\Big)^{\frac{1}{6}}+\exp\big(-cn_2\min\{\epsilon,\frac{\epsilon^2}{\bar{\alpha}}\}\big)}
\end{equation*}
with
\begin{equation}\label{statistical error unnormalized}
    \epsilon=\prth{\alpha-\bar{\alpha}-C\sqrt{\frac{\log( p/\beta)}{n_2}}}_+
\end{equation}
where $C$ and $c$ are universal constants.
\end{theorem}

\begin{theorem}[Finite-sample chance constraint feasibility guarantee for normalized validator]\label{feasibility:gs normalized}
Under the same conditions of Theorem \ref{feasibility:gs unnormalized}, the solution output by Algorithm \ref{algo:supremum_normalized} satisfies
\begin{equation*}
P_{\bm{\xi}_{1:n_2}}(x^*(\hat s^*)\text{ is feasible for \eqref{chance constraint}})\geq 1-\beta-C\prth{\Big(\frac{\log^7(pn_2)}{n_2\alpha}\Big)^{\frac{1}{6}}+\frac{\log^2 (pn_2)}{\sqrt{n_2\alpha}}+\exp\big(-cn_2\min\{\epsilon,\frac{\epsilon^2}{\bar{\alpha}}\}\big)}
\end{equation*}
with
\begin{equation}\label{statistical error normalized}
    \epsilon=\prth{\alpha-\bar{\alpha}-C\sqrt{\frac{\big(\bar{\alpha}+\log (n_2\alpha)/n_2\big)\log (p/\beta)}{n_2}}}_+
\end{equation}
where $C$ and $c$ are universal constants.
\end{theorem}
A comparison between Theorems \ref{feasibility:gs unnormalized} and \ref{feasibility:gs normalized} again reveals the higher statistical efficiency of Algorithm \ref{algo:supremum_normalized} than Algorithm \ref{algo:supremum_unnormalized} which, in the CCP context, applies to the case when the satisfaction probability is large (i.e., the common case). Suppose that $1-\alpha$ approaches $1$. In order to make $\epsilon>0$ in \eqref{statistical error unnormalized}, we need a sample size $n_2$ of order $(\alpha-\bar{\alpha})^{-2}$ (after ignoring the logarithmic factor $\log(p/\beta)$), whereas in \eqref{statistical error normalized} it can be seen to need only an $n_2$ of order $\alpha(\alpha-\bar{\alpha})^{-2}$, a much smaller size when $1-\alpha$ is close to $1$.

Lastly, we have the following asymptotic feasibility guarantee in the case of CCP in parallel to Corollary \ref{asymptotic feasibility:supremum + general constraint}:
\begin{corollary}[Asymptotic chance constraint feasibility guarantee]\label{asymptotic feasibility:supremum + chance constraint}
Let $\bar{\alpha}=1-\max_{1\leq j\leq p}P(x^*(s_j))$. For every solution set $\{x^*(s_j):1\leq j\leq p\}$ such that $\bar{\alpha}<\alpha$ and every $\beta\in(0,\frac{1}{2})$, the solution output by Algorithm \ref{algo:supremum_unnormalized} or \ref{algo:supremum_normalized} satisfies
\begin{equation*}
\liminf_{n_2\to\infty\text{ and }p\exp(-n_2^{1/7})\to 0}P_{\bm{\xi}_{1:n_2}}(x^*(\hat s^*)\text{ is feasible for \eqref{chance constraint}})\geq 1-\beta.
\end{equation*}
\end{corollary}

Appendix \ref{main proof: gaussian supremum validators} details the proofs of Theorem \ref{feasibility:gs unnormalized}, Theorem \ref{feasibility:gs normalized} and Corollary \ref{asymptotic feasibility:supremum + chance constraint}.

To close this section, we note that our Gaussian supremum validators also enjoy a notion of asymptotic solution-path optimality under additional assumptions. To streamline our presentation, we defer this discussion to the next section and combine it with the discussion of our next validator.

\section{Validation via Univariate Gaussian Margin}\label{sec:standard gaussian margin}

We offer an alternate validator that can perform more efficiently than Algorithms \ref{algo:supremum_unnormalized} and \ref{algo:supremum_normalized}, provided that further regularity assumptions are in place. This is a scheme that simply uses a standard univariate Gaussian critical value to calibrate the margin (Algorithm \ref{algo:marginal}).

Algorithm \ref{algo:marginal} outputs a solution with objective value no worse than Algorithms \ref{algo:supremum_unnormalized} and \ref{algo:supremum_normalized}. Comparing the criteria to choose $\hat s^*$, we see that, due to the stochastic dominance of the maximum among a multivariate Gaussian vector over each of its individual components, the margin in \eqref{validate:gs unnormalized} satisfies $q_{1-\beta}\geq z_{1-\beta}\hat\sigma_j$ for all $j$, and similarly the margin in \eqref{validate:gs normalized} satisfies $q_{1-\beta}\hat\sigma_j\geq z_{1-\beta}\hat\sigma_j$, so that both are bounded from below by the margin in \eqref{validate:marginal gaussian}. Consequently the solution from \eqref{validate:marginal gaussian} achieves an objective value no worse than the other two.



\begin{algorithm}
\caption{$V$: Univariate Gaussian Validator}
\label{algo:marginal}
\begin{algorithmic}
\STATE {\textbf{Input:} $\{\xi_{1},\ldots,\xi_{n_2}\},\{x^*(s_1),\ldots,x^*(s_p)\},1-\beta$}
\vspace{1ex}
\STATE \textbf{1.}~For each $j=1,\ldots,p$ compute the sample mean $\hat H_j=(1/n_2)\sum_{i=1}^{n_2}h(x^*(s_j),\xi_i)$ and sample variance $\hat \sigma_j^2=(1/n_2)\sum_{i=1}^{n_2}(h(x^*(s_j),\xi_i)-\hat H_j)^2$.
\STATE \textbf{2.}~Compute
\begin{equation}\label{validate:marginal gaussian}
   \hat s^* = \text{argmin}\Bigg\{f(x^*(s_j))\Big\vert \hat H_j\geq \gamma + \frac{z_{1-\beta}\hat \sigma_j}{\sqrt{n_2}},1\leq j\leq p\Bigg\}
\end{equation}
where $z_{1-\beta}$ is the $1-\beta$ quantile of the standard Gaussian distribution.
\vspace{1ex}
\STATE {\textbf{Output:} $\hat s^*,x^*(\hat s^*)$.}
\end{algorithmic}
\end{algorithm}


 The univariate Gaussian critical value used in the margin in Algorithm \ref{algo:marginal} hints that feasibility needs to be validated at only one value of $s$ instead of the solution path $S$. The validity of this procedure is based on the statistical consistency of the obtained solution $x^*(\hat s^*)$ to some limiting solution (correspondingly $\hat s^*$ to some limiting optimal parameter value) as $n_2$ increases. Intuitively, this implies that with sufficient sample size one can focus feasibility validation on a small neighborhood of $\hat s^*$, which further suggests that we need to control only the statistical error at effectively one solution parametrized at $\hat s^*$. For this argument to hold, however, we would need several additional technical assumptions including a low functional complexity of $h$, and a different line of derivations.

\subsection{Asymptotic Performance Guarantees}\label{sec:asymptotic guarantee univariate}
We present the statistical guarantees of Algorithm \ref{algo:marginal} as Phase two data size $n_2\to\infty$. We assume continuity for the objective of \eqref{stoc_opt}:
\begin{assumption}[Continuous objective]\label{continuous obj}
The objective function $f(x)$ is continuous on $\mathcal X$.
\end{assumption}
For the constraint, we assume the following:
\begin{assumption}[Functional complexity]\label{Donsker h}
The function class $\mathcal F:=\{h(x,\cdot)\vert x\in\mathcal X\}$ is $F$-Donsker.
\end{assumption}
\begin{assumption}[$L_2$-boundedness]\label{moment2}
$\mathbb E_F\big[\sup_{x\in\mathcal X}\abs{h(x,\xi)-H(x)}^2\big]<\infty$.
\end{assumption}
\begin{assumption}[$L_2$-continuity]\label{continuous constraint}
For any fixed $x\in\mathcal X$ and another $x'\in\mathcal X$, we have $\lim_{x'\to x}\mathbb E_F[(h(x',\xi)-h(x,\xi))^2]=0$.
\end{assumption}

To give a sense of the generality of the above assumptions, we identify two general classes of constraints for which these assumptions are guaranteed to hold, one suitable for general $h$, and another one for CCPs:
\begin{proposition}\label{two valid classes of constraints}
Assumptions \ref{Donsker h}-\ref{continuous constraint} hold in each of the following two cases:
\begin{enumerate}
\item[i.]There exists some $M(\xi)$ such that $\mathbb E_F[M(\xi)^2]<\infty$ and $\abs{h(x_1,\xi)-h(x_2,\xi)}\leq M(\xi)\norm{x_1-x_2}$ for all $x_1,x_2\in\mathcal X$, there exists some $\tilde x\in\mathcal X$ such that $\mathbb E_F[h(\tilde x,\xi)^2]<\infty$, and the decision space $\mathcal X$ is compact;
\item[ii.]$h(x,\xi)=\mathbf{1}(a'_kA_k(x)\leq b_k\text{ for }k=1,\ldots,K)$ for some $K<\infty$, where each $A_k(\cdot):\R^d\to \R^{m_k}$ is a continuous mapping and each $a_k\in \R^{m_k},b_k\in \R$ satisfies either (i) $a_k$ has a density on $\R^{m_k}$ and $b_k$ is a non-zero constant or (ii) $(a_k,b_k)$ has a density on $\R^{m_k+1}$.
\end{enumerate}
\end{proposition}
Case (i) in Proposition \ref{two valid classes of constraints} follows from standard results in empirical process theory,  including in particular the Jain-Marcus Theorem. The proof of Case (ii) involves checking the finite Vapnik-Chervonenkis (VC) dimension and pointwise separability of the function class in order to verify $F$-Donskerness. Appendix \ref{sec:proofs univariate} details the proof of Proposition \ref{two valid classes of constraints}.

We impose one more assumption on the constraint function regarding its variance:
\begin{assumption}[Non-degeneracy of the variance on the boundary]\label{non-degenerate variance}
$\sigma^2(x)>0$ for all $x\in\mathcal X$ such that $H(x)=\gamma$.
\end{assumption}
In Assumption \ref{non-degenerate variance}, non-zero variance is assumed only for those $x$'s at which the stochastic constraint is satisfied with equality, but not necessarily for other $x$. This is significant in the case of CCP \eqref{chance constraint}. While there could exist $100\%$ or $0\%$ safe solutions, i.e., $x$ such that $P(x)=1$ or $0$, and hence non-degeneracy may not be satisfied over the whole $\mathcal X$, it holds for those $x$'s with $P(x)=1-\alpha$ that have (the same) non-zero variance $\alpha(1-\alpha)$.


Now we present our assumptions on the data-driven reformulation $OPT(s),s\in S$. We focus on formulations with a single parameter (A separate set of results for formulations with multiple parameters can be found in Appendix \ref{sec:multidimensional parameter univariate gaussian}). We first assume that the solution path is piecewise continuous:
\begin{assumption}[Piecewise continuous solution path]\label{piecewise uniform continuous solution curve}
The parameter space $S$ is a finite interval $[s_l,s_u]$. The optimal solution $x^*(s)$ of $OPT(s)$ exists and is unique except for a finite number of parameter values $\tilde s_i, i=1,\ldots,M-1$ such that $s_l=\tilde s_0< \tilde s_1<\cdots<\tilde s_{M-1}<\tilde s_{M}=s_u$, and the parameter-to-solution mapping $x^*(s)$ is uniformly continuous on each piece $[\tilde s_0,\tilde s_1)$, $(\tilde s_{M-1},\tilde s_M]$, and $(\tilde s_{i-1},\tilde s_i)$ for $i=2,\ldots,M-1$.
\end{assumption}
Continuity of the solution path allows approximating the whole solution curve by discretizing the parameter space $S$. Also note that under Assumption \ref{piecewise uniform continuous solution curve} the solution $x^*(s)$ exists and is unique for almost surely every $s\in S$ with respect to the Lebesgue measure. Therefore, if one discretizes the parameter space by randomizing via a continuous distribution over $S$, then with probability one the solution $x^*(s)$ is unique at all sampled parameter values. This provides an easy way to ensure the assumption that none of the parameter values $\{s_1,\ldots,s_p\}$ used in Phase one of Algorithm \ref{two-phase framework} belongs to the discontinuity set $\{\tilde s_1,\ldots,\tilde s_{M-1}\}$.


To explain the superior performance of Algorithm \ref{algo:marginal}, we introduce a notion of optimality within the solution path $\{x^*(s):s\in S\}$. First, since the parameter-to-solution mapping $x^*(s)$ is not defined at the discontinuities under Assumption \ref{piecewise uniform continuous solution curve}, we need to fill in these holes in the solution path. Thanks to uniform continuity, the mapping $x^*(s)$ on each piece $(\tilde s_{i-1},\tilde s_i)$ can be continuously extended to the closure $[\tilde s_{i-1},\tilde s_i]$ by taking left and right limits. Specifically, we define:
\begin{definition}
Under Assumption \ref{piecewise uniform continuous solution curve}, the parameter-to-solution mapping $x^*(\cdot)$ at each discontinuity $\tilde s_i,i=1,\ldots,M-1$ is defined in an extended fashion as
\begin{equation*}
x^*(\tilde s_i)=\{x^*(\tilde s_i-),x^*(\tilde s_i+)\}\text{ where }x^*(\tilde s_i-):=\lim_{s\to\tilde s_i-}x^*(s)\text{ and }x^*(\tilde s_i+):=\lim_{s\to \tilde s_i+}x^*(s).
\end{equation*}
\end{definition}
Note that the two solutions $x^*(\tilde s_i-)$ and $x^*(\tilde s_i+)$ are different if the $i$-th and $(i+1)$-th pieces are disconnected. With the extended parameter-to-solution mapping $x^*(\cdot)$, we now introduce the notions of optimal solution and optimal parameter associated with the solution path:
\begin{definition}
Associated with the solution path $\{x^*(s):s\in S\}$, the optimal solution set is
\begin{equation}\label{optimal solution}
    \mathcal X_S^*:=\text{argmin}\{f(x):H(x)\geq \gamma,x=x^*(s)\text{ for }s\notin\{\tilde s_1,\ldots,\tilde s_{M-1}\}\text{ or }x\in x^*(\tilde s_i)\text{ for some }i=1,\ldots,M-1\}
\end{equation}
and the optimal parameter set is
\begin{equation}\label{optimal parameter}
    S^*:=\{s\notin\{\tilde s_1,\ldots,\tilde s_{M-1}\}:x^*(s)\in \mathcal X_S^*\}\cup\{\tilde s_i:x^*(\tilde s_i) \cap \mathcal X_S^*\neq \emptyset,i=1,\ldots,M-1\}.
\end{equation}
\end{definition}
We need several additional technical assumptions. The first is that the stochastic constraint is not binding at the endpoints of each piece of the solution path:
\begin{assumption}\label{non-binding}
$H(x^*(\tilde s_i-))\neq \gamma$ and $H(x^*(\tilde s_i+))\neq \gamma$ for all $i=1,\ldots,M-1$, $H(x^*( s_l))\neq \gamma$, $H(x^*( s_u))\neq \gamma$, and $\sup_{s\notin \{\tilde s_1,\ldots,\tilde s_{M-1}\}}H(x^*(s))> \gamma$.
\end{assumption}
Since the solution path $\{x^*(s):s\in S\}$ depends on Phase one data $\bm{\xi}_{n_2+1:n}$, the path and hence the endpoints $x^*(\tilde s_i-),x^*(\tilde s_i+)$ are random objects, and so the first part of Assumption \ref{non-binding} is expected to hold almost surely provided that the set $\{x\in\mathcal X:H(x)=\gamma\}$ is a null set under the Lebesgue measure on $\R^d$. The second part states that the solution path contains a strictly feasible solution which in turn ensures that the optimal solution set $\mathcal X_S^*$ is non-empty. Note that this can typically be achieved by simply including very conservative parameter values in $S$.

Another property we assume regards the monotonicity of the feasible set size with respect to the parameter $s$ in the reformulation $OPT(s)$:
\begin{assumption}\label{monotonicity}
Denote by $\mathrm{Sol}(s):=\mathcal X\cap \hat{\mathcal F}(s)$ the feasible set of $OPT(s)$. Assume $\mathrm{Sol}(s)$ is a closed set for all $s\in S$ and $\mathrm{Sol}(s_2)\subseteq \mathrm{Sol}(s_1)$ for all $s_1,s_2\in S$ such that $s_1<s_2$.
\end{assumption}
Assmption \ref{monotonicity} holds for all common reformulations (all examples in the beginning of Section \ref{sec:overview}) as $s$ controls the conservativeness level. For instance, in RO with ellipsoidal uncertainty set, the RO feasible region shrinks with the radius of the ellipsoid, and similar relations hold for DRO, SAA, and SO. A straightforward consequence of Assumption \ref{monotonicity} is the monotonicity of the parameter-to-objective mapping
\begin{equation*}
v(s):=\inf\{f(x):x\in \mathcal X\cup \hat{\mathcal F}(s)\}
\end{equation*}
as described in the following proposition (proof in Appendix \ref{sec:proofs univariate}):
\begin{proposition}\label{monotonic objective value}
Suppose Assumptions \ref{piecewise uniform continuous solution curve} and \ref{monotonicity} hold. For all $s_1,s_2\in S$ such that $s_1<s_2$ it holds $v(s_1)\leq v(s_2)$, and if additionally $s_1,s_2\notin \{\tilde s_1,\ldots,\tilde s_{M-1}\}$ then $v(s_1)< v(s_2)$ if and only if $x^*(s_1)\neq x^*(s_2)$.
\end{proposition}

The assumptions we have made for the formulation $OPT(s)$ give rise to the following uniqueness characterization of the optimal solution set $\mathcal X_S^*$ and the optimal parameter set $S^*$ within the solution path, which would be used to establish the feasibility guarantees for Algorithm \ref{algo:marginal}.
\begin{proposition}[Structure of solution-path optima]\label{unique optimum}
Under Assumptions \ref{continuous obj}, \ref{continuous constraint}, and \ref{piecewise uniform continuous solution curve}-\ref{monotonicity}, the optimal solution set $\mathcal X_S^*$ is a singleton $\{x_S^*\}$ and the optimal parameter set $S^*$ is a closed interval $[s_l^*,s_u^*]$ for $s_l^*,s_u^*\in S$. In addition, if $v(s)$ is strictly increasing on $S$, then $S^*$ is a singleton $\{s^*\}$.
\end{proposition}
The proof of Proposition \ref{unique optimum}, which is in Appendix \ref{sec:proofs univariate}, involves an exhaustion of all possible structures of the set $\mathcal X_S^*$ that contain more than one solution, and showing each of them contradict with our assumptions (especially Assumption \ref{monotonicity}).


Lastly, we assume the following technical assumption for the set of optima:
\begin{assumption}\label{unique optimal solution}
For any $\epsilon>0$ there exists an $s\notin\{\tilde s_1,\ldots,\tilde s_{M-1}\}$ such that $H(x^*(s))>\gamma$ and $\Vert x^*(s)-x_S^*\Vert_2<\epsilon$, where $x_S^*$ is the unique optimal solution from Proposition \ref{unique optimum}.
\end{assumption}
This assumption trivially holds if $\mathcal X_S^*=\{x_S^*\}$ as described in Proposition \ref{unique optimum} and $H(x_S^*)>\gamma$. Otherwise, if $H(x_S^*)=\gamma$, it rules out the case that the solution path $x^*(s)$ passes through $x_S^*$ without entering the interior of the feasible set of \eqref{stoc_opt}. The latter exceptional case typically happens with zero probability, in view of the fact that the solution path is itself random with respect to Phase one data.



Now we are ready to present the asymptotic performance guarantee for Algorithm \ref{algo:marginal}:
\begin{theorem}[Asymptotic joint feasibility$+$optimality guarantee]\label{asymptotic joint:gaussian}
Suppose Assumptions \ref{continuous obj}-\ref{non-degenerate variance} hold for \eqref{stoc_opt}. Also suppose Assumptions \ref{piecewise uniform continuous solution curve}-\ref{unique optimal solution} hold for the reformulation $OPT(s)$ constructed in Algorithm \ref{two-phase framework}, and $\{s_1,\ldots,s_p\}\cap \{\tilde s_1,\ldots,\tilde s_{M-1}\}=\emptyset$. Denote by $\epsilon_S = \sup_{s\in S}\inf_{1\leq j\leq p}\abs{s-s_j}$ the mesh size, and by $x_S^*$ the unique optimal solution from Proposition \ref{unique optimum}. Then, with respect to $\{\xi_{1},\ldots,\xi_{n_2}\}$, the solution and parameter output by Algorithm \ref{algo:marginal} satisfy
\begin{equation}\label{asymptotic optimality}
    \lim_{n_2\to\infty,\epsilon_S\to 0}x^*(\hat{s}^*)=x_S^* \text{ and } \lim_{n_2\to\infty,\epsilon_S\to 0}d(\hat s^*,S^*)=0
\end{equation}
almost surely. Moreover, if $H(x_S^*)=\gamma$ we have
\begin{equation}\label{feasibility gaussian:boundary}
    \liminf_{n_2\to\infty,\epsilon_S\to 0}P_{\bm{\xi}_{1:n_2}}(x^*(\hat{s}^*)\text{ is feasible for \eqref{stoc_opt}})\geq 1-\beta,
\end{equation}
otherwise if $H(x_S^*)>\gamma$ we have
\begin{equation}\label{feasibility gaussian:interior}
\lim_{n_2\to\infty,\epsilon_S\to 0}P_{\bm{\xi}_{1:n_2}}(x^*(\hat{s}^*)\text{ is feasible for \eqref{stoc_opt}})= 1.
\end{equation}
\end{theorem}
Theorem \ref{asymptotic joint:gaussian} states that as the mesh $\{s_1,\ldots,s_p\}$ gets increasingly fine and the data size grows, the solution given by Algorithm \ref{algo:marginal} enjoys performance guarantees concerning both feasibility and solution-path optimality. In particular, the estimated solution and the conservativeness parameter converge to the optimal solution $x_S^*$ and the optimal parameter set $S^*$ respectively, while simultaneously the obtained solution is feasible with the desired confidence level $1-\beta$.

The proof of Theorem \ref{asymptotic joint:gaussian} is in Appendix \ref{sec:proofs univariate}. The consistency result in \eqref{asymptotic optimality} is shown via a dense approximation of the discrete parameter set $\{s_1,\ldots,s_p\}$ on the continuum $S$, through the continuity of the solution path and a uniform law of large numbers. Then, based on this consistency, the feasibility guarantee \eqref{feasibility gaussian:boundary} is established by showing $P\big(H(x^*(\hat s^*))\geq \gamma\big)\geq P\big(H(x^*(\hat s^*))\geq \hat H(x^*(\hat s^*))-z_{1-\beta}\hat\sigma(x^*(\hat s^*))/\sqrt{n_2}\big)\approx P\big(H(x^*_S)\geq \hat H(x^*_S)-z_{1-\beta}\hat\sigma(x^*_S)/\sqrt{n_2}\big)\to 1-\beta$, where the ``$\geq$'' follows from our validation criterion \eqref{validate:marginal gaussian} whereas the ``$\approx$'' comes from the asymptotic tightness of the empirical process $\{\sqrt{n_2}(\hat H(x^*(s))-H(x^*(s))):s\in S\}$ and the $L_2$ continuity of the constraint function $h(x,\xi)$.

Furthermore, under additional smoothness conditions on the constraint function $h$ and the solution path $\{x^*(s):s\in S\}$, we also establish the finite-sample counterparts for the optimality guarantee \eqref{asymptotic optimality} and feasibility guarantee \eqref{feasibility gaussian:boundary} for Algorithm \ref{algo:marginal}. These are presented in Appendix \ref{sec:finite sample error univariate gaussian}.

Note that the confidence level \eqref{feasibility gaussian:boundary} at which Algorithm \ref{algo:marginal} outputs a feasible solution (and also Algorithms \ref{algo:supremum_unnormalized} and \ref{algo:supremum_normalized}, i.e., Corollaries \ref{asymptotic feasibility:supremum + general constraint} and \ref{asymptotic feasibility:supremum + chance constraint}) is generally not tight, i.e., a lower bound instead of an equality is guaranteed. However, with a strict monotonicity condition on the reformulation $OPT(s)$ and a finer discretization mesh for the conservativeness parameter, Algorithm \ref{algo:marginal} can give a tight confidence guarantee:


\begin{theorem}[Asymptotically tight feasibility guarantee]\label{asymptotic tight coverage:gaussian}
In addition to the conditions of Theorem \ref{asymptotic joint:gaussian}, further assume that the parameter-to-objective mapping $v(s)$ is strictly increasing on $S$, and consider the case that $H(x_S^*)=\gamma$. If the mesh $\{s_1,\ldots,s_p\}$ is fine enough so that \begin{equation}\label{high resolution}
    \max_{i=1,\ldots,M}\max_{j=1,\ldots,p_i-1}\abs{H(x^*(s^i_j))-H(x^*(s^i_{j+1}))}=o\big(\frac{1}{\sqrt{n_2}}\big)
\end{equation}
where $s^i_1<\cdots<s^i_{p_i}$ are the parameter values $\{s_j:s_j\in(\tilde s_{i-1},\tilde s_i),j=1,\ldots,p\}$ (so that $\sum_{i=1}^Mp_i=p$), then we must have $$\lim_{n_2\to\infty\text{ and }\epsilon_S\to 0\text{ s.t. \eqref{high resolution} holds}}P_{\bm{\xi}_{1:n_2}}(x^*(\hat{s}^*)\text{ is feasible for \eqref{stoc_opt}})= 1-\beta$$
for the solution output by Algorithm \ref{algo:marginal}.
\end{theorem}
Roughly speaking, the loose confidence guarantee in \eqref{feasibility gaussian:boundary} can be attributed to the one-sided nature of the inequality criterion used in \eqref{validate:marginal gaussian}. The monotonicity of $v(s)$ and the mesh condition \eqref{high resolution} give rise to a tight confidence guarantee by strengthening this inequality criterion to an equality (with a negligible error) at the chosen parameter value $\hat s^*$. Note that, when the expected constraint value $H(x^*(s))$ is Lipschitz continuous in the parameter, the mesh condition \eqref{high resolution} is guaranteed if $\epsilon_S=o\big(\frac{1}{\sqrt{n_2}}\big)$ or if $\frac{p}{\sqrt{n_2}}\to\infty$ and the mesh is equispaced. The proof of Theorem \ref{asymptotic tight coverage:gaussian} is in Appendix \ref{sec:proofs univariate}.

Relatedly, the following shows that, like Algorithm \ref{algo:marginal}, the supremum-based validators in Algorithms \ref{algo:supremum_unnormalized} and \ref{algo:supremum_normalized} also exhibit joint asymptotic feasibility and solution-path optimality guarantees. However, their confidence guarantees for feasibility are not as tight. This result complements our discussions at the end of Section \ref{sec:gaussian supremum margin} regarding the optimality property of the supremum-based validators, and also at beginning of Section \ref{sec:standard gaussian margin} regarding the better objective value of the solution obtained by Algorithm \ref{algo:marginal}, which is consistent with its tighter achievement of the feasibility confidence level.
\begin{theorem}[Asymptotic joint feasibility$+$optimality guarantee with Gaussian supremum validator]\label{asymptotic joint:gaussian supremum}
Under the same conditions as Theorem \ref{asymptotic joint:gaussian}, the solution and parameter from Algorithm \ref{algo:supremum_unnormalized} satisfy the consistency guarantee \eqref{asymptotic optimality}. In the case $H(x_S^*)=\gamma$ it holds
\begin{equation}\label{asymptotic feasibility gaussian supremum_unnormalized:boundary}
    \liminf_{n_2\to\infty,\epsilon_S\to 0}P_{\bm{\xi}_{1:n_2}}(x^*(\hat{s}^*)\text{ is feasible for \eqref{stoc_opt}})\geq \Phi\big(\frac{\bar{q}_{1-\beta}}{\sigma(x^*_S)}\big)\geq 1-\beta
\end{equation}
where $\bar{q}_{1-\beta}$ is the $1-\beta$ quantile of the supremum of the Gaussian process indexed by $s\in S\backslash\{\tilde s_1,\ldots,\tilde s_{M-1}\}$ with the covariance structure $\mathrm{Cov}(s,s')=\mathrm{Cov}_F(h(x^*(s),\xi),h(x^*(s'),\xi))$, and $\Phi$ is the distribution function of the standard normal.

If it is further assumed that $\inf_{x\in\mathcal X}\sigma^2(x)>0$, then \eqref{asymptotic optimality} also holds for Algorithm \ref{algo:supremum_normalized}, and in the case $H(x_S^*)=\gamma$ we have
\begin{equation}\label{asymptotic feasibility gaussian supremum_normalized:boundary}
    \liminf_{n_2\to\infty,\epsilon_S\to 0}P_{\bm{\xi}_{1:n_2}}(x^*(\hat{s}^*)\text{ is feasible for \eqref{stoc_opt}})\geq \Phi(\tilde{q}_{1-\beta})\geq 1-\beta
\end{equation}
where $\tilde{q}_{1-\beta}$ is the $1-\beta$ quantile of the supremum of the Gaussian process
on $S\backslash\{\tilde s_1,\ldots,\tilde s_{M-1}\}$ with covariance $\mathrm{Cov}(s,s')=\mathrm{Cov}_F(h(x^*(s),\xi),h(x^*(s'),\xi))/(\sigma(x^*(s))\sigma(x^*(s')))$.
\end{theorem}
In general, when the Gaussian processes involved in \eqref{asymptotic feasibility gaussian supremum_unnormalized:boundary} and \eqref{asymptotic feasibility gaussian supremum_normalized:boundary} have non-constant covariance structures, the asymptotic confidence levels rendered by Algorithms \ref{algo:supremum_unnormalized} and \ref{algo:supremum_normalized} are strictly higher than the nominal level $1-\beta$. This suggests that supremum-based margins tend to generate more conservative solutions than the univariate Gaussian margin does, although they all approach the same optimal solution $x_S^*$ in the limit.

The proof of Theorem \ref{asymptotic joint:gaussian supremum} (in Appendix \ref{sec:proofs univariate}) involves steps similar to that of Theorem \ref{asymptotic joint:gaussian}, but furthermore showing the statistical consistency of the critical value $q_{1-\beta}$ calibrated in Algorithm \ref{algo:supremum_unnormalized} or \ref{algo:supremum_normalized}. The latter utilizes the separability of the limiting Gaussian process and a control of errors coming from the associated multiplier bootstrap approximation.

\section{Applying Our Framework in Data-Driven Reformulations}\label{sec:formulation showcase}
In this section we showcase various data-driven reformulations of \eqref{stoc_opt} or \eqref{chance constraint} to which our proposed framework can be applied. We first comment that our Gaussian supremum validators (Algorithms \ref{algo:supremum_unnormalized} and \ref{algo:supremum_normalized}) are applicable to all formulations considered here, as long as the constraint function $h(x,\xi)$ is sufficiently light-tailed as described in Assumption \ref{sub-gaussian bound} or the constraint is a chance constraint. That is,
\begin{theorem}[Applicability of Gaussian supremum validators]
Consider the general stochastically constrained problem \eqref{stoc_opt} that satisfies Assumption \ref{sub-gaussian bound}, or CCP \eqref{chance constraint}. All the data-driven reformulations $OPT(s)$ presented below, namely SAA, DRO with $\phi$-divergence, Wasserstein and moment-based uncertainty sets, RO with polyhedral and ellipsoidal uncertainty sets, and SO (the last two approaches are for CCP only), can be validated by the Gaussian supremum validators in Algorithms \ref{algo:supremum_unnormalized} and \ref{algo:supremum_normalized} and elicit the conclusions in all theorems and corollaries in Section \ref{sec:gaussian supremum margin}.
\end{theorem}

The tighter univariate Gaussian validator (Algorithm \ref{algo:marginal}) however requires some extra regularity conditions from the data-driven formulation $OPT(s)$, but still works for many common formulations. We consider decision space $\mathcal X$ that has the form:
\begin{assumption}\label{convex decision space}
$\mathcal X=\{x\in \R^d:f_r(x)\leq 0\text{ for }r=1,\ldots,R\text{ and }Wx\leq z\}$, where each $f_r$ is continuous and convex, and $W=[w_1,w_2,\ldots,w_L]'\in \R^{L\times d},z\in\R^L$.
\end{assumption}
We consider optimization formulations that satisfy the following two assumptions:
\begin{assumption}[Slater's condition]\label{Slater's condition}
Slater's condition holds for $OPT(s_u):=\min\{f(x):x\in\mathcal X\cap\hat{\mathcal F}(s_u) \}$ where $s_u$ is the maximal parameter value.
\end{assumption}
\begin{assumption}[Non-empty and bounded level set]\label{bounded level set}
There exists a constant $c$ such that $\mathcal X\cap \hat{\mathcal F}(s_l)\cap \{x:f(x)\leq c\}$ is bounded and $\mathcal X\cap \hat{\mathcal F}(s_u)\cap \{x:f(x)\leq c\}$ is non-empty where $s_l,s_u$ are the minimal and maximal parameter values.
\end{assumption}
Slater's condition (Assumption \ref{Slater's condition}) is a common property that is expected to hold for most optimization problems in practice. \cite{dur2016genericity} states that Slater's condition is a generic property for linear conic programs by showing that it holds for all problem data except in a null set of Lebesgue measure. Assumption \ref{bounded level set} also trivially holds in many settings, e.g., when $\mathcal X\cap\hat{\mathcal F}(s_l)$ is compact or $f(x)$ is coercive. Under these two assumptions, stability results from parametric optimization (Proposition 4.4 in \cite{bonnans2013perturbation}) ensure that the solution path $x^*(s)$ is continuous when the optimal solution is unique for each $OPT(s)$, or piecewise continuous when uniqueness fails at only a finite number of parameter values, leading to Assumption \ref{piecewise uniform continuous solution curve}. Since other assumptions from Section \ref{sec:asymptotic guarantee univariate} regarding $OPT(s)$ can be readily verified to hold in general, for each considered formulation below we focus on identifying the conditions that guarantee the validity of Assumption \ref{piecewise uniform continuous solution curve} in order to ensure the asymptotic feasibility and optimality guarantees. The proofs of all results in this section are presented in Appendix \ref{sec:application proofs}.

We introduce a condition that will appear in the following discussion. Consider the linear objective $f(x)=c'x$ for some deterministic $c\in \R^d$. We say a finite collection of vectors $\{v_1,\ldots,v_k\}\subset \R^d$ with $k\leq d-1$ satisfies the strict cone inclusion (SCI) condition if
\begin{equation*}
\text{SCI: }v_1,\ldots,v_k\text{ are linearly independent, and there exist }\lambda_1,\ldots,
\lambda_k>0\text{ such that }\sum_{i=1}^k\lambda_iv_i=-c.
\end{equation*}

\smallskip

\noindent\underline{SAA:} First consider the SAA reformulation for the general stochastic constraint in \eqref{stoc_opt} in the form
\begin{equation}\label{general constraint:SAA}
\min_{x\in\mathcal X}\ f(x)\text{\ \ subject to\ \ }\frac{1}{n}\sum_{i=1}^nh(x,\xi_i)\geq \gamma +s
\end{equation}
where $s>0$ is the margin to be tuned (and for convenience, in this section only, we use $n$ to represent a generic sample size; in applications this typically refers to the Phase one data size). We have the following result concerning the applicability of Algorithm \ref{algo:marginal}:
\begin{theorem}[Applying univariate Gaussian validator to SAA]\label{continuity of SAA}
Consider $OPT(s)$ using \eqref{general constraint:SAA}. Suppose Assumptions \ref{convex decision space}-\ref{bounded level set} hold. In either of the following two cases:
\begin{enumerate}
\item[i.]$f(x)$ is continuous and strictly convex, $h(x,\xi)$ is continuous and concave in $x$ for every $\xi$.
\item[ii.]$f(x)=c'x$ for some non-zero $c\in\R^d$, the functions $f_r,r=1,\ldots,R$ in Assumption \ref{convex decision space} are strictly convex, any $k\leq d-1$ rows $\{w_{l(1)},\ldots,w_{l(k)}\}$ of $W$ do not satisfy the SCI condition, and $h(x,\xi)$ either is continuous and strictly concave in $x$ for every $\xi$ or has the form $h(x,\xi)=A(\xi)'x+b(\xi)$ where $A(\xi)\in\R^d$ has a density on $\R^d$.
\end{enumerate}
Assumption \ref{piecewise uniform continuous solution curve} holds with $M=1$ almost surely in the data $\{\xi_1,\ldots,\xi_n\}$.
\end{theorem}
The proof of Theorem \ref{continuity of SAA} (and theorems for other formulations below) mainly consists of establishing the joint continuity of the data-driven constraint \eqref{general constraint:SAA} in $x$ and $s$, and the uniqueness of $x^*(s)$, two main ingredients that enable us to apply the stability theory from \cite{bonnans2013perturbation} to conclude the continuity of $x^*(s)$. The former is shown by direct verification, whereas the latter is established from either strict convexity or the SCI condition when the formulation has linear objectives and constraints.

In the case of chance constraint \eqref{chance constraint}, the SAA formulation has the form
\begin{equation}\label{chance constraint:SAA}
\min_{x\in\mathcal X}\ f(x)\text{\ \ subject to\ \ }\frac{1}{n}\sum_{i=1}^n\mathbf{1}((x,\xi_i)\in A)\geq 1-\alpha +s.
\end{equation}
Note that the left hand side can only take values $\frac{j}{n},j=0,1,\ldots,n$, therefore all $s$ such that $1-\alpha +s\in(\frac{j-1}{n},\frac{j}{n}]$ lead to the same feasible region and hence the same solution $x^*(s)$. As a result, the solution path $\{x^*(s):s\in S\}$ consists of at most $n$ constant pieces and Assumption \ref{piecewise uniform continuous solution curve} holds automatically. Thus we have:

\begin{theorem}[Applying univariate Gaussian validator to SAA under chance constraint]
Consider $OPT(s)$ using \eqref{chance constraint:SAA}. Assumption \ref{piecewise uniform continuous solution curve} holds for some $M\leq n$.
\end{theorem}

\smallskip



\noindent\underline{$\phi$-divergence DRO:} Given a convex function $\phi$ on $[0,+\infty)$ such that $\phi(1)=0$, consider the $\phi$-divergence DRO formulation for \eqref{stoc_opt} in the form
\begin{equation}\label{DRO divergence general}
\min_{x\in\mathcal X}\ f(x)\text{\ \ subject to\ \ }\inf\Big\{\sum_{i=1}^nw_ih(x,\xi_i):\sum_{i=1}^n\frac{1}{n}\phi(nw_i)\leq s,\sum_{i=1}^nw_i=1,w_i\geq 0\text{ for all } i\Big\}\geq \gamma.
\end{equation}
We have the following result:
\begin{theorem}[Applying univariate Gaussian validator to $\phi$-divergence DRO]
\label{continuity of phi DRO}
Consider $OPT(s)$ using \eqref{DRO divergence general}. Suppose Assumptions \ref{convex decision space}-\ref{bounded level set} hold, and $\phi$ is continuous and convex on $(0,+\infty)$ with $\phi(1)=0$. In either of the following three cases:
\begin{enumerate}
\item[i.]$f(x)$ is continuous and strictly convex, $h(x,\xi)$ is continuous and concave in $x$ for every $\xi$.
\item[ii.]$f(x)=c'x$ for some non-zero $c\in\R^d$, the functions $f_r,r=1,\ldots,R$ in Assumption \ref{convex decision space} are strictly convex, any $k\leq d-1$ rows $\{w_{l(1)},\ldots,w_{l(k)}\}$ of $W$ do not satisfy the SCI condition, and $h(x,\xi)$ is continuous and strictly concave in $x$ for every $\xi$.
\item[iii.]Assume the same conditions as in (ii) except that $h(x,\xi)$ is only concave (instead of strictly concave) in $x$ for every $\xi$. In addition, $\phi$ is differentiable and strictly convex on $(0,+\infty)$ with $\lim_{x\to 0+}\phi(x)=+\infty$. For any $x_1,x_2\in \mathcal X$ let $\widehat{\mathrm{Corr}}(x_1,x_2)=\widehat{\mathrm{Cov}}(h(x_1,\xi),h(x_2,\xi))/(\hat\sigma(x_1)\hat\sigma(x_2))$ be the empirical correlation coefficient between $h(x_1,\xi)$ and $h(x_2,\xi)$ based on data $\{\xi_1,\ldots,\xi_n\}$. $\hat\sigma^2(x)>0$ for all $x\in\mathcal X$, and there exist no distinct $x_1, x_2$ such that $\widehat{\mathrm{Corr}}(\lambda x_1+(1-\lambda)x_2,x_1)=1$ for all $\lambda\in[0,1]$.
\end{enumerate}
Assumption \ref{piecewise uniform continuous solution curve} holds with $M=1$ for the $\phi$-divergence DRO conditioned on the data $\{\xi_1,\ldots,\xi_n\}$.
\end{theorem}

\smallskip

\noindent\underline{Wasserstein DRO:}
Consider the Wasserstein DRO reformulation for the constraint in \eqref{stoc_opt} given by
\begin{equation}\label{Wasserstein DRO general}
\min_{x\in\mathcal X}\ f(x)\text{\ \ subject to\ \ }\inf\Big\{\mathbb E_{G}[h(x,\xi)]:d_p(G,F_n )\leq s\Big\}\geq \gamma
\end{equation}
where $F_n=\frac{1}{n}\sum_{i=1}^n\delta_{\xi_i}$ is the empirical distribution and $d_p(G,F_n)$ is the Wasserstein distance between an arbitrary probability measure $G$ and $F_n$ which is defined as
\begin{equation*}
d_p^p(G,F_n)=\inf\Big\{\mathbb E_{\pi}[\norm{\xi-\xi'}^p]:\pi\text{ is a probability measure on }\Xi^2\text{ with marginals }G\text{ and }F_n\Big\}
\end{equation*}
where $\Xi$ is the known domain of $\xi$ and $\norm{\cdot}$ is an arbitrary norm.

The following theorem gives conditions under which Wasserstein DRO satisfies Assumption \ref{piecewise uniform continuous solution curve}:
\begin{theorem}[Applying univariate Gaussian validator to Wasserstein DRO]\label{continuity of Wasserstein DRO}
Consider $OPT(s)$ using \eqref{Wasserstein DRO general}. Suppose Assumptions \ref{convex decision space}-\ref{bounded level set} hold, the domain $\Xi$ of $\xi$ is compact, and $1\leq p<\infty$. In either of the following two cases:
\begin{enumerate}
\item[i.]$f(x)$ is continuous and strictly convex, $h(x,\xi)$ is jointly continuous in $x,\xi$ and concave in $x$ for every $\xi$.
\item[ii.]$f(x)=c'x$ for some non-zero $c\in\R^d$, the functions $f_r,r=1,\ldots,R$ in Assumption \ref{convex decision space} are strictly convex, any $k\leq d-1$ rows $\{w_{l(1)},\ldots,w_{l(k)}\}$ of $W$ do not satisfy the SCI condition, and $h(x,\xi)$ is jointly continuous in $x,\xi$ and strictly concave in $x$ for every $\xi$.
\end{enumerate}
Assumption \ref{piecewise uniform continuous solution curve} holds with $M=1$.
\end{theorem}
Proving Theorem \ref{continuity of Wasserstein DRO} requires utilizing the recently developed strong duality theory for Wasserstein DRO to show the joint continuity of the constraint \eqref{Wasserstein DRO general} and the existence of a worst-case distribution (e.g., \cite{blanchet2019quantifying}; \cite{gao2016distributionally}) to establish its strict convexity.

\smallskip

\noindent\underline{Moment-based DRO:}
We restrict our discussion in this case to individual linear chance constraints
\begin{equation}\label{joint linear chance constraint}
\mathbb P_F(a'_ix\leq b_i)\geq 1-\alpha_i,\text{\ \ for }i=1,\ldots,K
\end{equation}
where each $a_i$ is random and $b_i$ is a deterministic constant, and $\alpha_i$ is an individual tolerance level. This setup also applies to the case of joint linear chance constraint, say, $\mathbb P_F(a'_ix\leq b_i\text{ for }i=1,\ldots,K)\geq 1-\alpha_i$, where one uses the Bonferroni correction to safely approximate with $K$ single chance constraints $\mathbb P_F(a'_ix\leq b_i)\geq 1-\frac{\alpha}{K},i=1,\ldots,K$. We restrict our discussion to \eqref{joint linear chance constraint} as it enables the tractable use of moment-based DRO; other settings are possible, but would lead to much more elaborate technicality that we do not pursue here.

We consider for each single constraint the following moment-based distributionally robust counterpart
\begin{equation*}
    \inf_{a_i\sim Q\text{ s.t. }(\mathbb E_{Q}[a_i],\mathrm{Cov}_{Q}[a_i])\in \mathcal U_i(s)}\mathbb P_Q(a'_ix\leq b_i)\geq 1-\alpha_i\text{\ \ for }i=1,\ldots,K
\end{equation*}
where each $\mathcal U_i(s)$ is a joint uncertainty set for the mean and covariance of the uncertain quantity $a_i$, all parametrized by the same $s$. For a fixed mean $\mu$ and covariance $\Sigma$, the robust constraint $\inf_{a_i\sim Q\text{ s.t. }\mathbb E_Q[a_i]=\mu,\mathrm{Cov}_{Q}(a_i)=\Sigma}\mathbb P_Q(a'_ix\leq b_i)\geq 1-\alpha_i$ has an analytic expression $\mu'x+\sqrt{\frac{1-\alpha_i}{\alpha_i}}\sqrt{x'\Sigma x}\leq b_i$ (\cite{ghaoui2003worst}), therefore this moment-based DRO takes the form
\begin{equation}\label{moment DRO CCP}
\begin{aligned}
    &\min_{x\in\mathcal X}&&f(x)\\
    &\text{subject to}&&\sup_{(\mu,\Sigma)\in \mathcal U_i(s)}\mu'x+\sqrt{\frac{1-\alpha_i}{\alpha_i}}\sqrt{x'\Sigma x}\leq b_i\text{\ \ for }i=1,\ldots,K.
\end{aligned}
\end{equation}

\begin{theorem}[Applying univariate Gaussian validator to moment-based DRO]\label{continuity of moment DRO}
Consider $OPT(s)$ given by \eqref{moment DRO CCP}. Suppose Assumptions \ref{convex decision space}-\ref{bounded level set} hold, and that for each $i$ the uncertainty set $\mathcal U_i(s)$ satisfies either (1)(2)(3) or (1)(2)(4) among: (1) $\mathcal U_i(s)$ is compact for all $s$; (ii) $\mathcal U_i(s)\subseteq \mathcal U_i(s')$ whenever $s< s'$ and $\overline{\cup_{s'<s}\mathcal U_i(s')}=\cap_{s'>s}\mathcal U_i(s')=\mathcal U_i(s)$ for all $s$; (3) for every $s$ and every $(\mu,\Sigma)\in \mathcal U_i(s)$, $\Sigma$ is positive definite; (4) $\mathcal U_i(s)=\mathcal U^1_i(s)\times \mathcal U^2_i(s)$, where $\mathcal U^1_i(s)$ and $\mathcal U^2_i(s)$ are uncertainty sets for the mean and covariance respectively, and there is a positive definite $\Sigma_s\in \mathcal U^2_i(s)$ such that $\Sigma\preceq\Sigma_s$ for all $\Sigma\in \mathcal U^2_i(s)$ where $\preceq$ is the ordering with respect to the positive semi-definite cone. Then, in either of the following two cases:
\begin{enumerate}
\item[i.]$f(x)$ is continuous and strictly convex.
\item[ii.]$f(x)=c'x$ for some non-zero $c\in\R^d$, the functions $f_r,r=1,\ldots,R$ in Assumption \ref{convex decision space} are strictly convex, any $k\leq d-1$ rows $\{w_{l(1)},\ldots,w_{l(k)}\}$ of $W$ do not satisfy the SCI condition, and each $b_i\neq 0$.
\end{enumerate}
Assumption \ref{piecewise uniform continuous solution curve} holds with $M=1$.
\end{theorem}
Conditions (1) and (2) in Theorem \ref{continuity of moment DRO} hold for common choices of moment-based uncertainty sets. We discuss some examples where (3) and (4) arise. (3) holds when $\mathcal U_i(s)$ is constructed to be a joint confidence region from, e.g., the delta method (\cite{marandi2019extending}), for the mean and covariance of $a_i$ whose covariance component converges to the true positive definite covariance as data size grows. (4) happens if the mean and covariance are treated separately and the uncertainty set for covariance takes the form $\mathcal U^2_i(s)=\{\Sigma:\Sigma_l(s)\preceq \Sigma\preceq\Sigma_u(s)\}$ (e.g., \cite{delage2010distributionally}).

\smallskip

\noindent\underline{RO with polyhedral uncertainty set:}
Consider the same linear chance constraint \eqref{joint linear chance constraint}, and for each $i$ we use the robust counterpart $\sup_{a_i\in\mathcal U_i(s)}a'_ix\leq b_i$ where
\begin{equation*}
\mathcal U_i(s)=\{a_i:\mathcal W_ia_i\leq z_i+se_i\}
\end{equation*}
for some $\mathcal W_i\in \R^{l_i\times d},z_i\in\R^{l_i}$ and $e_i\in\R_+^{l_i}:=[0,\infty)^{l_i}$. This robust counterpart can be expressed as a set of linear constraints, leading to the following formulation
\begin{equation}\label{RO polyhedral chance constraint}
\begin{aligned}
    &\min_{x\in\mathcal X}&&f(x)\\
    &\text{subject to}&&(z_i+se_i)'y_i\leq b_i\\
    &&&W'_iy_i=x\\
    &&&y_i\geq 0\text{\ \  for all }i=1,\ldots,K
\end{aligned}
\end{equation}
where each $y_i\in\R^{l_i}$ is an auxiliary variable.

\begin{theorem}[Applying univariate Gaussian validator to polyhedral RO]\label{continuity of RO polyhedral}
Consider $OPT(s)$ given by \eqref{RO polyhedral chance constraint}. Suppose Assumptions \ref{convex decision space}-\ref{bounded level set} hold. If $f(x)$ is continuous and strictly convex, then Assumption \ref{piecewise uniform continuous solution curve} holds with $M=1$. Otherwise, if $f(x)=c'x$ for some non-zero $c\in\R^d$, $R=0$ in Assumption \ref{convex decision space}, the uncertainty set $\mathcal U_i(s_u)$ of maximal size is bounded for each $i$, and every $k\leq d-1$ element in $\{w_1,\ldots,w_L\}\cup\big( \cup_{i=1}^K\{\widetilde{\mathcal W}_i^{-1}z_i+s\widetilde{\mathcal W}_i^{-1}e_i:\widetilde{\mathcal W}_i\in\R^{d\times d}\text{ is an invertible submatrix of }\mathcal W_i\}\big)$ satisfies the SCI condition at only finitely many $s$ values, then Assumption \ref{piecewise uniform continuous solution curve} holds with some finite $M$.
\end{theorem}
The proof of Theorem \ref{continuity of RO polyhedral} involves some technical developments to show that $x^*(s)$ has left and right limits at each discontinuity. This consists of transforming \eqref{RO polyhedral chance constraint} into an equivalent parametric linear program whose constraints correspond to the vertices of the uncertainty sets, and then showing that its optimal basis stays constant in a neighborhood of each discontinuity. Lastly, we use the Jordan decomposition of the optimal basis matrix to establish the existence of left and right limits.

\smallskip



\noindent\underline{RO with ellipsoidal uncertainty set:}
Consider \eqref{joint linear chance constraint} again, and now for each constraint we consider using $\sup_{a_i\in\mathcal U_i(s)}a'_ix\leq b_i$ with
\begin{equation*}
\mathcal U_i(s)=\{a_i:a_i=\mu_i+\Sigma_iv,\norm{v}_2\leq s\}
\end{equation*}
for some positive definite $\Sigma_i\in \R^{d\times d}$, and $\mu_i\in\R^d$. This robust formulation has the following second-order cone representation
\begin{equation}\label{RO ellipsoid chance constraint}
\begin{aligned}
    &\min_{x\in\mathcal X}&&f(x)\\
    &\text{subject to}&&\mu'_ix+s\norm{\Sigma_ix}_2\leq b_i\text{\ \ for all }i=1,\ldots,K.
\end{aligned}
\end{equation}

\begin{theorem}[Applying univariate Gaussian validator to ellipsoidal RO]\label{continuity of RO ellipsoidal}
Consider $OPT(s)$ given by \eqref{RO ellipsoid chance constraint}. Suppose Assumptions \ref{convex decision space}-\ref{bounded level set} hold, and each $\Sigma_i$ is positive definite. In either of the following two cases:
\begin{enumerate}
\item[i.]$f(x)$ is continuous and strictly convex.
\item[ii.]$f(x)=c'x$ for some non-zero $c\in\R^d$, the functions $f_r,r=1,\ldots,R$ in Assumption \ref{convex decision space} are strictly convex, any $k\leq d-1$ rows $\{w_{l(1)},\ldots,w_{l(k)}\}$ of $W$ do not satisfy the SCI condition, and each $b_i\neq 0$.
\end{enumerate}
Assumption \ref{piecewise uniform continuous solution curve} holds with $M=1$.
\end{theorem}

\smallskip

\noindent\underline{SO:} Consider the CCP \eqref{chance constraint}. Given the data $\{\xi_1,\ldots,\xi_n\}$, consider the following sequence $OPT(s)$ of programs
\begin{equation}\label{SO formulation chance constraint}
\begin{aligned}
    &\min_{x\in\mathcal X}&&f(x)\\
    &\text{subject to}&&(x,\xi_i)\in A\text{ for all }i=1,\ldots,s
\end{aligned}
\end{equation}
for $1\leq s\leq n$, i.e., each $OPT(s)$ uses only the first $s$ sampled constraints. Although $s$ takes integer values only, we can artificially extend the solution path to the continuum $[1,n]$ without introducing new solutions, by letting $x^*(s)=x^*(i)$ for all $s\in[i,i+1)$. Like the SAA formulation for chance constraints, the solution path $x^*(s)$ can now be viewed as piecewise constant in $s\in[1,n]$ hence Assumption \ref{piecewise uniform continuous solution curve} holds. Therefore we have:
\begin{theorem}[Applying univariate Gaussian validator to SO]
Consider $OPT(s)$ given by \eqref{SO formulation chance constraint}. Assumption \ref{piecewise uniform continuous solution curve} holds for some $M\leq n$.
\end{theorem}

Lastly, our univariate Gaussian validator also works on a variant of SO called FAST (\cite{care2014fast}), in a sense that we will detail in Section \ref{sec:numerics on SO}. FAST differs from the formulations we have discussed so far in that its solution path does not come as solutions of a parametrized optimization problem, but from a line segment connecting two suitably chosen solutions. Nonetheless, the notion of solution-path optimality still applies. In particular, the solution-path optimum is unique if the objective is strictly convex or linear, and all the statistical guarantees in Theorem \ref{asymptotic joint:gaussian} can be established using the same proof.



\section{Numerical Experiments}\label{sec:numerics}
We present numerical results to demonstrate the performances of our framework in several data-driven reformulations. We consider the following linear CCP
\begin{equation}\label{ccp}
\min\  c'x\text{\ \ subject to\ \ }\mathbb P_F(\xi'x\leq b)\geq 1-\alpha
\end{equation}
where $c\in\R^d,b\in\R$ are deterministic, the distribution $F$ of the randomness $\xi\in\R^d$ is multivariate Gaussian with mean $\mu$ and covariance $\Sigma$, and the tolerance level $1-\alpha$ is set to $90\%$.

We consider a range of data-driven reformulations, including RO (or relatedly SCA), DRO (moment-based), and SO (including its variant FAST). In our experiments, we generate i.i.d.~data $\xi_1,\ldots,\xi_n$ from the underlying true distribution $F$. Then, using a chosen reformulation, we compute a solution $\hat x$ of \eqref{ccp} that attempts to satisfy the chance constraint with a $95\%$ confidence level, while attain an objective value $c'\hat x$ as low as possible. For each reformulation, we compare the performance of an existing benchmark with unnormalized and normalized Gaussian supremum validators (Algorithms \ref{algo:supremum_unnormalized} and \ref{algo:supremum_normalized}) and univariate Gaussian validator (Algorithm \ref{algo:marginal}), in terms of both feasibility and optimality. Moreover, we also test a naive validator that directly compares the sample mean to $\gamma$ when checking feasibility, i.e., without the Gaussian margin $\frac{z_{1-\beta}\hat{\sigma}_j}{\sqrt{n_2}}$ in \eqref{validate:marginal gaussian}, in addition to the three proposed validators, which serves to demonstrate the necessity of the proposed Gaussian margins in the validation procedure. The ``plain average'' column of each table displays results of this extra validator. ``unnorm. GS'' denotes the unnormalized Gaussian supremum validator (Algorithm \ref{algo:supremum_unnormalized}), ``norm. GS'' denotes the normalized Gaussian supremum validator (Algorithm \ref{algo:supremum_normalized}), and ``uni. Gaussian'' denotes the univariate Gaussian validator (Algorithm \ref{algo:marginal}). When applying these validators in all experiments, we use the simple allocation rule of dividing the overall data size into Phases 1 and 2 equally, except only in the case of basic SO where a too small Phase 1 data size is provably subpar in guaranteeing feasibility.

To collect statistically meaningful estimates, for each formulation we repeat the experiments $1000$ times each with an independently generated data set and a data-driven solution output. We take down the average objective value achieved by these solutions (the ``mean obj. val.'' row of each table) and the proportion of feasible solutions as the empirical feasibility coverage (the ``feasibility level'' row of each table). Therefore, the smaller the ``mean obj. val.'' is, the better is the solution in terms of optimality, and ``feasibility level'' $\geq 95\%$ indicates that the desired feasibility confidence level is achieved and otherwise not.

\subsection{RO and SCA}
We first test the proposed framework on RO. We use the ellipsoid uncertainty set that leads to a robust counterpart in the form described in Example \ref{example:RO and SCA}, i.e., $\hat{\mu}'x+\sqrt{s}\Vert\hat{\Sigma}^{1/2}x\Vert_2\leq b$ where $\hat{\mu}$ and $\hat{\Sigma}$ are the sample mean and covariance for $\xi$ computed from Phase one data. The benchmark (``SCA'' in the tables) is set to an SCA (equation 2.4.11 of \cite{ben2009robust}) for unbounded $\xi$, which in our case can be expressed as
\begin{equation}\label{SCA in our case}
\mu'x+\sqrt{2\log\frac{1}{\alpha}}\sqrt{\sum_{k=1}^d({z^k}'x)^2}=\mu'x+\sqrt{2\log\frac{1}{\alpha}}\Vert\Sigma^{1/2}x\Vert_2\leq b
\end{equation}
where $\mu$ is the true mean, and $z^k$ is the $k$-th column of the square root $\Sigma^{1/2}$ of the true covariance matrix $\Sigma$. Note that \eqref{SCA in our case} is equivalent to the RO formulation with true mean and covariance and parameter value $s=2\log\frac{1}{\alpha}$. Here, we give this SCA or RO the advantage of knowing the true mean and covariance of the randomness.

To implement our validator, we need to provide a set of parameter values $\{s_1,\ldots,s_p\}$ at which the RO is solved. We take the $(1-\alpha)n_1$-th order statistic $\hat s_{1-\alpha}$ of $\{(\xi_{n_2+i}-\hat\mu)'\hat{\Sigma}^{-1}(\xi_{n_2+i}-\hat\mu):i=1,\ldots,n_1\}$, where $\xi_{n_2+i},i=1,\ldots,n_1$ are the Phase one data, so that $\{\xi:(\xi-\hat\mu)'\hat{\Sigma}^{-1}(\xi-\hat\mu)\leq \hat s_{1-\alpha}\}$ is roughly a $(1-\alpha)$-content set for $\xi$ (such type of quantile-based selection has been used in \cite{hong2017learning}). We then set the values $s_j=(\hat s_{1-\alpha}+20)\frac{j}{50}$ for $j=1,\ldots,50$ ($p=50$). Tables \ref{RO_10_200}, \ref{RO_10_500} and \ref{RO_50_500} summarize the results under different problem dimensions and data sizes.

%

\begin{table}[h]
\caption{RO with ellipsoidal uncertainty set. $d=10,n=200$. Data are split to $n_1=100,n_2=100$.}
\begin{center}
\begin{tabular}{|c|c|c|c|c|c|}
\hline
&SCA&unnorm. GS&norm. GS&uni. Gaussian&plain average\\\hline
mean obj. val.&$-3.57$&$-3.68$&$-4.20$&$-4.43$&$-5.15$\\\hline
feasibility level&$100\%$&$99.9\%$&$98.5\%$&$97.5\%$&$76.9\%$\\\hline
\end{tabular}
\end{center}
\label{RO_10_200}
\end{table}%

\begin{table}[h]
\caption{RO with ellipsoidal uncertainty set. $d=10,n=500$. Data are split to $n_1=250,n_2=250$.}
\begin{center}
\begin{tabular}{|c|c|c|c|c|c|}
\hline
&SCA&unnorm. GS&norm. GS&uni. Gaussian&plain average\\\hline
mean obj. val.&$-3.57$&$-4.42$&$-4.58$&$-4.80$&$-5.34$\\\hline
feasibility level&$100\%$&$99.8\%$&$99.6\%$&$98.8\%$&$77.9\%$\\\hline
\end{tabular}
\end{center}
\label{RO_10_500}
\end{table}%

\begin{table}[h]
\caption{RO with ellipsoidal uncertainty set. $d=50,n=500$. Data are split to $n_1=250,n_2=250$.}
\begin{center}
\begin{tabular}{|c|c|c|c|c|c|}
\hline
&SCA&unnorm. GS&norm. GS&uni. Gaussian&plain average\\\hline
mean obj. val.&$-16.70$&$-17.59$&$-17.33$&$-17.71$&$-20.31$\\\hline
feasibility level&$100\%$&$98.4\%$&$99.6\%$&$98.4\%$&$82.7\%$\\\hline
\end{tabular}
\end{center}
\label{RO_50_500}
\end{table}%

We highlight a few observations. First, our framework with the three proposed validators outperforms the SCA benchmark. In terms of the objective performance, both our unnormalized and normalized Gaussian supremum validators, and univariate Gaussian validators, achieve lower objective value than SCA (with a difference $\geq 0.6$), while at the same time retain the feasibility confidence to above $95\%$ in all the three tables. In particular, as the dimension grows from $10$ (Tables \ref{RO_10_200} and \ref{RO_10_500}) to $50$ (Table \ref{RO_50_500}), the feasibility confidence level remains above $95\%$, consistent with the dimension-free feasibility guarantee of our methods. Second, among the three proposed validators, the univariate Gaussian validator appears less conservative than the Gaussian supremum counterparts in achieving better objective values, and relatedly tighter feasibility confidence levels (i.e., closer to $95\%$). Specifically, the univariate Gaussian validator gives a feasibility confidence level around $98\%$ in all the three tables, whereas the Gaussian supremum validators give a level between $99\%$-$100\%$ (and also $0.1$-$0.4$ higher mean objective values). Finally, we comment that the ``plain average'' scheme does not have the desired feasibility confidence level even when the data size is as large as $500$ (Table \ref{RO_10_500}), which shows that margin adjustments to the naive sample average in the validators is necessary to ensure feasibility.

\subsection{Moment-based DRO}
The second formulation we consider is a moment-based DRO. We use the formulation
\begin{equation}\label{moment DRO with confidence}
\inf_{\xi\sim Q\text{ s.t. }(\mathbb E_Q[\xi],\mathrm{Cov}_{Q}(\xi))\in\mathcal U_{s}}\mathbb P_Q(\xi'x\leq b)\geq 1-\alpha
\end{equation}
where $\mathcal U_s$ is a confidence region for the true mean and covariance of $\xi$ obtained via the delta method described in Example \ref{example:DRO} (see Section 6 of \cite{marandi2019extending} for details). According to \eqref{moment DRO CCP}, \eqref{moment DRO with confidence} can be expressed as $\sup_{(\mu,\Sigma)\in\mathcal U_{s}}\mu'x+\sqrt{\frac{1-\alpha}{\alpha}}\norm{\Sigma^{1/2}x}_2\leq b$, which can be further reformulated as a conic constraint (see Theorem 1 of \cite{marandi2019extending}). In the benchmark case ``DRO ($\chi^2$ quantile)'' we choose $s$ to be the $95\%$ quantile of the limiting $\chi^2$ distribution as suggested in \cite{marandi2019extending} so that $\mathcal U_s$ is a valid $95\%$ confidence region. In our framework, we solve the DRO formulation at parameter values $s_j=1.5\hat s_{0.95}\frac{j}{50}$ for $j=1,\ldots,50$ where $\hat s_{0.95}$ is the $\chi^2$ quantile used in the benchmark. Tables \ref{DRO_10_200} and \ref{DRO_10_500} show the experimental results under different data sizes.
\begin{table}[h]
\caption{Moment-based DRO. $d=10,n=200$. Data are split to $n_1=100,n_2=100$.}
\begin{center}
\begin{tabular}{|c|c|c|c|c|c|}
\hline
&DRO ($\chi^2$ quantile)&unnorm. GS&norm. GS&uni. Gaussian&plain average\\\hline
mean obj. val.&$-1.83$&$-2.73$&$-2.73$&$-2.73$&$-2.73$\\\hline
feasibility level&$100\%$&$100\%$&$100\%$&$100\%$&$100\%$\\\hline
\end{tabular}
\end{center}
\label{DRO_10_200}
\end{table}%

\begin{table}[h]
\caption{Moment-based DRO. $d=10,n=500$. Data are split to $n_1=250,n_2=250$.}
\begin{center}
\begin{tabular}{|c|c|c|c|c|c|}
\hline
&DRO ($\chi^2$ quantile)&unnorm. GS&norm. GS&uni. Gaussian&plain average\\\hline
mean obj. val.&$-2.00$&$-2.62$&$-2.62$&$-2.62$&$-2.62$\\\hline
feasibility level&$100\%$&$100\%$&$100\%$&$100\%$&$100\%$\\\hline
\end{tabular}
\end{center}
\label{DRO_10_500}
\end{table}%

The comparisons between the benchmark and our framework here share similarities with the RO setting. The solutions output from our validators possess superior objective performance (with a difference of $0.6$-$0.9$) than simply setting $s$ to be the $95\%$-level $\chi^2$ quantile, while still attain the desired feasibility confidence level. Note that all validators (including the ``plain average'') give the same objective value ($-2.73$ in Table \ref{DRO_10_200} and $-2.62$ in Table \ref{DRO_10_500}), and have a $100\%$ feasibility confidence. This is because the chosen parameter $s$ turns out to be 0 for all of them. In other words, setting the moment constraints as equalities (to the estimated moments from Phase one data) is statistically feasible and achieves the best objective value, and any relaxation from this would lead to a deterioration of solution quality. This hints that the conventional choices of moment set size suggested in the literature could be very conservative.

\subsection{SO}\label{sec:numerics on SO}
Given the Phase one data $\{\xi_{n_2+1},\ldots,\xi_{n}\}$, we consider the data-driven feasible region specified by the first $s$ sampled constraints, $\xi'_{n_2+i}x\leq b\text{ for }i=1,\ldots,s$, and tune the number of satisfied constraints $s\in\{1,2,\ldots,n_1\}$. The benchmark ``SO'' in this case is to impose all the constraints given by the whole data set $\{\xi_1,\ldots,\xi_n\}$. Tables \ref{SO_10_200} and \ref{SO_10_500} summarize the results for data size $n=200,500$ respectively.

\begin{table}[h]
\caption{SO. $d=10,n=200$. Data are split to $n_1=150,n_2=50$.}
\begin{center}
\begin{tabular}{|c|c|c|c|c|c|}
\hline
&SO&unnorm. GS&norm. GS&uni. Gaussian&plain average\\\hline
mean obj. val.&$-3.90$&$-4.24$&$-4.31$&$-4.46$&$-4.91$\\\hline
feasibility level&$99.7\%$&$95.2\%$&$94.0\%$&$85.1\%$&$44.7\%$\\\hline
\end{tabular}
\end{center}
\label{SO_10_200}
\end{table}%
\begin{table}[h]
\caption{SO. $d=10,n=500$. Data are split to $n_1=250,n_2=250$.}
\begin{center}
\begin{tabular}{|c|c|c|c|c|c|}
\hline
&SO&unnorm. GS&norm. GS&uni. Gaussian&plain average\\\hline
mean obj. val.&$-3.28$&$-3.86$&$-4.10$&$-4.30$&$-4.69$\\\hline
feasibility level&$100\%$&$99.7\%$&$98.7\%$&$95.6\%$&$62.0\%$\\\hline
\end{tabular}
\end{center}
\label{SO_10_500}
\end{table}%

We observe the gain in objective performance of our validators compared to SO (a difference of $0.3$-$0.6$ in Table \ref{SO_10_200} and $0.6$-$1.0$ in Table \ref{SO_10_500}). We also note the drastic failure of ``plain average'' in rendering the desired $95\%$ feasibility confidence, thus showing that a margin adjustment to the validators is necessary. Our validators maintain feasibility in all cases, except the univariate Gaussian validator for $n=200$. This deficiency is attributed to two potential reasons. First is that with $n_1=150$ there is a non-negligible chance that none of the $n_1$ solutions $x^*(s),s=1,\ldots,n_1$, produced in Phase one is feasible, thus violating Assumption \ref{non-binding}. In fact, the infeasibility probability of the solution derived by an SO using all the $n_1$ constraints can be computed to be $6\%$ (\cite{campi2008exact}), leaving the actual confidence of obtaining a feasible solution at most $94\%$. The second possible cause is the finite-sample coverage error of the univariate Gaussian validator, seeing that the validation data size $n_2=50$ is relatively small. When both $n_1$ and $n_2$ increase to $250$ in Table \ref{SO_10_500}, the desired feasibility confidence level is recovered for the univariate Gaussian validator as the chance of all solution candidates being infeasible decreases to $<0.2\%$ and the finite-sample error is reduced due to a larger validation data size. Finally, although we do not pursue here, we should mention that the performances of the basic SO considered in the tables can plausibly be boosted by using techniques such as sampling-and-discarding (\cite{campi2011sampling}) and wait-and-judge (\cite{campi2018wait}). Comprehensive comparisons with these enhanced techniques would be left as important future work.

Lastly, we consider a variant of SO called FAST (\cite{care2014fast}), designed originally to tone down the sample size requirement in basic SO. Our comparison with FAST here is motivated by its similarity with our framework in that it also splits the data into two portions and uses a validation-based idea. With the first portion of data $\{\xi_{n_2+i},i=1,\ldots,n_1\}$, FAST computes a solution $\hat x$ by imposing all the $n_1$ constraints $\xi'_{n_2+i}x\leq b$ as in the basic SO, and then uses the second portion to obtain the final solution $\hat x^*$ by solving the following program
\begin{equation*}
\min\  c'((1-s)x_o+s\hat x)\text{\ \ subject to\ \ }\xi'_i((1-s)x_o+s\hat x)\leq b\text{ for all }i=1,\ldots,n_2\text{ and }0\leq s\leq 1
\end{equation*}
where $x_o$ is a feasible solution of \eqref{ccp} with $\mathbb P_F(\xi'x_o\leq b)=1$. One particular choice of $x_o$ for problem \eqref{ccp} is the vector of all zeros and is used in the experiment. When applying our framework to FAST, we search for the best feasible solution along the line segment $\{x^*(s)=(1-s)x_o+s\hat x:s\in[0,1]\}$ by validating solutions $x^*(s_j)$ at parameter values $s_j=\frac{j-1}{10}$ for $j=1,\ldots,11$ ($p=11$). Tables \ref{FAST_10_200} and \ref{FAST_50_500} show the results under different dimensions and data sizes.

\begin{table}[h]
\caption{FAST. $d=10,n=200$. Data are split to $n_1=100,n_2=100$.}
\begin{center}
\begin{tabular}{|c|c|c|c|c|c|}
\hline
&FAST&unnorm. GS&norm. GS&uni. Gaussian&plain average\\\hline
mean obj. val.&$-2.54$&$-3.55$&$-3.68$&$-3.87$&$-4.44$\\\hline
feasibility level&$100\%$&$98.9\%$&$98.9\%$&$97.3\%$&$79.6\%$\\\hline
\end{tabular}
\end{center}
\label{FAST_10_200}
\end{table}%

\begin{table}[h]
\caption{FAST. $d=50,n=500$. Data are split to $n_1=250,n_2=250$.}
\begin{center}
\begin{tabular}{|c|c|c|c|c|c|}
\hline
&FAST&unnorm. GS&norm. GS&uni. Gaussian&plain average\\\hline
mean obj. val.&$-8.92$&$-14.11$&$-15.06$&$-15.80$&$-18.14$\\\hline
feasibility level&$100\%$&$99.8\%$&$99.3\%$&$98.0\%$&$76.7\%$\\\hline
\end{tabular}
\end{center}
\label{FAST_50_500}
\end{table}%

Similar phenomena persist from our previous settings. Our three validators give tighter feasibility confidence levels and better objective performances (with a difference of $\geq 1$ in Table \ref{FAST_10_200} and $\geq 5$ in Table \ref{FAST_50_500}) compared to FAST. Among them, univariate Gaussian validator gives the tightest feasibility confidence level and best objective value. The naive ``plain average" validator fails in attaining the desired feasibility confidence. Here we have used a rather coarse mesh with only $11$ parameter values, and expect a sharper improvement should a finer mesh be used.

\section{Conclusion}
We have studied a validation-based framework to combat the conservativeness in data-driven optimization with uncertain constraints. We have demonstrated how the conventional approaches in several optimization paradigms, including SAA, RO and DRO, implicitly estimate the whole feasible region. This in turn leads to over-conservativeness caused by the need to control huge simultaneous estimation errors, especially for high-dimensional problems. On the other hand, we have also demonstrated that the solution output from these reformulation classes can often be represented in a low-dimensional manifold parametrized by key conservativeness parameters.
Our framework leverages this low dimensionality by extracting the parametrized solution path and selecting the best parameter value. We have proposed two types of validators for this parameter selection, one utilizing a multivariate Gaussian supremum (unnormalized or normalized) and another utilizing a univariate Gaussian, to set the margin in a sample average constraint when optimizing over the solution path. We have shown that the obtained solutions enjoy asymptotic and finite-sample performance guarantees on feasibility that scale lightly with the problem dimension, and asymptotic optimality within the reformulation class. The Gaussian supremum validator requires less regularity conditions and is applicable more generally, whereas the univariate Gaussian validator provides tighter guarantees when applicable. Our numerical results support these findings and show that our framework and validators consistently provide better solutions  compared to several benchmarks in terms of better objective values and tighter feasibility confidence. Our study provides a first rigorous validation-based framework to tackle over-conservativeness in data-driven constrained optimization, and is foreseen to open up follow-up investigations on more powerful validation strategies and refined statistical guarantees regarding joint feasibility and optimality.

\ACKNOWLEDGMENT{We gratefully acknowledge support from the National Science Foundation under grants CAREER CMMI-1653339/1834710 and IIS-1849280.}


\bibliographystyle{informs2014} 
\bibliography{reference} 

\ECSwitch


\ECHead{Proofs of Statements}
In all the proofs, for universal constants which are usually denoted $C$ or $c$, we abuse notation slightly to allow $C$ or $c$ to take a different value in each appearance. For example, consider three quantities $x,y,z$ such that $x\leq Cy$ and $z\leq 2x$. This implies $z\leq 2Cy$, but we would write as $z\leq Cy$ to simplify the notation.

\section{Existing Central Limit Theorems in High Dimensions}\label{sec:high-dimensional CLTs}
This section reviews some results on high-dimensional central limit theorems that are needed subsequently in our proofs.
We start with some notations. Let $\mathbf X_{i}:=(X_{i,1},\ldots,X_{i,p}), i=1,\ldots,n$ be $n$ i.i.d.~copies of the random vector $\mathbf X:=(X_1,\ldots,X_p)\in\R^p$, and $\mu_j:=E[X_j]$ for $j=1,\ldots,p$. Let $\bar{X}_j=\sum_{i=1}^nX_{i,j}/n$ be the sample mean of the $j$-th component. We denote by $\mathbf Z:=(Z_1,\ldots,Z_p)$ a $p$-dimensional Gaussian random vector with $E[Z_j]=0$ and covariance structure $Cov(Z_{j},Z_{j'})=\Sigma(j,j'):=Cov(X_{j},X_{j'})$ for $j,j'=1,\ldots,p$, and by $\widehat{\mathbf Z}:=(\widehat{Z}_1,\ldots,\widehat{Z}_p)$ a $p$-dimensional centered Gaussian random vector with covariance $\widehat{\Sigma}$, where
$$\widehat{\Sigma}(j,j')=\frac{1}{n}\sum_{i=1}^nX_{i,j}X_{i,j'}-\bar{X}_j\bar{X}_{j'}$$
is the sample covariance of all $\mathbf X_i$'s. We also denote $\sigma^2_j=\Sigma(j,j)$ and $\hat\sigma^2_j=\widehat\Sigma(j,j)$.

We make the following assumption:
\begin{assumption}\label{condition:clt1}
There exist constants $b>0$ and $B\geq 1$ such that
\begin{align*}
&Var[X_{j}]\geq b\text{ and }E[\exp(\abs{X_{j}-\mu_j}^2/B^2)]\leq 2\text{ for all }j=1,\ldots,p\\
&E[\abs{X_{j}-\mu_j}^{2+k}]\leq B^k\text{ for all }j=1,\ldots,p\text{ and }k=1,2.
\end{align*}
\end{assumption}
Note that, since the sub-exponential norm of a random variable is always bounded above by its sub-Gaussian norm up to some universal constant $C$, the exponential condition in Assumption \ref{condition:clt1} implies $E[\exp(\abs{X_{j}-\mu_j}/(CB))]\leq 2$. \cite{chernozhukov2017central} proved the following CLT:
\begin{theorem}[First half of Proposition 2.1 in \cite{chernozhukov2017central}]\label{clt1}
If Assumption \ref{condition:clt1} holds, then
\begin{align*}
\sup_{a_j\leq b_j,j=1,\ldots,p}\abs{P(a_j\leq \sqrt n(\bar{X}_j-\mu_j)\leq b_j\text{ for all }j)-P(a_j\leq Z_j\leq b_j\text{ for all }j)}\leq C_1\Big(\frac{B^2\log^7(pn)}{n}\Big)^{\frac{1}{6}}
\end{align*}
where the constant $C_1$ depends only on $b$.
\end{theorem}
To derive confidence bounds based on the CLT, one needs to properly estimate the quantile of the limit Gaussian vector $\mathbf Z\sim N_p(0,\Sigma)$. One common approach is to use the Gaussian vector $\widehat{\mathbf Z}\sim N_p(0,\widehat{\Sigma})$, where $\widehat{\Sigma}$ is the sample covariance matrix, to approximate $\mathbf Z$. This approach is also called the multiplier bootstrap. \cite{chernozhukov2017central} gave the following result concerning the statistical accuracy of the multiplier bootstrap:
\begin{theorem}[First half of Corollary 4.2 in \cite{chernozhukov2017central}]\label{plugin1}
If Assumption \ref{condition:clt1} holds, then for any constant $0<\alpha<\frac{1}{e}$ we have
\begin{align*}
\sup_{a_j\leq b_j,j=1,\ldots,p}\Big\lvert P(a_j\leq \widehat{Z}_j\leq b_j\text{ for all }j\vert \{\mathbf X_i\}_{i=1}^n)-P(a_j\leq Z_j\leq b_j\text{ for all }j)\Big\rvert\leq C_2\Big(\frac{B^2\log^5(pn)\log^2(1/\alpha)}{n}\Big)^{\frac{1}{6}}
\end{align*}
with probability at least $1-\alpha$, where the constant $C_2$ depends only on $b$.
\end{theorem}

\section{Proofs of Results in Section \ref{sec:gaussian supremum margin}}
This section proves the performance guarantees of our Gaussian supremum validators. Section \ref{clt:adaptation} adapts the high-dimensional CLTs in Appendix \ref{sec:high-dimensional CLTs} to handle small-variance situations that potentially arise in our optimization context. Section \ref{clt:extension to normalized case} extends them to the case where the sample means are normalized by their standard deviations, a key step to justify our normalized Gaussian supremum validator. Section \ref{multiplier bootstrap:unnormalized and normalized} presents results on the consistency of the multiplier bootstrap to approximate the limiting Gaussian distributions. Section \ref{main proof: gaussian supremum validators} puts together all these results to synthesize the main proofs for Section \ref{sec:gaussian supremum margin}.
\subsection{A CLT for Random Vectors with Potentially Small Variances}\label{clt:adaptation}
Note that in both Theorems \ref{clt1} and \ref{plugin1}, the constants $C_1,C_2$ depend on the minimal variance $b$. By tracing the proof in \cite{chernozhukov2017central}, the constant $C_1$ is of the form $c_1(b^{-1}\vee c_2)$ where $c_1,c_2$ are two universal constants. Due to such a dependence on the minimal variance, the bound can deteriorate when the noise levels across different components of $\mathbf X$ are not of the same scale, e.g., in the case of CCPs. To resolve this issue, we derive an alternate CLT that applies to normalized random vectors.
We assume:
\begin{assumption}\label{condition:clt2}
$Var[X_{j}]>0$ for all $j=1,\ldots,p$ and there exists some constant $D_1\geq 1$ such that
\begin{align}
&E\big[\exp\Big(\frac{\abs{X_{j}-\mu_j}^2}{D_1^2Var[X_j]}\Big)\big]\leq 2\text{ for all }j=1,\ldots,p\label{orlicz}\\
&E\big[\Big(\frac{\abs{X_{j}-\mu_j}}{\sqrt{Var[X_j]}}\Big)^{2+k}\big]\leq D_1^k\text{ for all }j=1,\ldots,p\text{ and }k=1,2.\label{skewness-kurtosis}
\end{align}
\end{assumption}
Note that rectangles in $\R^p$ are invariant with respect to component-wise rescaling, i.e., for any rectangle $R=\{(x_1,\ldots,x_p):a_j\leq x_j\leq b_j,j=1,\ldots,p\}$, the rescaled set $R':=\{(\lambda_1x_1,\ldots,\lambda_px_p):(x_1,\ldots,x_p)\in R\}$ with each $\lambda_j>0$ is still a rectangle which can be represented as $R'=\{(x_1,\ldots,x_p):\lambda_ja_j\leq x_j\leq \lambda_jb_j,j=1,\ldots,p\}$. Hence one can show the following CLT by applying Theorem \ref{clt1} to the rescaled data:
\begin{theorem}\label{clt2}
If Assumption \ref{condition:clt2} holds, then
\begin{align*}
\sup_{a_j\leq b_j,j=1,\ldots,p}\Big\lvert P(a_j\leq \sqrt n(\bar{X}_j-\mu_j)\leq b_j\text{ for all }j)-P(a_j\leq Z_j\leq b_j\text{ for all }j)\Big\rvert\leq C\Big(\frac{D_1^2\log^7(pn)}{n}\Big)^{\frac{1}{6}}
\end{align*}
where $C$ is a universal constant.
\end{theorem}
\proof{Proof of Theorem \ref{clt2}.} Consider the rescaled data $Y_{i,j}=(X_{i,j}-\mu_j)/\sqrt{Var[X_j]}$. Due to Assumption \ref{condition:clt2}, $Y_{i,j}$'s satisfy Assumption \ref{condition:clt1} with $b=1$ and $B=D_1$, and has covariance structure $\Sigma_Y(j,j')=\Sigma(j,j')/\sqrt{\Sigma(j,j)\Sigma(j',j')}$. Let $\bar{Y}_j=\sum_{i=1}^nY_{i,j}/n$. By Theorem \ref{clt1} we have
\begin{align*}
\sup_{a_j\leq b_j,j=1,\ldots,p}\Big\lvert P(a_j\leq \sqrt n\bar{Y}_j\leq b_j\text{ for all }j)-P(a_j\leq \frac{Z_j}{\sqrt{Var[X_j]}}\leq b_j\text{ for all }j)\Big\rvert\leq C\Big(\frac{D_1^2\log^7(pn)}{n}\Big)^{\frac{1}{6}}.
\end{align*}
The theorem follows from
\begin{align*}
&\sup_{a_j\leq b_j,j=1,\ldots,p}\Big\lvert P(a_j\leq \sqrt n(\bar{X}_j-\mu_j)\leq b_j\text{ for all }j)-P(a_j\leq Z_j\leq b_j\text{ for all }j)\Big\rvert\\
=&\sup_{a_j\leq b_j,j=1,\ldots,p}\Big\lvert P(\sqrt{Var[X_j]}a_j\leq \sqrt n(\bar{X}_j-\mu_j)\leq \sqrt{Var[X_j]}b_j\text{ for all }j)\\
&\hspace{14ex}-P(\sqrt{Var[X_j]}a_j\leq Z_j\leq \sqrt{Var[X_j]}b_j\text{ for all }j)\Big\rvert\\
=&\sup_{a_j\leq b_j,j=1,\ldots,p}\Big\lvert P(a_j\leq \sqrt n\bar{Y}_j\leq b_j\text{ for all }j)-P(a_j\leq \frac{Z_j}{\sqrt{Var[X_j]}}\leq b_j\text{ for all }j)\Big\rvert.
\end{align*}\hfill\Halmos

Similarly, we have the following result regarding the multiplier bootstrap:
\begin{theorem}\label{plugin2}
If Assumption \ref{condition:clt2} holds, then for any constant $0<\alpha<\frac{1}{e}$ we have
\begin{align*}
\sup_{a_j\leq b_j,j=1,\ldots,p}\Big\lvert P(a_j\leq \widehat{Z}_j\leq b_j\text{ for all }j\vert \{\mathbf X_i\}_{i=1}^n)-P(a_j\leq Z_j\leq b_j\text{ for all }j)\Big\rvert\leq C\Big(\frac{D_1^2\log^5(pn)\log^2(1/\alpha)}{n}\Big)^{\frac{1}{6}}
\end{align*}
with probability at least $1-\alpha$, where $C$ is a universal constant.
\end{theorem}
\proof{Proof of Theorem \ref{plugin2}.} Again, consider the rescaled data $Y_{i,j}=(X_{i,j}-\mu_j)/\sqrt{Var[X_j]}$. Note that the sample covariance of $Y_{i,j}$ is the same as the covariance structure of $(\widehat{Z}_1/\sqrt{Var[X_1]},\ldots,\widehat{Z}_m/\sqrt{Var[X_p]})$. Theorem \ref{plugin2} entails that
\begin{align*}
&\sup_{a_j\leq b_j,j=1,\ldots,p}\Big\lvert P(a_j\leq \frac{\widehat{Z}_j}{\sqrt{Var[X_j]}}\leq b_j\text{ for all }j\vert \{\mathbf X_i\}_{i=1}^n)-P(a_j\leq \frac{Z_j}{\sqrt{Var[X_j]}}\leq b_j\text{ for all }j)\Big\rvert\\
\leq& C\Big(\frac{D_1^2\log^5(pn)\log^2(1/\alpha)}{n}\Big)^{\frac{1}{6}}
\end{align*}
with probability at least $1-\alpha$. The desired conclusion then follows by invariance of the class of rectangles under component-wise rescaling.\hfill\Halmos

Theorems \ref{clt2} and \ref{plugin2} rely on conditions more pertinent to our context than those in Theorems \ref{clt1} and \ref{plugin1}. The first condition \eqref{orlicz} of Assumption \ref{condition:clt2} measures the ratio of the sub-exponential norm to the $L_2$ norm of each component of the vector, whereas the second condition \eqref{skewness-kurtosis} concerns the kurtosis of each component. Therefore, to guarantee a valid CLT, we need the vector to be sufficiently light-tailed after being normalized to have unit variance.

\subsection{CLTs for Sample Means Normalized by Standard Deviations}\label{clt:extension to normalized case}

We establish CLTs for sample means normalized by sample standard deviations, needed to prove results regarding our normalized supremum validator. Note that when the dimension $p$ is fixed, such CLTs can be established by Slutsky's theorem, but when $p$ is huge or grows with the data size $n$ this is no longer applicable. Instead, we need to develop concentration inequalities for sample variances, which we state below.
\begin{lemma}[Concentration inequalities for sample variances]\label{concentration:variance}
Let $\xi_1,\ldots,\xi_n$ be $n$ i.i.d.~copies of the random variable $\xi\in\R$, $\sigma^2:=Var[\xi]$ be the true variance, and $\hat\sigma^2=\sum_{i=1}^n(\xi_i-\bar{\xi})^2/n$, where $\bar{\xi}=\sum_{i=1}^n\xi_i/n$ is the sample mean, be the sample variance. We have the following concentration inequalities:
\begin{enumerate}
\item[1.]if $\xi$ is $[0,1]$-valued, then there exists a universal constant $C$ such that for any $t>0$
\begin{equation}\label{case:ccp}
P(\lvert\hat\sigma^2-\sigma^2\rvert>t)\leq 2\exp\big(-\frac{Cnt^2}{\sigma^2+t}\big).
\end{equation}
\item[2.]if $\xi-E[\xi]$ has a sub-Gaussian norm at most $K$, i.e.~$E[\exp((\xi-E[\xi])^2/K^2)]\leq 2$, then there exists a universal constant $C$ such that for any $t>0$
\begin{equation}\label{case:generalcvx}
P(\lvert\hat\sigma^2-\sigma^2\rvert>t)\leq 4\exp(-\frac{Cnt^2}{K^4+K^2t}).
\end{equation}
\end{enumerate}

\end{lemma}
\proof{Proof of Lemma \ref{concentration:variance}.}\textbf{Case 1:} Since the unbiased sample variance, $(n/(n-1))\hat\sigma^2$, is a U-statistic of degree $2$, Hoeffding's inequality for U-statistics (see, e.g., \cite{hoeffding1963probability,peel2010empirical}) entails that with probability at least $1-\alpha$
\begin{align*}
\lvert\frac{n}{n-1}\hat\sigma^2-\sigma^2\rvert\leq \sqrt{\frac{4Var[(\xi-\xi')^2]}{n}\log \frac{2}{\alpha}}+\frac{4}{3n}\log\frac{2}{\alpha}
\end{align*}
where $\xi,\xi'$ are i.i.d.~copies. Note that $Var[(\xi-\xi')^2]\leq E[(\xi-\xi')^4]\leq E[(\xi-\xi')^2]=2\sigma^2$ because $\abs{\xi-\xi'}\leq 1$, and $\hat\sigma^2\leq 1$ for the same reason. Therefore with probability at least $1-\alpha$
\begin{eqnarray*}
\lvert\hat\sigma^2-\sigma^2\rvert&\leq& \sqrt{\frac{8\sigma^2}{n}\log \frac{2}{\alpha}}+\frac{4}{3n}\log\frac{2}{\alpha}+\frac{\hat\sigma^2}{n-1}\\
&\leq& \sqrt{\frac{8\sigma^2}{n}\log \frac{2}{\alpha}}+\frac{5}{n}\log\frac{2}{\alpha}
\end{eqnarray*}
and the conclusion easily follows by fixing the right hand side and solving for $\alpha$.

\textbf{Case 2:} Since the sub-Gaussian norm of $\xi-E[\xi]$ is at most $K$, $(\xi-E[\xi])^2$ has a sub-exponential norm of at most $K^2$ by definition. Centering a variable can only inflate its sub-exponential norm by a constant factor (Remark 5.18 in \cite{vershynin2010introduction}), that is, $(\xi-E[\xi])^2-\sigma^2$ must have a sub-exponential norm of at most $CK^2$ where $C$ is a universal constant. By Hoeffding's inequality and Bernstein's inequality for sums of independent variables (Propositions 5.10 and 5.16 in \cite{vershynin2010introduction}) we have for some universal constant $C$ and any $t>0$
\begin{align*}
&P\big(\big\lvert\frac{1}{n}\sum_{i=1}^n(\xi_i-E[\xi])^2-\sigma^2\big\rvert>t\big)\leq 2\exp(-\frac{Cnt^2}{K^4+K^2t})\\
&P\big(\big\lvert\frac{1}{n}\sum_{i=1}^n\xi_i-E[\xi]\big\rvert>t\big)\leq 2\exp(-\frac{Cnt^2}{K^2}).
\end{align*}
Note that the sample variance can be expressed as
\begin{align*}
\hat\sigma^2=\frac{1}{n}\sum_{i=1}^n(\xi_i-E[\xi])^2-(\frac{1}{n}\sum_{i=1}^n\xi_i-E[\xi])^2.
\end{align*}
Hence by a union bound
\begin{align*}
P(\abs{\hat\sigma^2-\sigma^2}>t)&\leq P\big(\big\lvert\frac{1}{n}\sum_{i=1}^n(\xi_i-E[\xi])^2-\sigma^2\big\rvert>t/2\big)+P\big(\big\lvert\frac{1}{n}\sum_{i=1}^n\xi_i-E[\xi]\big\rvert>\sqrt{t/2}\big)\\
&\leq 2\exp(-\frac{Cnt^2}{K^4+K^2t})+2\exp(-\frac{Cnt}{K^2})\\
&\leq 4\exp(-\frac{Cnt^2}{K^4+K^2t}).
\end{align*}
This completes the proof.\hfill\Halmos

Note that inequality \eqref{case:ccp} cannot be deduced from inequality \eqref{case:generalcvx} as a special case because of the appearance of $\sigma^2$ in the bound. In fact \eqref{case:ccp} is a sharper bound than \eqref{case:generalcvx} when the variable $\xi$ is Bernoulli, e.g., in the case of chance constrained optimization, because the sub-Gaussian norm of $\xi-E[\xi]$ is of order $K^2=\Theta(\log^{-1}(1/\epsilon))\gg \Theta(\epsilon)=\sigma^2$ when the success probability $\epsilon$ is small.

We also need the following anti-concentration inequality for Gaussian distribution:
\begin{lemma}[Nazarov's inequality]\label{anti-concentration}
Let $(Y_1,\ldots,Y_p)$ be an $p$-dimensional centered Gaussian random vector such that $Var[Y_j]\geq b$ for all $j=1,\ldots,p$ and some constant $b>0$. Then for every $-\infty\leq a_j\leq b_j\leq +\infty,j=1,\ldots,p$ and every $\delta>0$ it holds
\begin{equation*}
P(a_j-\delta\leq Y_j\leq b_j+\delta\text{ for all }j)-P(a_j\leq Y_j\leq b_j\text{ for all }j)\leq C_3\delta\sqrt{\log p}
\end{equation*}
where $C_3$ is a constant that depends only on $b$.
\end{lemma}
A special case of this inequality where $a_j=-\infty$ for all $j=1,\ldots,p$ has appeared in \cite{chernozhukov2017central}. Establishing a similar inequality for the case of possibly finite $a_j$'s involves a routine application of union bound. For completeness we provide a proof here.
\proof{Proof of Lemma \ref{anti-concentration}.}Lemma A.1 in \cite{chernozhukov2017central} states that for every $\delta>0$ and every $b_1,\ldots,b_p$ the following bound holds
\begin{equation*}
P(Y_j\leq b_j+\delta\text{ for all }j)-P(Y_j\leq b_j\text{ for all }j)\leq C\delta\sqrt{\log p}
\end{equation*}
where $C$ depends on $b$ only. Applying the same bound to $-Y_j,j=1,\ldots,p$ and $-a_j,j=1,\ldots,p$ gives
\begin{equation*}
P(a_j-\delta\leq Y_j\text{ for all }j)-P(a_j\leq Y_j\text{ for all }j)\leq C\delta\sqrt{\log p}.
\end{equation*}
Therefore
\begin{eqnarray*}
&&P(a_j-\delta\leq Y_j\leq b_j+\delta\text{ for all }j)-P(a_j\leq Y_j\leq b_j\text{ for all }j)\\
&=&P(a_j-\delta\leq Y_j\leq b_j+\delta\text{ for all }j)-P(a_j-\delta\leq Y_j\leq b_j\text{ for all }j)\\
&&\hspace{2em}+P(a_j-\delta\leq Y_j\leq b_j\text{ for all }j)-P(a_j\leq Y_j\leq b_j\text{ for all }j)\\
&\leq&P(Y_j\leq b_j+\delta\text{ for all }j)-P( Y_j\leq b_j\text{ for all }j)\\
&&\hspace{2em}+P(a_j-\delta\leq Y_j\text{ for all }j)-P(a_j\leq Y_j\text{ for all }j)\\
&\leq&2C\delta\sqrt{\log p}.
\end{eqnarray*}
This completes the proof.\hfill\Halmos


We have the following CLT with componentwise normalization. Recall that $\hat\sigma_j^2$ is the sample variance of $X_j$ computed from the data $\{X_{1,j},\ldots,X_{n,j}\}$.
\begin{theorem}\label{clt3:generalcvx}
Under Assumption \ref{condition:clt2} we have
\begin{align*}
&\sup_{a_j\leq b_j,j=1,\ldots,p}\Big\lvert P(\hat\sigma_ja_j\leq \sqrt n(\bar{X}_j-\mu_j)\leq \hat\sigma_jb_j\text{ for all }j)-P(\sigma_ja_j\leq Z_j\leq \sigma_jb_j\text{ for all }j)\Big\rvert\\
\leq&C\Big(\frac{D_1^2\log^7(pn)}{n}\Big)^{\frac{1}{6}}+Cp\exp\big(-\frac{cn^{2/3}}{D_1^{10/3}}\big)
\end{align*}
for some universal constants $C,c$.
\end{theorem}

If each component of the random vector is $[0,1]$-valued, we assume:
\begin{assumption}\label{condition:clt3_ccp}
Each $X_j$ is $[0,1]$-valued and $\sigma_j^2:=Var[X_j]\geq \delta$ for all $j=1,\ldots,p$ and some constant $\delta>0$.
\end{assumption}
Then we have an alternate CLT:
\begin{theorem}\label{clt3:ccp}
Under Assumptions \ref{condition:clt2} and \ref{condition:clt3_ccp} we have
\begin{align*}
&\sup_{a_j\leq b_j,j=1,\ldots,p}\Big\lvert P(\hat\sigma_ja_j\leq \sqrt n(\bar{X}_j-\mu_j)\leq \hat\sigma_jb_j\text{ for all }j)-P(\sigma_ja_j\leq Z_j\leq \sigma_jb_j\text{ for all }j)\Big\rvert\\
\leq&C\Big(\frac{D_1^2\log^7(pn)}{n}\Big)^{\frac{1}{6}}+Cp\exp\big(-c\delta D_1^{2/3}n^{2/3}\big)
\end{align*}
for some universal constants $C,c$.
\end{theorem}
\proof{Proof of Theorems \ref{clt3:generalcvx} and \ref{clt3:ccp}.} For any $a_j\leq b_j,j=1,\ldots,p$ and $0<\epsilon<1/2$
\begin{align*}
&P(\hat\sigma_ja_j\leq \sqrt n(\bar{X}_j-\mu_j)\leq\hat\sigma_j b_j\text{ for all }j)\\
=&P(\frac{\hat\sigma_j}{\sigma_j}a_j\leq \frac{\sqrt n(\bar{X}_j-\mu_j)}{\sigma_j}\leq\frac{\hat\sigma_j}{\sigma_j} b_j\text{ for all }j)\\
=&P(\frac{\hat\sigma_j}{\sigma_j}a_j\leq \frac{\sqrt n(\bar{X}_j-\mu_j)}{\sigma_j}\leq\frac{\hat\sigma_j}{\sigma_j} b_j,\big\lvert\frac{\hat\sigma_j}{\sigma_j}-1\big\rvert\leq \epsilon\text{ for all }j)\\
&+P(\frac{\hat\sigma_j}{\sigma_j}a_j\leq \frac{\sqrt n(\bar{X}_j-\mu_j)}{\sigma_j}\leq\frac{\hat\sigma_j}{\sigma_j} b_j\text{ for all }j,\big\lvert\frac{\hat\sigma_j}{\sigma_j}-1\big\rvert> \epsilon\text{ for some }j)\\
\leq &P(a_j-\epsilon \abs{a_j}\leq \frac{\sqrt n(\bar{X}_j-\mu_j)}{\sigma_j}\leq b_j+\epsilon \abs{b_j}\text{ for all }j)+P(\big\lvert\frac{\hat\sigma_j}{\sigma_j}-1\big\rvert> \epsilon\text{ for some }j)\\
\leq & P(a_j-\epsilon \abs{a_j}\leq \frac{Z_j}{\sigma_j}\leq b_j+\epsilon \abs{b_j}\text{ for all }j)+C\Big(\frac{D_1^2\log^7(pn)}{n}\Big)^{\frac{1}{6}}+\sum_{j=1}^pP(\big\lvert\hat\sigma_j-\sigma_j\big\rvert> \epsilon\sigma_j)
\end{align*}
where in the last inequality we use Theorem \ref{clt2} for the first probability and apply the union bound to the second probability.
Note that $\epsilon<1/2$ hence for any constant $M>0$ if we denote by $f_M(x)=-\infty\cdot \mathbf{1}(x<-M)+\infty\cdot \mathbf{1}(x>M)$ then we have
\begin{align*}
&P(a_j-\epsilon \abs{a_j}\leq \frac{Z_j}{\sigma_j}\leq b_j+\epsilon \abs{b_j}\text{ for all }j)\\
\leq &P(f_M(a_j)+a_j-\epsilon \abs{a_j}\leq \frac{Z_j}{\sigma_j}\leq f_M(b_j)+b_j+\epsilon \abs{b_j}\text{ for all }j)+2p\exp(-cM^2)\\
\leq &P(f_M(a_j)+a_j\leq \frac{Z_j}{\sigma_j}\leq f_M(b_j)+b_j\text{ for all }j)+C\epsilon M\sqrt{\log p}+2p\exp(-cM^2)\\
\leq &P(a_j\leq \frac{Z_j}{\sigma_j}\leq b_j\text{ for all }j)+4p\exp(-cM^2)+C\epsilon M\sqrt{\log p}
\end{align*}
where in the second inequality we use Lemma \ref{anti-concentration} (note that $f_M(a_j)+a_j-\epsilon \abs{a_j}$ is either $\infty$ or its absolute value $\leq \frac{3}{2}M$, so is $f_M(b_j)+b_j+\epsilon \abs{b_j}$), the term $\exp(-cM^2)$ is the tail bound of the univariate standard normal distribution, and $C,c$ are universal constants. Therefore we have derived the following upper bound
\begin{align*}
&P(\hat\sigma_ja_j\leq \sqrt n(\bar{X}_j-\mu_j)\leq\hat\sigma_j b_j\text{ for all }j)\\
\leq&P(a_j\leq \frac{Z_j}{\sigma_j}\leq b_j\text{ for all }j)+4p\exp(-cM^2)+C\epsilon M\sqrt{\log p}+C\Big(\frac{D_1^2\log^7(pn)}{n}\Big)^{\frac{1}{6}}+\sum_{j=1}^pP(\big\lvert\hat\sigma_j-\sigma_j\big\rvert> \epsilon\sigma_j)
\end{align*}
which holds true for $\epsilon<1/2$ and $M>0$. Similarly, one can show
\begin{align*}
&P(\hat\sigma_ja_j\leq \sqrt n(\bar{X}_j-\mu_j)\leq\hat\sigma_j b_j\text{ for all }j)\\
\geq &P(a_j+\epsilon \abs{a_j}\leq \frac{Z_j}{\sigma_j}\leq b_j-\epsilon \abs{b_j}\text{ for all }j)-C\Big(\frac{D_1^2\log^7(pn)}{n}\Big)^{\frac{1}{6}}-\sum_{j=1}^pP(\big\lvert\hat\sigma_j-\sigma_j\big\rvert> \epsilon\sigma_j)
\end{align*}
along with
\begin{align*}
&P(a_j+\epsilon \abs{a_j}\leq \frac{Z_j}{\sigma_j}\leq b_j-\epsilon \abs{b_j}\text{ for all }j)\\
\geq &P(f_M(a_j)+a_j+\epsilon \abs{a_j}\leq \frac{Z_j}{\sigma_j}\leq f_M(b_j)+b_j-\epsilon \abs{b_j}\text{ for all }j)-2p\exp(-cM^2)\\
\geq &P(f_M(a_j)+a_j\leq \frac{Z_j}{\sigma_j}\leq f_M(b_j)+b_j\text{ for all }j)-C\epsilon M\sqrt{\log p}-2p\exp(-cM^2)\\
\geq &P(a_j\leq \frac{Z_j}{\sigma_j}\leq b_j\text{ for all }j)-C\epsilon M\sqrt{\log p}-4p\exp(-cM^2).
\end{align*}
This leads to
\begin{align}
\nonumber&\big\lvert P(\hat\sigma_ja_j\leq \sqrt n(\bar{X}_j-\mu_j)\leq\hat\sigma_j b_j\text{ for all }j)-P(a_j\leq \frac{Z_j}{\sigma_j}\leq b_j\text{ for all }j)\big\rvert\\
\nonumber\leq&4p\exp(-cM^2)+C\epsilon M\sqrt{\log p}+C\Big(\frac{D_1^2\log^7(pn)}{n}\Big)^{\frac{1}{6}}+\sum_{j=1}^pP(\big\lvert\hat\sigma_j-\sigma_j\big\rvert> \epsilon\sigma_j)\\
\nonumber\leq&4p\exp(-cM^2)+C\epsilon M\sqrt{\log p}+C\Big(\frac{D_1^2\log^7(pn)}{n}\Big)^{\frac{1}{6}}+\sum_{j=1}^pP(\big\lvert\hat\sigma^2_j-\sigma^2_j\big\rvert>\epsilon\sigma^2_j)\\
\leq&4p\exp(-cM^2)+C\epsilon M\sqrt{\log p}+C\Big(\frac{D_1^2\log^7(pn)}{n}\Big)^{\frac{1}{6}}+4p\exp(-\frac{cn\epsilon^2}{D_1^4+D_1^2\epsilon})\label{uniform_bound}
\end{align}
where the last inequality holds because Assumption \ref{condition:clt2} guarantees that the sub-Gaussian norm of $X_j-\mu_j$ is at most $D_1\sigma_j$ and one then applies Lemma \ref{concentration:variance}. Now set
\begin{equation*}
M=\sqrt{\frac{1}{c}\log(pn)},\;\epsilon=\Big(\frac{D_1^2}{n}\Big)^{\frac{1}{6}}
\end{equation*}
and note that $\epsilon< 1/2$ can be assumed since otherwise the first term in the desired bound is already greater than $1$ (by enlarging the universal constant if necessary) and the bound is trivial. We get the uniform bound
\begin{align*}
&\big\lvert P(\hat\sigma_ja_j\leq \sqrt n(\bar{X}_j-\mu_j)\leq\hat\sigma_j b_j\text{ for all }j)-P(a_j\leq \frac{Z_j}{\sigma_j}\leq b_j\text{ for all }j)\big\rvert\\
\leq &C\Big(\frac{D_1^2\log^7(pn)}{n}\Big)^{\frac{1}{6}}+Cp\exp\Big(-\frac{cD_1^{2/3}n^{2/3}}{D_1^4+D_1^2\epsilon}\Big)\\
\leq &C\Big(\frac{D_1^2\log^7(pn)}{n}\Big)^{\frac{1}{6}}+Cp\exp\Big(-\frac{cn^{2/3}}{D_1^{10/3}}\Big)
\end{align*}
where the second inequality holds because $D_2\geq 1$ and $\epsilon< 1/2$. In particular, if $X_j$'s are $[0,1]$-valued, we use the concentration inequality \eqref{case:ccp} instead of \eqref{case:generalcvx} to refine the bound \eqref{uniform_bound} to be
\begin{align*}
&\big\lvert P(\hat\sigma_ja_j\leq \sqrt n(\bar{X}_j-\mu_j)\leq\hat\sigma_j b_j\text{ for all }j)-P(a_j\leq \frac{Z_j}{\sigma_j}\leq b_j\text{ for all }j)\big\rvert\\
\leq&4p\exp(-cM^2)+C\epsilon M\sqrt{\log p}+C\Big(\frac{D_1^2\log^7(pn)}{n}\Big)^{\frac{1}{6}}+4p\exp(-\frac{cn\delta\epsilon^2}{1+\epsilon}).
\end{align*}
Likewise, letting $M$ and $\epsilon$ take the same values as before, we obtain
\begin{align*}
&\big\lvert P(\hat\sigma_ja_j\leq \sqrt n(\bar{X}_j-\mu_j)\leq\hat\sigma_j b_j\text{ for all }j)-P(a_j\leq \frac{Z_j}{\sigma_j}\leq b_j\text{ for all }j)\big\rvert\\
\leq &C\Big(\frac{D_1^2\log^7(pn)}{n}\Big)^{\frac{1}{6}}+Cp\exp\Big(-c\delta D_1^{2/3} n^{2/3}\Big).
\end{align*}
This completes the proof of Theorems \ref{clt3:generalcvx} and \ref{clt3:ccp}.\hfill\Halmos

The following are corresponding results for the multiplier bootstrap:
\begin{theorem}\label{plugin3:general_cvx}
If Assumption \ref{condition:clt2} holds, for any constant $0<\alpha<\frac{1}{e}$ we have
\begin{align*}
&\sup_{a_j\leq b_j,j=1,\ldots,p}\Big\lvert P(\hat\sigma_ja_j\leq \widehat{Z}_j\leq \hat\sigma_jb_j\text{ for all }j\vert \{\mathbf X_i\}_{i=1}^n)-P(\sigma_ja_j\leq Z_j\leq \sigma_jb_j\text{ for all }j)\Big\rvert\\
\leq &C\Big(\frac{D_1^2\log^5(pn)\log^2(1/\alpha)}{n}\Big)^{\frac{1}{6}}+\frac{CD_1^2\log(pn)\log(p/\alpha)}{\sqrt n}
\end{align*}
with probability at least $1-\alpha$, where $C$ is a universal constant.
\end{theorem}
\begin{theorem}\label{plugin3:ccp}
If Assumptions \ref{condition:clt2} and \ref{condition:clt3_ccp} hold, for any constant $0<\alpha<\frac{1}{e}$ we have
\begin{align*}
&\sup_{a_j\leq b_j,j=1,\ldots,p}\Big\lvert P(\hat\sigma_ja_j\leq \widehat{Z}_j\leq \hat\sigma_jb_j\text{ for all }j\vert \{\mathbf X_i\}_{i=1}^n)-P(\sigma_ja_j\leq Z_j\leq \sigma_jb_j\text{ for all }j)\Big\rvert\\
\leq &C\Big(\frac{D_1^2\log^5(pn)\log^2(1/\alpha)}{n}\Big)^{\frac{1}{6}}+\frac{C\log(pn)\log(p/\alpha)}{\sqrt{n\delta}}
\end{align*}
with probability at least $1-\alpha$, where $C$ is a universal constant.
\end{theorem}
\proof{Proof of Theorems \ref{plugin3:general_cvx} and \ref{plugin3:ccp}.} For any $a_j\leq b_j,j=1,\ldots,p$ and $0<\epsilon<1/2$
\begin{align*}
&P(\hat\sigma_ja_j\leq \widehat{Z}_j\leq\hat\sigma_j b_j\text{ for all }j\vert \{\mathbf X_i\}_{i=1}^n)\\
=&P(\frac{\hat\sigma_j}{\sigma_j}a_j\leq \frac{\widehat{Z}_j}{\sigma_j}\leq\frac{\hat\sigma_j}{\sigma_j} b_j\text{ for all }j\vert \{\mathbf X_i\}_{i=1}^n)\\
=&P(\frac{\hat\sigma_j}{\sigma_j}a_j\leq \frac{\widehat{Z}_j}{\sigma_j}\leq\frac{\hat\sigma_j}{\sigma_j} b_j,\big\lvert\frac{\hat\sigma_j}{\sigma_j}-1\big\rvert\leq \epsilon\text{ for all }j\vert \{\mathbf X_i\}_{i=1}^n)\\
&+P(\frac{\hat\sigma_j}{\sigma_j}a_j\leq \frac{\widehat{Z}_j}{\sigma_j}\leq\frac{\hat\sigma_j}{\sigma_j} b_j\text{ for all }j,\big\lvert\frac{\hat\sigma_j}{\sigma_j}-1\big\rvert> \epsilon\text{ for some }j\vert \{\mathbf X_i\}_{i=1}^n)\\
\leq &P(a_j-\epsilon \abs{a_j}\leq \frac{\widehat{Z}_j}{\sigma_j}\leq b_j+\epsilon \abs{b_j}, \big\lvert\frac{\hat\sigma_j}{\sigma_j}-1\big\rvert\leq \epsilon\text{ for all }j\vert \{\mathbf X_i\}_{i=1}^n)\\
&\text{with probability at least }1-\sum_{j=1}^pP(\big\lvert\hat\sigma_j-\sigma_j\big\rvert> \epsilon\sigma_j)\\
\leq &P(a_j-\epsilon \abs{a_j}\leq \frac{\widehat{Z}_j}{\sigma_j}\leq b_j+\epsilon \abs{b_j}\text{ for all }j\vert \{\mathbf X_i\}_{i=1}^n)\\
\leq & P(a_j-\epsilon \abs{a_j}\leq \frac{Z_j}{\sigma_j}\leq b_j+\epsilon \abs{b_j}\text{ for all }j)+C\Big(\frac{D_1^2\log^5(pn)\log^2(4/\alpha)}{n}\Big)^{\frac{1}{6}}\\
&\text{with probability at least }1-\frac{\alpha}{4}.
\end{align*}
In the proof of Theorems \ref{clt3:generalcvx} and \ref{clt3:ccp} we show
\begin{align*}
P(a_j-\epsilon \abs{a_j}\leq \frac{Z_j}{\sigma_j}\leq b_j+\epsilon \abs{b_j}\text{ for all }j)
\leq &P(a_j\leq \frac{Z_j}{\sigma_j}\leq b_j\text{ for all }j)+4p\exp(-cM^2)+C\epsilon M\sqrt{\log p}.
\end{align*}
Similarly we can show the other direction
\begin{align*}
&P(\hat\sigma_ja_j\leq \widehat{Z}_j\leq\hat\sigma_j b_j\text{ for all }j\vert \{\mathbf X_i\}_{i=1}^n)\\
\geq& P(a_j+\epsilon \abs{a_j}\leq \frac{Z_j}{\sigma_j}\leq b_j-\epsilon \abs{b_j}\text{ for all }j)-C\Big(\frac{D_1^2\log^5(pn)\log^2(4/\alpha)}{n}\Big)^{\frac{1}{6}}\\
&\text{with probability at least }1-\frac{\alpha}{4}-\sum_{j=1}^pP(\big\lvert\hat\sigma_j-\sigma_j\big\rvert> \epsilon\sigma_j)\\
\geq &P(a_j\leq \frac{Z_j}{\sigma_j}\leq b_j\text{ for all }j)-C\epsilon M\sqrt{\log p}-4p\exp(-cM^2)-C\Big(\frac{D_1^2\log^5(pn)\log^2(4/\alpha)}{n}\Big)^{\frac{1}{6}}.
\end{align*}
Therefore the following uniform bound holds with probability at least $1-\frac{\alpha}{2}-2\sum_{j=1}^pP(\big\lvert\hat\sigma_j-\sigma_j\big\rvert> \epsilon\sigma_j)$
\begin{align*}
&\big\lvert P(\hat\sigma_ja_j\leq \widehat{Z}_j\leq\hat\sigma_j b_j\text{ for all }j\vert \{\mathbf X_i\}_{i=1}^n)-P(a_j\leq \frac{Z_j}{\sigma_j}\leq b_j\text{ for all }j)\big\rvert\\
\leq &C\Big(\frac{D_1^2\log^5(pn)\log^2(4/\alpha)}{n}\Big)^{\frac{1}{6}}+C\epsilon M\sqrt{\log p}+4p\exp(-cM^2).
\end{align*}
Note that
\begin{align*}
2\sum_{j=1}^pP(\big\lvert\hat\sigma_j-\sigma_j\big\rvert> \epsilon\sigma_j)\leq 8p\exp(-\frac{cn\epsilon^2}{D_1^4+D_1^2\epsilon}).
\end{align*}
By setting the right hand side of the above inequality to be $\alpha/2$ and $M=\sqrt{\frac{1}{c}\log(pn)}$ we get
\begin{align*}
&\big\lvert P(\hat\sigma_ja_j\leq \widehat{Z}_j\leq\hat\sigma_j b_j\text{ for all }j\vert \{\mathbf X_i\}_{i=1}^n)-P(a_j\leq \frac{Z_j}{\sigma_j}\leq b_j\text{ for all }j)\big\rvert\\
\leq &C\Big(\frac{D_1^2\log^5(pn)\log^2(4/\alpha)}{n}\Big)^{\frac{1}{6}}+C\Big(\sqrt{\frac{D_1^4}{n}\log\frac{16p}{\alpha}}+\frac{D_1^2}{n}\log\frac{16p}{\alpha}\Big)\log(pn)\\
\leq &C\Big(\frac{D_1^2\log^5(pn)\log^2(4/\alpha)}{n}\Big)^{\frac{1}{6}}+C\frac{D_1^2}{\sqrt{n}}\log\frac{16p}{\alpha}\log(pn)
\end{align*}
with probability at least $1-\alpha$.

In case of $[0,1]$-valued variables, we use \eqref{case:ccp} instead of \eqref{case:generalcvx} to get
\begin{align*}
&\big\lvert P(\hat\sigma_ja_j\leq \widehat{Z}_j\leq\hat\sigma_j b_j\text{ for all }j\vert \{\mathbf X_i\}_{i=1}^n)-P(a_j\leq \frac{Z_j}{\sigma_j}\leq b_j\text{ for all }j)\big\rvert\\
\leq &C\Big(\frac{D_1^2\log^5(pn)\log^2(4/\alpha)}{n}\Big)^{\frac{1}{6}}+C\Big(\sqrt{\frac{1}{n\delta}\log\frac{16p}{\alpha}}+\frac{1}{n\delta}\log\frac{16p}{\alpha}\Big)\log(pn)\\
\leq &C\Big(\frac{D_1^2\log^5(pn)\log^2(4/\alpha)}{n}\Big)^{\frac{1}{6}}+C\Big(\frac{1}{\sqrt{n\delta}}+\frac{1}{n\delta}\Big)\log\frac{16p}{\alpha}\log(pn)
\end{align*}
with probability at least $1-\alpha$. We can assume that $n\delta\geq 1$ to get the desired bound since otherwise the second term in the bound is already greater than $1$.\hfill\Halmos

\subsection{Coverage Probability through Multiplier Bootstrap}\label{multiplier bootstrap:unnormalized and normalized}
\begin{theorem}[Coverage probability for unnormalized supremum]\label{consistency:unnormalized}
Under Assumption \ref{condition:clt2}, for every $0<\beta<1$ we have
\begin{align*}
\lvert P(\sqrt n(\bar{X}_j-\mu_j)\leq \hat q_{1-\beta}\text{ for all }j)- (1-\beta)\rvert\leq C\Big(\frac{D_1^2\log^7(pn)}{n}\Big)^{\frac{1}{6}}
\end{align*}
where $\hat q_{1-\beta}$ is such that
\begin{align*}
P(\widehat Z_j\leq \hat q_{1-\beta}\text{ for all }j\vert \{\mathbf X_i\}_{i=1}^n)=1-\beta
\end{align*}
and $C$ is a universal constant.
\end{theorem}
\proof{Proof of Theorem \ref{consistency:unnormalized}.}Denote by $\epsilon=C\Big(\frac{D_1^2\log^5(pn)\log^2(1/\alpha)}{n}\Big)^{1/6}$ and by $A_{\alpha}$ the event that
\begin{align*}
\sup_{a_j\leq b_j,j=1,\ldots,p}\Big\lvert P(a_j\leq \widehat{Z}_j\leq b_j\text{ for all }j\vert \{\mathbf X_i\}_{i=1}^n)-P(a_j\leq Z_j\leq b_j\text{ for all }j)\Big\rvert\leq \epsilon.
\end{align*}
Then we can rewrite
\begin{align*}
& P(\sqrt n(\bar{X}_j-\mu_j)\leq \hat q_{1-\beta}\text{ for all }j)\\
=&P(\sqrt n(\bar{X}_j-\mu_j)\leq \hat q_{1-\beta}\text{ for all }j,\text{ and }A_{\alpha})+P(\sqrt n(\bar{X}_j-\mu_j)\leq \hat q_{1-\beta}\text{ for all }j,\text{ and }A^c_{\alpha})
\end{align*}
The second term is bounded by $\alpha$ because of Theorem \ref{plugin2}. To study the first term, denote by $q_{1-\beta}$ the true ($1-\beta$)-level quantile of the limit distribution, i.e., $q_{1-\beta}$ is such that $P(Z_j\leq  q_{1-\beta}\text{ for all }j)=1-\beta$. On event $A_{\alpha}$ we have $q_{1-\beta-\epsilon}\leq\hat q_{1-\beta} \leq q_{1-\beta+\epsilon}$, therefore
\begin{eqnarray*}
P(\sqrt n(\bar{X}_j-\mu_j)\leq q_{1-\beta-\epsilon}\text{ for all }j,\text{ and }A_{\alpha})&\leq& P(\sqrt n(\bar{X}_j-\mu_j)\leq \hat q_{1-\beta}\text{ for all }j,\text{ and }A_{\alpha})\\
&\leq&P(\sqrt n(\bar{X}_j-\mu_j)\leq q_{1-\beta+\epsilon}\text{ for all }j,\text{ and }A_{\alpha}).
\end{eqnarray*}
From this two-sided bound we get
\begin{align*}
& P(\sqrt n(\bar{X}_j-\mu_j)\leq \hat q_{1-\beta}\text{ for all }j)\\
\leq &P(\sqrt n(\bar{X}_j-\mu_j)\leq q_{1-\beta+\epsilon}\text{ for all }j)+\alpha\\
\leq &P(Z_j\leq q_{1-\beta+\epsilon}\text{ for all }j)+C\Big(\frac{D_1^2\log^7(pn)}{n}\Big)^{1/6}+\alpha\\
=&1-\beta+\epsilon+C\Big(\frac{D_1^2\log^7(pn)}{n}\Big)^{1/6}+\alpha.
\end{align*}
Similarly the lower bound can be derived as
\begin{align*}
& P(\sqrt n(\bar{X}_j-\mu_j)\leq \hat q_{1-\beta}\text{ for all }j)\\
\geq&P(\sqrt n(\bar{X}_j-\mu_j)\leq q_{1-\beta-\epsilon}\text{ for all }j,\text{ and }A_{\alpha})\\
= &P(\sqrt n(\bar{X}_j-\mu_j)\leq q_{1-\beta+\epsilon}\text{ for all }j)-P(\sqrt n(\bar{X}_j-\mu_j)\leq q_{1-\beta-\epsilon}\text{ for all }j,\text{ and }A^c_{\alpha})\\
\geq &P(Z_j\leq q_{1-\beta-\epsilon}\text{ for all }j)-C\Big(\frac{D_1^2\log^7(pn)}{n}\Big)^{1/6}-\alpha\\
=&1-\beta-\epsilon-C\Big(\frac{D_1^2\log^7(pn)}{n}\Big)^{1/6}-\alpha.
\end{align*}
This gives the following bound for any $\alpha<1/e$
\begin{align*}
\lvert P(\sqrt n(\bar{X}_j-\mu_j)\leq \hat q_{1-\beta}\text{ for all }j)- (1-\beta)\rvert\leq C\Big(\frac{D_1^2\log^7(pn)}{n}\Big)^{\frac{1}{6}}+\epsilon+\alpha.
\end{align*}
Set $\alpha=1/n$ and note that $1/n$ is less than the leading term, thus we have shown the desired conclusion.\hfill\Halmos

\begin{theorem}[Coverage probability for normalized supremum]\label{consistency:normalized_generalcvx}
Under Assumptions \ref{condition:clt2}, for every $0<\beta<1$ we have
\begin{align*}
&\lvert P(\sqrt n(\bar{X}_j-\mu_j)\leq \hat\sigma_j\hat q_{1-\beta}\text{ for all }j)-(1-\beta)\rvert\\
\leq  &C\Big(\Big(\frac{D_1^2\log^7(pn)}{n}\Big)^{\frac{1}{6}}+\frac{D_1^2\log^2(pn)}{\sqrt n}+p\exp\big(-\frac{cn^{2/3}}{D_1^{10/3}}\big)\Big).
\end{align*}
If Assumption \ref{condition:clt3_ccp} also holds, then
\begin{align*}
&\lvert P(\sqrt n(\bar{X}_j-\mu_j)\leq \hat\sigma_j\hat q_{1-\beta}\text{ for all }j)-(1-\beta)\rvert\\
\leq& C\Big(\Big(\frac{D_1^2\log^7(pn)}{n}\Big)^{\frac{1}{6}}+\frac{\log^2 (pn)}{\sqrt{n\delta}}+p\exp\big(-c\delta D_1^{2/3}n^{2/3}\big)\Big).
\end{align*}
Here $\hat q_{1-\beta}$ is such that
\begin{align*}
P(\widehat Z_j\leq \hat\sigma_j\hat q_{1-\beta}\text{ for all }j\vert \{\mathbf X_i\}_{i=1}^n)=1-\beta
\end{align*}
and $C,c$ are universal constants.
\end{theorem}
\proof{Proof of Theorem \ref{consistency:normalized_generalcvx}.}Let $\epsilon=C\Big(\frac{D_1^2\log^5(pn)\log^2(1/\alpha)}{n}\Big)^{1/6}+\frac{CD_1^2\log(pn)\log(p/\alpha)}{\sqrt n}$ and $A_{\alpha}$ be the event that
\begin{equation*}
\sup_{a_j\leq b_j,j=1,\ldots,p}\Big\lvert P(\hat\sigma_ja_j\leq \widehat{Z}_j\leq \hat\sigma_jb_j\text{ for all }j\vert \{\mathbf X_i\}_{i=1}^n)-P(\sigma_ja_j\leq Z_j\leq \sigma_jb_j\text{ for all }j)\Big\rvert\leq \epsilon.
\end{equation*}
We know that $P(A_{\alpha})\leq \alpha$ from Theorem \ref{plugin3:general_cvx}. Following the same line of the proof for Theorem \ref{consistency:unnormalized} and using the CLT in Theorem \ref{clt3:generalcvx} we can derive that
\begin{align*}
\lvert P(\sqrt n(\bar{X}_j-\mu_j)\leq \hat\sigma_j\hat q_{1-\beta}\text{ for all }j)-(1-\beta)\rvert\leq C\Big(\frac{D_1^2\log^7(pn)}{n}\Big)^{\frac{1}{6}}+Cp\exp\big(-\frac{cn^{2/3}}{D_1^{10/3}}\big)+\epsilon+\alpha.
\end{align*}
Again setting $\alpha=\frac{1}{n}$ leads to the first bound.

The second bound can be derived similarly. Define $\epsilon=C\Big(\frac{D_1^2\log^5(pn)\log^2(1/\alpha)}{n}\Big)^{1/6}+\frac{C\log(pn)\log(p/\alpha)}{\sqrt{n\delta}}$, and now Theorem \ref{plugin3:ccp} entails that $P(A_{\alpha})\leq \alpha$ again. Using the CLT in Theorem \ref{clt3:ccp} gives
\begin{align*}
\lvert P(\sqrt n(\bar{X}_j-\mu_j)\leq \hat\sigma_j\hat q_{1-\beta}\text{ for all }j)-(1-\beta)\rvert\leq C\Big(\frac{D_1^2\log^7(pn)}{n}\Big)^{\frac{1}{6}}+Cp\exp\big(-c\delta D_1^{2/3}n^{2/3}\big)+\epsilon+\alpha.
\end{align*}
The second bound follows from setting $\alpha=\frac{1}{n}$.\hfill\Halmos


\subsection{Proofs of Main Statistical Guarantees}\label{main proof: gaussian supremum validators}
We now put together all the previous results to prove the statistical guarantees of our validators. For convenience, we suppress the subscript $\bm{\xi}_{1:n_2}$ in the probability notation.
\proof{Proof of Theorem \ref{feasibility_supremum_unnormalized:general constraint}.}We bound the probability as follows
\begin{eqnarray}
\notag&&P(x^*(\hat s^*)\text{ is feasible for \eqref{stoc_opt}})\\
\notag&\geq& P(\hat H_j\geq \gamma + \frac{ q_{1-\beta}}{\sqrt{n_2}}\text{ for some }j=1,\ldots,p\text{ in \eqref{validate:gs unnormalized} and }H(x^*(s_j))\geq\hat H_j-  \frac{ q_{1-\beta}}{\sqrt{n_2}}\text{ for all } j=1,\ldots,p)\\
\notag&\geq&P(H(x^*(s_j))\geq\hat H_j-  \frac{ q_{1-\beta}}{\sqrt{n_2}}\text{ for all } j=1,\ldots,p)-P(\hat H_j< \gamma + \frac{ q_{1-\beta}}{\sqrt{n_2}}\text{ for all }j=1,\ldots,p)\\
&\geq&1-\beta-C\Big(\frac{D_1^2\log^7(pn_2)}{n_2}\Big)^{\frac{1}{6}}-P(\hat H_j< \gamma + \frac{ q_{1-\beta}}{\sqrt{n_2}}\text{ for all }j=1,\ldots,p)\label{bound1 for feasibility}
\end{eqnarray}
where we use Theorem \ref{consistency:unnormalized} for the first probability by letting $X_{i,j}=h(x^*(s_j),\xi_i)$. To bound the second probability, we recall that $\hat\sigma_j^2$ is the sample variance computed from $\{h(x^*(s_j),\xi_1),\ldots,h(x^*(s_j),\xi_{n_2})\}$ and write
\begin{eqnarray*}
&&P(\hat H_j< \gamma + \frac{ q_{1-\beta}}{\sqrt{n_2}}\text{ for all }j=1,\ldots,p)\\
&\leq &P(\hat H_j< \gamma + \frac{ q_{1-\beta}}{\sqrt{n_2}}\text{ and }\hat{\sigma}^2_j\leq 2\sigma^2(x^*(s_j))\text{ for all }j=1,\ldots,p)+P(\hat{\sigma}^2_j> 2\sigma^2(x^*(s_j))\text{ for some }j=1,\ldots,p)\\
&\leq &P(\hat H_j< \gamma + \frac{C\bar{\sigma}\sqrt{\log(p/\beta)}}{\sqrt{n_2}}\text{ for all }j=1,\ldots,p)+P(\hat{\sigma}^2_j> 2\sigma^2(x^*(s_j))\text{ for some }j=1,\ldots,p)\\
&&\text{\ \ because of the fact that $q_{1-\beta}\leq C\max_{j}\hat{\sigma}_j\sqrt{\log(p/\beta)}$ for some universal constant $C$}\\
&\leq &P(\hat H_{\bar{j}}-\overline{H}< \gamma + \frac{C\bar{\sigma}\sqrt{\log(p/\beta)}}{\sqrt{n_2}}-\overline{H})+Cp\exp(-\frac{cn_2}{D_1^4})\\
&&\text{\ \ where $\bar{j}$ is the index such that $H(x^*(s_{\bar{j}}))=\overline{H}$ and the concentration \eqref{case:generalcvx} is used}\\
&\leq&C\exp\big(-\frac{cn_2\epsilon^2}{D_1^2\bar{\sigma}^2}\big)+Cp\exp(-\frac{cn_2}{D_1^4})\text{\ \ because the sub-Gaussian norm of }h(x^*(s_{\bar{j}}),\xi)\text{ is at most }D_1\bar{\sigma}.
\end{eqnarray*}
Substituting this bound into \eqref{bound1 for feasibility} gives the desired conclusion.\hfill\Halmos
\proof{Proof of Theorem \ref{feasibility_supremum_normalized:general constraint}.}Similar to the proof of Theorem \ref{feasibility_supremum_unnormalized:general constraint}, we have the bound
\begin{eqnarray*}
&&P(x^*(\hat s^*)\text{ is feasible for \eqref{stoc_opt}})\\
&\geq&P(H(x^*(s_j))\geq\hat H_j-  \frac{ q_{1-\beta}\hat\sigma_j}{\sqrt{n_2}}\text{ for all } j=1,\ldots,p)-P(\hat H_j< \gamma + \frac{ q_{1-\beta}\hat\sigma_j}{\sqrt{n_2}}\text{ for all }j=1,\ldots,p)\\
&\geq&1-\beta-C\prth{\Big(\frac{D_1^2\log^7(pn_2)}{n_2}\Big)^{\frac{1}{6}}+\frac{D_1^2\log^2(pn_2)}{\sqrt{n_2}}+p\exp\big(-\frac{cn_2^{2/3}}{D_1^{10/3}}\big)}-P(\hat H_j< \gamma + \frac{ q_{1-\beta}\hat\sigma_j}{\sqrt{n_2}}\text{ for all }j=1,\ldots,p)\\
&&\text{\ \ where the first bound is due to Theorem \ref{consistency:normalized_generalcvx}}.
\end{eqnarray*}
For the second probability we write
\begin{eqnarray*}
&&P(\hat H_j< \gamma + \frac{ q_{1-\beta}\hat\sigma_j}{\sqrt{n_2}}\text{ for all }j=1,\ldots,p)\\
&\leq&P(\hat H(x^*(\bar{s}))< \gamma + \frac{ q_{1-\beta}\hat\sigma(x^*(\bar{s}))}{\sqrt{n_2}})\\
&\leq &P(\hat H(x^*(\bar{s}))< \gamma + \frac{ q_{1-\beta}\sqrt{2}\sigma(x^*(\bar{s}))}{\sqrt{n_2}})+P(\hat\sigma^2(x^*(\bar{s}))>2\sigma^2(x^*(\bar{s})))\\
&\leq &P(\hat H(x^*(\bar{s}))< \gamma + \frac{C\sqrt{\log(p/\beta)}\sigma(x^*(\bar{s}))}{\sqrt{n_2}})+P(\hat\sigma^2(x^*(\bar{s}))>2\sigma^2(x^*(\bar{s})))\text{\ \ since }q_{1-\beta}\leq C\sqrt{\log(p/\beta)}\\
&\leq&C\exp\big(-\frac{cn_2\epsilon^2}{D_1^2\sigma^2(x^*(\bar{s}))}\big)+C\exp(-\frac{cn_2}{D_1^4})
\end{eqnarray*}
Combining the two probability bounds and noting that $p\exp\big(-\frac{cn_2^{2/3}}{D_1^{10/3}}\big)$ dominates $\exp(-\frac{cn_2}{D_1^4})$ (because $\frac{n_2}{D_1^4}=\frac{n_2^{2/3}}{D_1^{10/3}}\cdot \big(\frac{n_2}{D_1^2}\big)^{1/3}$ and $\frac{n_2}{D_1^2}\geq 1$ can be assumed), we obtain the desired conclusion. \hfill\Halmos

\proof{Proof of Corollary \ref{asymptotic feasibility:supremum + general constraint}.}When $\overline{H}>\gamma$, we have $\epsilon\to \overline{H}-\gamma>0$ in Theorems \ref{feasibility_supremum_unnormalized:general constraint} and \ref{feasibility_supremum_normalized:general constraint}, therefore the exponential error term with $\epsilon$ vanishes as $n_2\to\infty$. Under the condition that $p\exp(-n_2^{1/7})\to 0$ it is straightforward to check that other error terms also vanish.\hfill\Halmos

\proof{Proof of Theorem \ref{feasibility:gs unnormalized}.}Unlike the proof of Theorem \ref{feasibility_supremum_unnormalized:general constraint}, we use the Bernoulli structure to derive the error bound. Note that in this case $\gamma=1-\alpha$. Define events
\begin{align*}
E_1&=\big\{\hat H_j\geq 1-\alpha + \frac{ q_{1-\beta}}{\sqrt{n_2}}\text{ for some }j=1,\ldots,p\text{ in \eqref{validate:gs unnormalized}}\big\}\\
E_2&=\big\{H(x^*(s_j))\geq\hat H_j-  \frac{ q_{1-\beta}}{\sqrt{n_2}}\text{ for all $j$ such that }H(x^*(s_j))\in(\alpha,1-\alpha)\big\}\\
E_3&=\big\{\hat H_j<1-\alpha + \frac{ q_{1-\beta}}{\sqrt{n_2}}\text{ for all $j$ such that }H(x^*(s_j))\leq \alpha\big\}.
\end{align*}
Then we have
\begin{eqnarray}
\notag P(x^*(\hat s^*)\text{ is feasible for \eqref{stoc_opt}})&\geq& P(E_1\cap E_2\cap E_3)\\
\notag&\geq&1-P(E_1^c)-P(E_2^c)-P(E_3^c)\\
&=&P(E_2)-P(E_1^c)-P(E_3^c).\label{three probability bounds}
\end{eqnarray}
We bound the three probabilities. Let $q^{\alpha}_{1-\beta}$ be the $1-\beta$ quantile of $\max\{Z_j:H(x^*(s_j))\in(\alpha,1-\alpha),1\leq j\leq p\}$ where $(Z_1,\ldots,Z_p)\sim N_p(0,\hat\Sigma)$. By stochastic dominance it is clear that $q^{\alpha}_{1-\beta}\leq q_{1-\beta}$ almost surely, therefore
\begin{eqnarray*}
P(E_2)\geq P\big(H(x^*(s_j))&\geq&\hat H_j-  \frac{ q^{\alpha}_{1-\beta}}{\sqrt{n_2}}\text{ for all $j$ such that }H(x^*(s_j))\in(\alpha,1-\alpha)\big)\\
&\geq&1-\beta-C\Big(\frac{\log^7(pn_2)}{n_2\alpha}\Big)^{\frac{1}{6}}
\end{eqnarray*}
by applying Theorem \ref{consistency:unnormalized} to $\{h(x^*(s_j),\xi):H(x^*(s_j))\in(\alpha,1-\alpha),1\leq j\leq p\}$ and noticing that Assumption \ref{condition:clt2} is satisfied with $D_1=\frac{C}{\sqrt \alpha}$ for some universal constant $C$.

We then bound the second probability
\begin{eqnarray*}
P(E_1^c)&=&P(\hat H_j< 1-\alpha + \frac{ q_{1-\beta}}{\sqrt{n_2}}\text{ for all }j=1,\ldots,p)\\
&\leq&P(\hat H_{\bar{j}}< 1-\alpha + \frac{ q_{1-\beta}}{\sqrt{n_2}})\text{\ \ where $\bar{j}$ is the index such that }H(x^*(s_{\bar{j}}))=1-\bar{\alpha}\\
&\leq &P(\hat H_{\bar{j}}< 1-\alpha + \frac{C\sqrt{\log(p/\beta)}}{\sqrt{n_2}})\text{\ \ because }q_{1-\beta}\leq C\max_{j}\hat{\sigma}_j\sqrt{\log(p/\beta)}\leq C\sqrt{\log(p/\beta)}\\
&\leq&\exp\big(-\frac{n_2\epsilon^2}{2(\bar{\alpha}(1-\bar{\alpha})+\epsilon/3)}\big)
\end{eqnarray*}
where in the last line we use a Bernstein's inequality for sums of bounded random variables (see equation (2.10) in \cite{boucheron2013concentration}). Note that this is further bounded by $\exp\big(-cn_2\min\{\epsilon,\frac{\epsilon^2}{\bar{\alpha}}\}\big)$ if $\bar{\alpha}\leq 1/2$.

The third probability can be bounded as
\begin{eqnarray*}
P(E_3^c)&\leq&P\big(\hat H_j\geq 1-\alpha \text{ for some $j$ such that }H(x^*(s_j))\leq \alpha\big)\\
&\leq&\sum_{j:H(x^*(s_j))\leq \alpha}P(\hat H_j\geq 1-\alpha)\\
&\leq&p\exp(-2n_2(1-2\alpha)^2)\leq p\exp(-cn_2)\text{\ \ by Hoeffding's inequality}.
\end{eqnarray*}
Substituting the bounds into \eqref{three probability bounds} leads to
\begin{equation*}
P(x^*(\hat s^*)\text{ is feasible for \eqref{stoc_opt}})\leq 1-\beta-C\Big(\frac{\log^7(pn_2)}{n_2\alpha}\Big)^{\frac{1}{6}}-\exp\big(-cn_2\min\{\epsilon,\frac{\epsilon^2}{\bar{\alpha}}\}\big)-p\exp(-cn_2).
\end{equation*}

It remains to show that $p\exp(-cn_2)$ is negligible relative to other error terms. Since $\alpha<1$ it is clear that $\big(\frac{1}{n_2}\big)^{1/6}\leq \Big(\frac{\log^7(pn_2)}{n_2\alpha}\Big)^{1/6}$, and we argue that $\big(\frac{1}{n_2}\big)^{1/6}\geq p\exp(-cn_2)$ can be assumed so that $p\exp(-cn_2)\leq \Big(\frac{\log^7(pn_2)}{n_2\alpha}\Big)^{1/6}$. If $\big(\frac{1}{n_2}\big)^{1/6}< p\exp(-cn_2)$, then $p>\exp(cn_2)n_2^{-1/6}$, and $\frac{\log^7(pn_2)}{n_2\alpha}\geq \frac{(cn_2)^7}{n_2\alpha}\geq c^7n_2^6$, hence the first error term already exceeds $1$ (enlarge the universal constant $C$ if necessary) and the error bound holds true trivially.\hfill\Halmos

\proof{Proof of Theorem \ref{feasibility:gs normalized}.}The proof follows the one for Theorem \ref{feasibility:gs unnormalized}, and we focus on the modifications. The events are now defined as
\begin{align*}
E_1&=\big\{\hat H_j\geq 1-\alpha + \frac{ q_{1-\beta}\hat\sigma_j}{\sqrt{n_2}}\text{ for some }j=1,\ldots,p\text{ in \eqref{validate:gs normalized}}\big\}\\
E_2&=\big\{H(x^*(s_j))\geq\hat H_j-  \frac{ q_{1-\beta}\hat\sigma_j}{\sqrt{n_2}}\text{ for all $j$ such that }H(x^*(s_j))\in(\alpha,1-\alpha)\big\}\\
E_3&=\big\{\hat H_j<1-\alpha + \frac{ q_{1-\beta}\hat\sigma_j}{\sqrt{n_2}}\text{ for all $j$ such that }H(x^*(s_j))\leq \alpha\big\}.
\end{align*}
Again we have $P(x^*(\hat s^*)\text{ is feasible for \eqref{stoc_opt}})\geq P(E_2)-P(E_1^c)-P(E_3^c)$.

The first probability bound becomes
\begin{eqnarray*}
P(E_2)\geq 1-\beta-C\Big(\Big(\frac{\log^7(pn_2)}{n_2\alpha}\Big)^{\frac{1}{6}}+\frac{\log^2 (pn_2)}{\sqrt{n_2\alpha}}+p\exp\big(-c(n_2\alpha)^{2/3}\big)\Big)
\end{eqnarray*}
by using the second half of Theorem \ref{consistency:normalized_generalcvx} and noting that $\delta=\alpha(1-\alpha)\geq \frac{1}{2}\alpha$ if $\alpha<\frac{1}{2}$ and $D_1=\frac{C}{\sqrt \alpha}$. For the second probability we have
\begin{eqnarray*}
P(E_1^c)&\leq&P(\hat H_{\bar{j}}< 1-\alpha + \frac{ q_{1-\beta}\hat\sigma_{\bar{j}}}{\sqrt{n_2}})\text{\ \ where $\bar{j}$ is the index such that }H(x^*(s_{\bar{j}}))=1-\bar{\alpha}\\
&\leq&P(\hat H_{\bar{j}}< 1-\alpha + \frac{ q_{1-\beta}t}{\sqrt{n_2}})+P(\hat\sigma_{\bar{j}}>t)\text{\ \  where }t=\sqrt{\bar{\alpha}(1-\bar{\alpha})}+\sqrt{2\log(n_2\alpha)/n_2}\\
&\leq &P(\hat H_{\bar{j}}< 1-\alpha + \frac{ q_{1-\beta}t}{\sqrt{n_2}})+\frac{1}{n_2\alpha}\\
&&\text{\ where the bound $\frac{1}{n_2\alpha}$ is derived from \eqref{case:ccp} (see Theorem 10 in \cite{maurer2009empirical})}\\
&\leq &P(\hat H_{\bar{j}}< 1-\alpha + \frac{C\sqrt{(\bar{\alpha}+\log(n_2\alpha)/n_2)\log(p/\beta)}}{\sqrt{n_2}})+\frac{1}{n_2\alpha}\text{\ \ because }q_{1-\beta}\leq C\sqrt{\log(p/\beta)}\\
&\leq&\exp\big(-\frac{n_2\epsilon^2}{2(\bar{\alpha}(1-\bar{\alpha})+\epsilon/3)}\big)+\frac{1}{n_2\alpha}\text{\ \ by Bernstein's inequality}.
\end{eqnarray*}
Whereas for the third probability we still have $P(E_3^c)\leq p\exp(-cn_2)$.

Finally, using a similar argument in the proof of Theorem \ref{feasibility:gs unnormalized}, we can show that $\frac{1}{n_2\alpha}$, $p\exp(-cn_2)$, and $p\exp\big(-c(n_2\alpha)^{2/3}\big)$ are all dominated by $\Big(\frac{\log^7(pn_2)}{n_2\alpha}\Big)^{1/6}$ when $\Big(\frac{\log^7(pn_2)}{n_2\alpha}\Big)^{1/6}<1$, therefore the desired conclusion follows from combining the three probability bounds.\hfill\Halmos

\proof{Proof of Corollary \ref{asymptotic feasibility:supremum + chance constraint}.}Like Corollary \ref{asymptotic feasibility:supremum + general constraint}, this is a direct consequence of the finite sample result, Theorem \ref{feasibility:gs unnormalized} or \ref{feasibility:gs normalized}.\hfill\Halmos

\section{Proofs of Results in Section \ref{sec:standard gaussian margin}}\label{sec:proofs univariate}
\proof{Proof of Proposition \ref{two valid classes of constraints}.}Case (i): Assumption \ref{Donsker h} follows from the Jain-Marcus theorem (see Example 2.11.13 in \cite{van1996weak}). Assumption \ref{moment2} holds because $\sup_{x\in\mathcal X}\lvert h(x,\xi)\rvert\leq \lvert h(\tilde x,\xi)\rvert+\mathrm{diam}(\mathcal X)M(\xi)$, where $\mathrm{diam}(X)$ denotes the (finite) diameter of $\mathcal X$, and finiteness of second moments of $h(\tilde x,\xi)$ and $M(\xi)$. Assumption \ref{continuous constraint} then follows from the Lipschitz continuity of $h(x,\xi)$ in $x$ and an application of the dominated convergence theorem.

Case (ii): We need two results from empirical process theory to verify Assumption \ref{Donsker h}:
\begin{theorem}[Theorem 2.6.8 and its proof in \cite{van1996weak}]\label{VC implies Donsker}
If a class $\mathcal F$ of measurable functions satisfies:
\begin{enumerate}
\item[i.]there exists a countable subset $\mathcal F_c\subseteq \mathcal F$ such that for every $f\in\mathcal F$ there exists a sequence $f_n\in \mathcal F_c$ such that $\lim_{n\to\infty}f_n(\xi)=f(\xi)$ for all $\xi$;
\item[ii.]the envelope $\mathcal E(\xi):=\sup_{f\in \mathcal F}\lvert f(\xi)\rvert$ satisfies $\mathbb E_F[(\mathcal E(\xi))^2]<\infty$;
\item[iii.]$\mathcal F$ is a VC-subgraph class (see Section 2.6.2 of \cite{van1996weak}),
\end{enumerate}
then $\mathcal F$ is $F$-Donsker.
\end{theorem}
\begin{theorem}[Theorem 2.10.1 in \cite{van1996weak}]\label{Donsker permanance:subclass}
If a function class $\mathcal F$ is $F$-Donsker, then any subclass $\mathcal G\subseteq \mathcal F$ is also $F$-Donsker.
\end{theorem}
\begin{theorem}[Example 2.10.8 in \cite{van1996weak}]\label{Donsker permanance:product}
If $\mathcal F$ and $\mathcal G$ both are uniformly bounded $F$-Donsker classes, then $\mathcal F\cdot \mathcal G:=\{fg:f\in\mathcal F,g\in\mathcal G\}$ is also $F$-Donsker.
\end{theorem}
In order to show $F$-Donskerness of the class of constraint functions, it suffices to show $F$-Donskerness for the larger function class $\tilde{\mathcal F}:=\{\mathbf{1}(a'_kx_k\leq b_ky_k+z_k\text{ for }k=1,\ldots,K):x_k\in\R^{m_k},y_k,z_k\in\R,\text{ for }k=1,\ldots,K\}$ according to Theorem \ref{Donsker permanance:subclass}. Moreover, note that $\tilde{\mathcal F}=\tilde{\mathcal F}_1\cdot\tilde{\mathcal F}_2 \cdots \tilde{\mathcal F}_K$ where each $\tilde{\mathcal F}_k:=\{\mathbf{1}(a'_kx\leq b_ky+z):x\in\R^{m_k},y,z\in\R\}$, therefore by applying Theorem \ref{Donsker permanance:product} recursively we see that $F$-Donskerness for all $\tilde{\mathcal F}_k$'s implies $F$-Donskerness of $\tilde{\mathcal F}$. It remains to prove $F$-Donskerness of each $\tilde{\mathcal F}_k$ using Theorem \ref{VC implies Donsker}. Among the conditions of Theorem \ref{VC implies Donsker}, (ii) is trivially satisfied since the family of indicator functions is uniformly bounded by $1$. By writing $a'_kx\leq b_ky+z$ as $(a'_k, -b_k,-1) (x',y,z)' \leq 0$ and noting that the collection of all half-spaces on $\R^{m_k+2}$ has a $VC$ dimension $m_k+4$ (Problem 14 in Section 2.6 in \cite{van1996weak}), we have that $\tilde{\mathcal F}_k$ is a VC-subgraph class (Problem 9 in Section 2.6 in \cite{van1996weak}) therefore (iii) holds. To verify condition (i), consider the countable subclass $\tilde{\mathcal F}_k^c=\{\mathbf{1}(a'_kx\leq b_ky+z):x\in\mathbb Q^{m_k},y,z\in\mathbb Q\}$ where $\mathbb Q$ denotes the set of all rationals. Given $x_o\in \R^{m_k},y_o,z_o\in \R$, one can pick a sequence $x^i_o\in\mathbb Q^{m_k},y^i_o,z^i_o\in\mathbb Q$ such that $z^i_o>z_o$ for all $i$, $\lim_{i\to\infty }x^i_o\to x_o,\lim_{i\to\infty }y^i_o\to y_o,\lim_{i\to\infty }z^i_o\to z_o$ and
\begin{equation}\label{choice of approximation sequence}
\lim_{i\to\infty}\frac{\Vert x^i_o-x_o\Vert_2+\lvert y^i_o-y_o\rvert}{z^i_o-z_o}=0.
\end{equation}
For every fixed $a_k,b_k$ such that $a'_kx_o< b_ky_o+z_o$ ($a'_kx_o> b_ky_o+z_o$) we have $a'_kx^i_o< b_ky^i_o+z^i_o$ ($a'_kx^i_o> b_ky^i_o+z^i_o$) for sufficiently large $i$ because of the convergence of $x^i_o,y^i_o,z^i_o$ to $x_o,y_o,z_o$. For $a_k,b_k$ such that $a'_kx_o= b_ky_o+z_o$ we have $a'_kx^i_o\leq b_ky^i_o+z^i_o$ for sufficiently large $i$ thanks to \eqref{choice of approximation sequence}. Therefore $\mathbf{1}(a'_kx^i_o\leq b_ky^i_o+z^i_o)$ converges to $\mathbf{1}(a'_kx_o\leq b_ky_o+z_o)$ pointwise as $i\to\infty$, giving rise to condition (i). Theorem \ref{VC implies Donsker} then implies that each $\tilde{\mathcal F}_k$ is $F$-Donsker.

Assumption \ref{moment2} trivially holds since indicator functions are uniformly bounded by $1$. It remains to prove Assumption \ref{continuous constraint}. For any $x,x'$ we write
\begin{eqnarray*}
&&\lvert \mathbf{1}(a'_kA_k(x)\leq b_k\text{ for }k=1,\ldots,K)-\mathbf{1}(a'_kA_k(x')\leq b_k\text{ for }k=1,\ldots,K)\rvert\\
&\leq&\sum_{k=1}^K\mathbf{1}(a'_kA_k(x)\leq b_k < a'_kA_k(x')\text{ or }a'_kA_k(x')\leq b_k < a'_kA_k(x))\\
&\leq&\sum_{k=1}^K\mathbf{1}(\lvert a'_kA_k(x)- b_k\rvert \leq \Vert a_k\Vert_2\Vert A_k(x)-A_k(x')\Vert_2).
\end{eqnarray*}
Therefore
\begin{eqnarray*}
&&\mathbb E_F[\lvert \mathbf{1}(a'_kA_k(x)\leq b_k\text{ for }k=1,\ldots,K)-\mathbf{1}(a'_kA_k(x')\leq b_k\text{ for }k=1,\ldots,K)\rvert^2]\\
&\leq&\big(\sum_{k=1}^K\mathbb P_F(\lvert a'_kA_k(x)- b_k\rvert \leq \Vert a_k\Vert_2\Vert A_k(x)-A_k(x')\Vert_2)\big)^2
\end{eqnarray*}
hence it suffices to show each $\mathbb P_F(\lvert a'_kA_k(x)- b_k\rvert \leq \Vert a_k\Vert_2\Vert A_k(x)-A_k(x')\Vert_2)\to 0$ as $x'\to x$. We use the bound $\mathbb P_F(\lvert a'_kA_k(x)- b_k\rvert \leq \Vert a_k\Vert_2\Vert A_k(x)-A_k(x')\Vert_2)\leq \mathbb P_F(\lvert a'_kA_k(x)- b_k\rvert \leq \epsilon)+\mathbb P_F( \Vert a_k\Vert_2\Vert A_k(x)-A_k(x')\Vert_2>\epsilon)$ for any $\epsilon>0$. On one hand we have $\mathbb P_F(\lvert a'_kA_k(x)- b_k\rvert \leq \epsilon)\to 0$ as $\epsilon\to 0$. To explain, if $a_k$ has a density and $b_k\neq 0$, then $a'_kA_k(x)- b_k$ either has a density on $\R$ or is a point mass at $b_k$ (when $A_k(x)$ is the zero vector), either of which implies $\mathbb P_F(\lvert a'_kA_k(x)- b_k\rvert \leq \epsilon)\to 0$. Otherwise if $(a_k,b_k)$ has a joint density, $a'_kA_k(x)- b_k$ has a density hence $\mathbb P_F(\lvert a'_kA_k(x)- b_k\rvert \leq \epsilon)\to 0$ again. On the other hand, by the continuity of $A_k$ it holds $A_k(x')\to A_k(x)$ hence $\Vert a_k\Vert_2\Vert A_k(x)-A_k(x')\Vert_2=o_p(1)$, leading to $\mathbb P_F( \Vert a_k\Vert_2\Vert A_k(x)-A_k(x')\Vert_2>\epsilon)\to 0$ as $x'\to x$ for each fixed $\epsilon$. By sending $\epsilon$ to $0$, we show $\mathbb P_F(\lvert a'_kA_k(x)- b_k\rvert \leq \Vert a_k\Vert_2\Vert A_k(x)-A_k(x')\Vert_2)\to 0$.\hfill\Halmos

\proof{Proof of Proposition \ref{monotonic objective value}.}For any $s_1<s_2$, $v(s_1)\leq v(s_2)$ follows trivially from the monotonicity property $\mathrm{Sol}(s_2)\subseteq\mathrm{Sol}(s_1)$. Furthermore, if $x^*(s_1)$ and $x^*(s_2)$ are the unique optimal solutions for $OPT(s_1)$ and $OPT(s_2)$ respectively and they are distinct, then we have $v(s_1)=f(x^*(s_1))<f(x^*(s_2))=v(s_2)$ because $x^*(s_2)$ is feasible but not optimal for $OPT(s_1)$. Otherwise if $x^*(s_1)=x^*(s_2)$ then obviously $v(s_1)=f(x^*(s_1))=f(x^*(s_2))=v(s_2)$.\hfill\Halmos

\proof{Proof of Proposition \ref{unique optimum}.}A consequence of Assumption \ref{continuous constraint} is the continuity of $H(x)$ on $\mathcal X$ because for every $x',x$ it holds $\lvert H(x')-H(x)\rvert\leq \mathbb E_F[\lvert h(x',\xi)-h(x,\xi)\rvert]\leq \sqrt{\mathbb E_F[\lvert h(x',\xi)-h(x,\xi)\rvert^2]}$.

We prove the uniqueness of optimal solution by contradiction. Suppose there are $x_1\neq x_2$ and both $x_1,x_2\in\mathcal X_S^*$. Case (i): there are $s_1,s_2\notin\{\tilde s_1,\ldots,\tilde s_{M-1}\}$ such that $x_1=x^*(s_1),x_2=x^*(s_2)$. In this case we must have $s_1\neq s_2$ hence $f(x_1)\neq f(x_2)$ by Proposition \ref{monotonic objective value}, contradicting with the fact that both $x_1,x_2$ are optimal. Case (ii): there exists some $s_1\notin\{\tilde s_1,\ldots,\tilde s_{M-1}\}$ such that $x_1=x^*(s_1)$, and $x_2\in x^*(\tilde s_{i^*})$ for some $1\leq i^*\leq M-1$ but $x_2\neq x^*(s)$ for all $s\notin\{\tilde s_1,\ldots,\tilde s_{M-1}\}$. Since $x_2$ is feasible we have $H(x_2)\geq\gamma$ on one hand. On the other hand, $H(x_2)\neq \gamma$ due to Assumption \ref{non-binding}, therefore we must have $H(x_2)>\gamma$. We argue that it must be the case that $x_2=x^*(\tilde s_{i^*}+)$. If $x_2=x^*(\tilde s_{i^*}-)$ then as $s\to\tilde s_{i^*}-$ we must have $x^*(s)\neq x_2$ and $x^*(s)\to x_2$, therefore by the continuity of $H(x)$ there exist $s_1'<s_2'<\tilde s_{i^*}$ such that $H(x^*(s_1'))>\gamma,H(x^*(s_2'))>\gamma$ and $x^*(s_1')\neq x^*(s_2')$. For such $s_1',s_2'$ we have $f(x^*(s_1'))<f(x^*(s_2'))\leq f(x_2)$ from Proposition \ref{monotonic objective value}, i.e., $x^*(s_1')$ is a feasible solution with strictly less objective value than $x_2$, contradicting with the optimality of $x_2$. Hence $x_2=x^*(\tilde s_{i^*}+)$ must hold. If $s_1<\tilde s_{i^*}$, we argue that $f(x_1)<f(x_2)$ hence arrive at a contradiction. Note that the feasible set $\mathrm{Sol}(s_1)$ is closed, that $\mathrm{Sol}(s)\subseteq\mathrm{Sol}(s_1)$ for all $s>s_1$, and that $x_2=\lim_{s\to \tilde s_{i^*}}x^*(s)$ with each $x^*(s)\in \mathrm{Sol}(s_1)$, hence $x_2\in \mathrm{Sol}(s_1)$. Since $x_2\neq x^*(s_1)=x_1$ we must have $f(x_1)<f(x_2)$ by the uniqueness of $x^*(s_1)$ for $OPT(s_1)$. Otherwise if $s_1>\tilde s_{i^*}$, we take an $s\in (\tilde s_{i^*},s_1)$ sufficiently close to $\tilde s_{i^*}$ so that $x^*(s)$ is sufficiently close to $x_2$ and $x^*(s)\neq x_1=x^*(s_1)$, then from Proposition \ref{monotonic objective value} we have $f(x_2)\leq f(x^*(s))<f(x_1)$, a contradiction again. Case (iii): there are $\tilde s_{i^*_1},\tilde s_{i^*_2}$ such that $x_1\in x^*(\tilde s_{i^*_1})$ and $x_2\in x^*(\tilde s_{i^*_2})$, but there is no $s\notin \{\tilde s_1,\ldots,\tilde s_{M-1}\}$ such that $x_1=x^*(s)$ or $x_2=x^*(s)$. By the same argument in Case (ii), we can show that it must be the case that $x_1=x^*(\tilde s_{i^*_1}+)$ and $x_2=x^*(\tilde s_{i^*_2}+)$, therefore $H(x_1),H(x_2)>\gamma$. Assume $\tilde s_{i^*_1}<\tilde s_{i^*_2}$ without loss of generality, and consider an $s\in (\tilde s_{i^*_1},\tilde s_{i^*_2})$ that is sufficiently close to $\tilde s_{i^*_1}$ so that $H(x^*(s))>\gamma$, then by Proposition \ref{monotonic objective value} we have $f(x_1)=f(x^*(s))=f(x_2)$ hence $x^*(s)\in\mathcal X_S^*$, and we are in Case (ii) again. This proves that $\mathcal X_S^*$ must be a singleton.

To show that the optimal parameter set $S^*$ must be a closed interval, we first observe that $S^*$ must be a closed set due to the continuity of the solution path. Let $s_l^*=\min\{s:s\in S^*\}$ and $s_u^*=\max\{s:s\in S^*\}$, then we have $S^*\subseteq [s_l^*,s_u^*]$. Case (i): both $s_l^*,s_u^*\notin \{\tilde s_1,\ldots,\tilde s_{M-1}\}$. Note that $x^*(s_l^*)=x^*(s_u^*)=x^*_S$ and $v(s_l^*)=v(s_u^*)$, hence Proposition \ref{monotonic objective value} then forces $x^*(s)=x^*_S$ hence $s\in S^*$ for all $s\in [s_l^*,s_u^*]-\{\tilde s_1,\ldots,\tilde s_{M-1}\}$ because otherwise $v(s_l^*)<v(s)<v(s_u^*)$. This further implies $x^*(\tilde s_i-)=x^*(\tilde s_i+)=\{x^*_S\}$ and subsequently $\tilde s_i\in S^*$ for any $\tilde s_i\in [s_l^*,s_u^*]$. Altogether we have $[s_l^*,s_u^*]=S^*$. Case (ii): $s_l^*=\tilde s_{i^*}$ for some $1\leq i^*\leq M-1$ and $s_u^*\notin \{\tilde s_1,\ldots,\tilde s_{M-1}\}$. Using Proposition \ref{monotonic objective value} as in Case (i), one can show that for every $s\in [s_l^*,s_u^*]-\{\tilde s_1,\ldots,\tilde s_{M-1}\}$ we have $v(s_l^*)=v(s)=v(s_u^*)$ and $x^*(s)=x_S^*$, therefore $S^*=[s_l^*,s_u^*]$ again. Case (iii): $s_u^*=\tilde s_{i^*}$ for some $1\leq i^*\leq M-1$ and $s_l^*\notin \{\tilde s_1,\ldots,\tilde s_{M-1}\}$. This case resembles Case (ii) and $S^*=[s_l^*,s_u^*]$ can be shown using the same argument. Case (iv): $s_l^*=\tilde s_{i_1^*},s_u^*=\tilde s_{i_2^*}$ for some $1\leq i_1^*\leq i_2^*\leq M-1$. If $\tilde s_{i_1^*}=\tilde s_{i_2^*}$ then $x_S^*$ is a singleton and the interval representation trivially holds, so we focus on the case $\tilde s_{i_1^*}<\tilde s_{i_2^*}$. We argue that $x^*(\tilde s_{i_1^*}+)=x_S^*$. Otherwise if $x^*(\tilde s_{i_1^*}-)=x_S^*$, then Assumption \ref{non-binding} forces $H(x^*(\tilde s_{i_1^*}-))>\gamma$, and the continuity of the solution path and the constraint $H$ imply that $H(x^*(s'))>\gamma$ for some $s'$ sufficiently close to $\tilde s_{i_1^*}$ but $s'<\tilde s_{i_1^*}$. Note that such an $s'$ corresponds to an objective value $f(x^*(s'))=v(s')\leq \lim_{s\to \tilde s_{i_1^*}-}v(s)= f(x_S^*)$ by Proposition \ref{monotonic objective value}, therefore $x^*(s')=x_S^*$ by the uniqueness of $\mathcal X_S^*$, contradicting with the definition of $s_l^*$. Therefore it must be the case that $x^*(\tilde s_{i_1^*}+)=x_S^*$. Because $H(x^*(\tilde s_{i_1^*}+))>\gamma$, there exists a $\delta>0$ so that $\tilde s_{i_1^*}+\delta<\tilde s_{i_2^*}$, $\tilde s_{i_1^*}+\delta\notin \{\tilde s_1,\ldots,\tilde s_{M-1}\}$, and $H(x^*( s))>\gamma$ for all $s\in (\tilde s_{i_1^*},\tilde s_{i_1^*}+\delta]$. Since Proposition \ref{monotonic objective value} implies $v(s)=f(x_S^*)$ for such $s$, we must have $x^*(s)=x_S^*$ hence $[s_l^*,s_l^*+\delta]\subseteq S^*$. The rest part $[s_l^*+\delta,s_u^*]\subseteq S^*$ can be shown by treating $s_l^*+\delta$ as the $s_l^*$ in Case (iii). Altogether we still have $S^*=[s_l^*,s_u^*]$. In particular, when $v(s)$ is strictly monotonic, it is clear that there can be at most one optimal parameter hence $S^*$ becomes a singleton.\hfill\Halmos

\proof{Proof of Theorem \ref{asymptotic joint:gaussian}.}For any function class $\mathcal G$ and $g\in\mathcal G$, we write $P(g)=\mathbb E_F[g(\xi)]$ (or just $Pg$) and $P_{n_2}(g) = \frac{1}{n_2}\sum_{i=1}^{n_2}g(\xi_i)$ (or just $P_{n_2}g$), as functions from $\mathcal G\to \R$. For any function $\phi:\mathcal G\to \R$, define $\Vert \phi \Vert_{\mathcal G}=\sup_{g\in G}\lvert \phi(g)\rvert$. For example $\Vert P_{n_2}-P \Vert_{\mathcal G}$ denotes the maximal deviation of the sample mean.

First we show uniform convergence of the standard-deviation-adjusted sample mean to the expected constraint value. The function class $\{h(x,\cdot)\vert x\in\mathcal X\}$ is $F$-Donsker by Assumption \ref{Donsker h}, hence is $F$-Glivenko-Cantelli(GC). By Lemma 2.10.14 from \cite{van1996weak}, the squared class $\{h^2(x,\cdot)\vert x\in\mathcal X\}$ is also $F$-GC under Assumptions \ref{Donsker h} and \ref{moment2}. Define $\mathcal X_S=\{x^*(s):s\in S\backslash\{\tilde s_1,\ldots,\tilde s_{M-1}\}\}\cup \big(\cup_{i=1}^{M-1}x^*(\tilde s_i)\big)$. As sub-classes, $\mathcal F_{\mathcal X_S}:=\{h(x,\cdot)\vert x\in\mathcal X_S\}$ and $\mathcal F_{\mathcal X_S}^2:=\{h^2(x,\cdot)\vert x\in\mathcal X_S\}$ are both $F$-GC, i.e.,
\begin{align*}
&\norm{P_{n_2}-P}_{\mathcal F_{\mathcal X_S}}\to 0\text{\ a.s.}\\
&\norm{P_{n_2}-P}_{\mathcal F_{\mathcal X_S}^2}\to 0\text{\ a.s..}
\end{align*}
Letting $\hat{\sigma}^2(h)=P_{n_2}(h^2)-(P_{n_2}(h))^2$ and $\sigma^2(h)=P(h^2)-(P(h))^2$ be the sample and true variances, we have
\begin{eqnarray}
\notag\norm{\hat{\sigma}^2-\sigma^2}_{\mathcal F_{\mathcal X_S}}&\leq &\norm{P_{n_2}-P}_{\mathcal F_{\mathcal X_S}^2}+\norm{P_{n_2}-P}_{\mathcal F_{\mathcal X_S}}^2+2\norm{P}_{\mathcal F_{\mathcal X_S}}\norm{P_{n_2}-P}_{\mathcal F_{\mathcal X_S}}\\
\notag&=&\norm{P_{n_2}-P}_{\mathcal F_{\mathcal X_S}^2}+\norm{P_{n_2}-P}_{\mathcal F_{\mathcal X_S}}^2+2\sup_{x\in\mathcal X_S}\abs{H(x)}\norm{P_{n_2}-P}_{\mathcal F_{\mathcal X_S}}\\
&\to& 0\text{\ a.s.}\label{uniform convergence sample variance}
\end{eqnarray}
where the limit comes from the fact that $\sup_{x\in\mathcal X_S}\abs{H(x)}<\infty$ because $H$ is continuous (implied by Assumption \ref{continuous constraint}) and $\mathcal X_S$ is compact (implied by the piecewise uniform continuity condition, i.e., Assumption \ref{piecewise uniform continuous solution curve}). By Assumption \ref{moment2} we have $\norm{\sigma^2}_{\mathcal F_{\mathcal X_S}}=\sup_{x\in\mathcal X_S}\mathrm{Var}(h(x,\xi))\leq \mathbb E[\sup_{x\in\mathcal X}h^2(x,\xi)]<\infty$, and arrive at
\begin{equation*}
\norm{P_{n_2}-\frac{z_{1-\beta}}{\sqrt{n_2}}\hat{\sigma}-P}_{\mathcal F_{\mathcal X_S}}\leq \norm{P_{n_2}-P}_{\mathcal F_{\mathcal X_S}}+\frac{z_{1-\beta}}{\sqrt{n_2}}\sqrt{\norm{\hat{\sigma}^2}_{\mathcal F_{\mathcal X_S}}}\to 0\text{\ a.s.}.
\end{equation*}
When we use a discrete mesh $\{s_1,\ldots,s_p\}$, it is clear that, using the notations from Algorithm \ref{algo:marginal} and $H_j:=H(x^*(s_j))$
\begin{equation}\label{uniform convergence with margin}
\max_{1\leq j\leq p} \abs{\hat H_j-\frac{z_{1-\beta}}{\sqrt{n_2}}\hat{\sigma}_j-H_j}\leq \norm{P_{n_2}-\frac{z_{1-\beta}}{\sqrt{n_2}}\hat{\sigma}-P}_{\mathcal F_{\mathcal X_S}}\to 0\text{\ a.s.}.
\end{equation}

Secondly, we prove convergence of the estimated solution $x^*(\hat s^*)$ to the optimum $x_S^*$. Fixing any $\epsilon>0$, we argue that almost surely there exists a finite $N$ and $\overline{\epsilon_S}>0$ such that for all $n_2\geq N$ and $\epsilon_S\leq\overline{\epsilon_S}$ it holds $\norm{x^*(\hat s^*)-x_S^*}<\epsilon$. To proceed, define
\begin{equation}\label{growth}
\delta:=\min_{x\in\mathcal X_S}\{f(x)-f(x_S^*)\vert H(x)\geq\gamma,\norm{x-x_S^*}\geq \epsilon\}.
\end{equation}
Since the objective $f$ is continuous and $\{x\vert x\in\mathcal X_S,H(x)\geq\gamma,\norm{x-x_S^*}\geq \epsilon\}$ is a compact set, by a compactness argument we must have $\delta>0$. By Assumption \ref{unique optimal solution}, for any $\epsilon'\leq \epsilon$ there exists some $s'\notin \{\tilde s_1,\ldots,\tilde s_{M-1}\}$ such that $H(x^*(s'))>\gamma$ and $\norm{x^*(s')-x_S^*}<\epsilon'$. By continuity of $f$, one can set $\epsilon'$ small enough so that $f(x^*(s'))-f(x_S^*)< \frac{\delta}{2}$. Moreover, due to the continuity of $x^*(s)$ at $s'$ and the continuity of $f$ and $H$, there exists an $\overline{\epsilon_S}>0$ such that $\min_{s\text{ s.t. }\lvert s-s'\rvert\leq \overline{\epsilon_S}}H(x^*(s))>\gamma$ and $\max_{s\text{ s.t. }\lvert s-s'\rvert\leq \overline{\epsilon_S}}f(x^*(s))<f(x_S^*)+ \frac{\delta}{2}$. Therefore, when the mesh size $\epsilon_S\leq \overline{\epsilon_S}$, there must exist some $s_{j'}\in \{s_1,\ldots,s_p\}$ such that
\begin{align}
 H(x^*(s_{j'}))&\geq \min_{s\text{ s.t. }\lvert s-s'\rvert\leq \overline{\epsilon_S}}H(x^*(s))>\gamma\label{uniform lower bound of H}\\
f(x^*(s_{j'}))&\leq \max_{s\text{ s.t. }\lvert s-s'\rvert\leq \overline{\epsilon_S}}f(x^*(s))<f(x_S^*)+ \frac{\delta}{2}.\label{feasible near-optimal}
\end{align}
For the given $\epsilon$ define for $\Delta\geq 0$
\begin{equation}\label{optimality gap with relaxed constraint}
\delta_{\Delta}:=\min_{x\in\mathcal X_S}\{f(x)-f(x_S^*)\vert H(x)\geq\gamma-\Delta,\norm{x-x_S^*}\geq \epsilon\}.
\end{equation}
We argue that $\lim_{\Delta\to 0+}\delta_{\Delta}\to\delta$ by contradiction. Clearly $\delta_{\Delta}$ is non-increasing in $\Delta$ and $\delta_{\Delta}\leq \delta$, hence the limit $\lim_{\Delta\to 0+}\delta_{\Delta}$ must exists and is finite. Suppose $\lim_{\Delta\to 0+}\delta_{\Delta}<\delta$, then there exist $\tilde{\delta}<\delta$ and a sequence $\{x_k\}_{k=1}^{\infty}\subset \mathcal X_S$ such that $\norm{x_k-x_S^*}\geq \epsilon$, $H(x_k)\geq \gamma-\Delta_k$ with $\Delta_k\to 0+$, and $f(x_k)-f(x_S^*)\leq \tilde{\delta}$. By the compactness of $\mathcal X_S$, there must exist a subsequence $\{x_{k_s}\}_{s=1}^{\infty}$ converging to some $x_{\infty}\in\mathcal X_S$, and by continuity $x_{\infty}$ must satisfy $\norm{x_{\infty}-x_S^*}\geq \epsilon$, $H(x_{\infty})\geq \gamma$ and $f(x_{\infty})-f(x_S^*)\leq \tilde{\delta}$. From the definition \eqref{growth} of $\delta$ this implies $\delta\leq \tilde{\delta}$, a contradiction. Now pick a small enough $\Delta'$ so that $\delta_{\Delta'}>\frac{\delta}{2}$. From the uniform convergence \eqref{uniform convergence with margin} we know that almost surely there exists some $N$ such that for any $n_2\geq N$
\begin{equation}\label{uniform error}
\max_{1\leq j\leq p} \abs{\hat H_j-\frac{z_{1-\beta}}{\sqrt{n_2}}\hat{\sigma}_j-H_j}< \min\{\Delta',\min_{s\text{ s.t. }\lvert s-s'\rvert\leq \overline{\epsilon_S}}H(x^*(s))-\gamma\}\leq \min\{\Delta',H(x^*(s_{j'}))-\gamma\}.
\end{equation}
where the second inequality is due to \eqref{uniform lower bound of H}. In particular, \eqref{uniform error} implies that for such $n_2$ we have $\hat H(x^*(s_{j'}))-\frac{z_{1-\beta}}{\sqrt{n_2}}\hat{\sigma}(x^*(s_{j'}))>H(x^*(s_{j'}))-\min\{\Delta',H(x^*(s_{j'}))-\gamma\}\geq \gamma$,
therefore on one hand we must have
\begin{equation}\label{bound:objective value}
f(x^*(\hat s^*))\leq f(x^*(s_{j'}))<f(x_S^*)+\frac{\delta}{2}
\end{equation}
where the first inequality holds due to the way $\hat s^*$ is chosen and the second results from \eqref{feasible near-optimal}. On the other hand it also follows from \eqref{uniform error} that
\begin{eqnarray}
\notag H(x^*(\hat s^*))&>&\hat H(x^*(\hat s^*))-\frac{z_{1-\beta}}{\sqrt{n_2}}\hat{\sigma}(x^*(\hat s^*))-\min\{\Delta',H(x^*(s_{j'}))-\gamma\}\\
\notag&\geq& \gamma-\min\{\Delta',H(x^*(s_{j'}))-\gamma\}\\
&\geq& \gamma-\Delta'.\label{bound:constraint}
\end{eqnarray}
The bounds \eqref{bound:objective value} and \eqref{bound:constraint} on the objective value and the constraint value at the estimated solution $x^*(\hat s^*)$, together with the fact that $\delta_{\Delta'}>\frac{\delta}{2}$ due to the way $\Delta'$ is chosen, imply that $\norm{x^*(\hat s^*)-x_S^*}<\epsilon$ by the definition \eqref{optimality gap with relaxed constraint} of $\delta_{\Delta}$. Since $\epsilon$ can be arbitrarily small, we have $\lim_{n_2\to\infty}x^*(\hat{s}^*)=x_S^*$ a.s.. Convergence of $\hat s^*$ to the optimal parameter set $\mathcal S^*$ is then a consequence of the convergence of $x^*(\hat s^*)$ to $x_S^*$. Suppose $\hat s^*$ does not converge to $S^*$, then by compactness of $S$ there exists a subsequence $\hat s^*_{k}$ converging to some $s_{\infty}\notin S^*$. Since the corresponding $x^*(\hat s^*_k)\to x_S^*$ we have either $s_{\infty}\notin\{\tilde s_1,\ldots,\tilde s_{M-1}\}$ with $x^*(s_{\infty})=x_S^*$ or $s_{\infty}\in\{\tilde s_1,\ldots,\tilde s_{M-1}\}$ with $x_S^*\in x^*(s_{\infty})$, however in either case $s_{\infty}\in S^*$, a contradiction.


Then we prove the feasibility guarantees. The case $H(x_S^*)>\gamma$ is relatively straightforward. By the continuity of $H$ and that a.s. $x^*(\hat s^*)\to x_S^*$ we have $H(x^*(\hat s^*))\to H(x_S^*)>\gamma$ a.s.. Almost surely convergence implies convergence in probability, thus $H(x^*(\hat s^*))\to H(x_S^*)$ in probability and, in particular, $P_{\bm{\xi}_{1:n_2}}(H(x^*(\hat s^*))\geq \gamma)\to 1$. If $H(x_S^*)=\gamma$ we denote by
$$\mathbb G_{n_2}(x):=\sqrt{n_2}(P_{n_2}(h(x,\cdot))-P(h(x,\cdot)))$$
the empirical process indexed by the decision $x$ and let $\hat\sigma^2(x),\sigma^2(x)$ represent the sample and true variance of $h(x,\xi)$, and then write
\begin{eqnarray}
H(x^*(\hat s^*))&\geq& H(x^*(\hat s^*))-\big(\hat H(x^*(\hat s^*))-z_{1-\beta}\frac{\hat \sigma(x^*(\hat s^*))}{\sqrt{n_2}}-\gamma\big)\label{towards tight coverage}\\
\notag&=&\gamma+(H(x^*(\hat s^*))-\hat H(x^*(\hat s^*)))+z_{1-\beta}\frac{\hat \sigma(x^*(\hat s^*))}{\sqrt{n_2}}\\
\notag&=&\gamma-\frac{1}{\sqrt{n_2}}\mathbb G_{n_2}(x^*(\hat s^*))+z_{1-\beta}\frac{\hat \sigma(x^*(\hat s^*))}{\sqrt{n_2}}\\
&=&\gamma-\frac{1}{\sqrt{n_2}}\mathbb G_{n_2}(x_S^*)+z_{1-\beta}\frac{\hat \sigma (x_S^*)}{\sqrt{n_2}}+E_1+E_2\label{representation of constraint value}
\end{eqnarray}
where the errors
\begin{equation*}
E_1= \frac{1}{\sqrt{n_2}}\big(\mathbb G_{n_2}(x_S^*)-\mathbb G_{n_2}(x^*(\hat s^*))\big),\ E_2=\frac{z_{1-\beta}}{\sqrt{n_2}}\big(\hat \sigma(x^*(\hat s^*))-\hat\sigma(x_S^*)\big).
\end{equation*}

We need to show that $E_1=o_p\big(\frac{1}{\sqrt{n_2}}\big),E_2=o_p\big(\frac{1}{\sqrt{n_2}}\big)$. We deal with $E_2$ first. $E_2$ can be bounded as
\begin{equation*}
\abs{E_2}\leq \frac{z_{1-\beta}}{\sqrt{n_2}}\big(2\norm{\hat{\sigma}-\sigma}_{\mathcal F_{\mathcal X_S}}+\abs{\sigma(x^*(\hat s^*))-\sigma(x_S^*)}\big).
\end{equation*}
On one hand we have already shown that $x^*(\hat s^*)\to x_S^*$ a.s.. On the other hand, $\sigma^2(x)$ is continuous in $x$. Therefore $\abs{\sigma(x^*(\hat s^*))-\sigma(x_S^*)}\to 0$ a.s. as $n_2\to\infty$. By uniform convergence \eqref{uniform convergence sample variance} and the relation $\norm{\hat{\sigma}-\sigma}_{\mathcal F_{\mathcal X_S}}\leq\sqrt{\norm{\hat{\sigma}^2-\sigma^2}_{\mathcal F_{\mathcal X_S}}}$, we have $\norm{\hat{\sigma}-\sigma}_{\mathcal F_{\mathcal X_S}}\to 0$ a.s.. Consequently it holds $\sqrt{n_2}\abs{E_2}=o(1)$ a.s. and, in particular, $E_2=o_p\big(\frac{1}{\sqrt{n_2}}\big)$.

To bound the error $E_1$, let $\rho(x,x')=\sqrt{\mathrm{Var}(h(x,\xi)-h(x',\xi))}$ denote the intrinsic semimetric of the tight Gaussian process $\mathbb G$ indexed by $x\in\mathcal X_S$ with mean zero and covariance structure $\mathrm{Cov}(\mathbb G(x),\mathbb G(x'))=\mathrm{Cov}_F(h(x,\xi),h(x',\xi))$, and for any $\epsilon>0$ let $\delta(\epsilon) = \sup\{\rho(x,x_S^*)\vert x\in\mathcal X_S,\norm{x-x_S^*}<\epsilon\}$. Note that Assumption \ref{continuous constraint} entails $\delta(\epsilon)\to 0$ as $\epsilon\to 0$. We have for any $\epsilon>0$
\begin{eqnarray*}
\abs{E_1}&=&\abs{E_1}\mathbf{1}\{\norm{x^*(\hat s^*)-x_S^*}<\epsilon\}+\abs{E_1}\mathbf{1}\{\norm{x^*(\hat s^*)-x_S^*}\geq\epsilon\}\\
&\leq& \sup_{x\in\mathcal X_S\text{ s.t. }\norm{x-x_S^*}< \epsilon}\frac{1}{\sqrt{n_2}}\abs{\mathbb G_{n_2}(x_S^*)-\mathbb G_{n_2}(x)}\mathbf{1}\{\norm{x^*(\hat s^*)-x_S^*}<\epsilon\}+\infty \cdot\mathbf{1}\{\norm{x^*(\hat s^*)-x_S^*}\geq\epsilon\}\\
&&\text{where }\infty\cdot 0=0\\
&\leq& \sup_{x\in\mathcal X_S\text{ s.t. }\norm{x-x_S^*}< \epsilon}\frac{1}{\sqrt{n_2}}\abs{\mathbb G_{n_2}(x_S^*)-\mathbb G_{n_2}(x)}+\infty \cdot\mathbf{1}\{\norm{x^*(\hat s^*)-x_S^*}\geq\epsilon\}\\
&\leq&\sup_{x,x'\in\mathcal X_S\text{ s.t. }\rho(x,x')\leq \delta(\epsilon)}\frac{1}{\sqrt{n_2}}\abs{\mathbb G_{n_2}(x)-\mathbb G_{n_2}(x')}+\infty \cdot\mathbf{1}\{\norm{x^*(\hat s^*)-x_S^*}\geq\epsilon\}.
\end{eqnarray*}
We have already shown that $\norm{x^*(\hat s^*)-x_S^*}\to 0$ a.s., hence $P_{\bm{\xi}_{1:n_2}}(\norm{x^*(\hat s^*)-x_S^*}\geq\epsilon)\to 0$ for any fixed $\epsilon>0$. Therefore we can choose an $n_2$-dependent $\epsilon:=\epsilon_{n_2}$ such that both $P_{\bm{\xi}_{1:n_2}}(\norm{x^*(\hat s^*)-x_S^*}\geq\epsilon_{n_2})\to 0$ and $\epsilon_{n_2}\to 0$ as $n_2\to\infty$, and get
\begin{equation}
\abs{E_1}\leq \sup_{x,x'\in\mathcal X_S\text{ s.t. }\rho(x,x')\leq \delta(\epsilon_{n_2})}\frac{1}{\sqrt{n_2}}\abs{\mathbb G_{n_2}(x)-\mathbb G_{n_2}(x')}+\infty \cdot\mathbf{1}\{\norm{x^*(\hat s^*)-x_S^*}\geq\epsilon_{n_2}\}.\label{error1 bound}
\end{equation}
By the way $\epsilon_{n_2}$ is chosen, the second term on the right hand side of \eqref{error1 bound} is of arbitrarily small order, in particular, $o_p\big(\frac{1}{\sqrt{n_2}}\big)$. To control the first term, note that $\delta(\epsilon_{n_2})\to 0$ as $n_2\to\infty$. Since the function class $\mathcal F_{\mathcal X_S}$ is $F$-Donsker, the empirical process $\mathbb G_{n_2}$ on $\mathcal F_{\mathcal X_S}$ is asymptotically tight, hence by Theorem 1.5.7 and Addendum 1.5.8 from \cite{van1996weak} $\mathbb G_{n_2}$ is asymptotically uniformly equicontinuous in probability with respect to the intrinsic semimetric $\rho$ of the limit Gaussian process $\mathbb G$, i.e., for any $\epsilon>0$
\begin{equation}\label{uniform equicontinuity in probability}
\lim_{\delta\to 0}\limsup_{n_2\to\infty}P_{\bm{\xi}_{1:n_2}}\Big(\sup_{x,x'\in\mathcal X_S\text{ s.t. }\rho(x,x')\leq \delta}\abs{\mathbb G_{n_2}(x)-\mathbb G_{n_2}(x')}> \epsilon\Big)=0.
\end{equation}
Note that $\sup_{x,x'\in\mathcal X_S\text{ s.t. }\rho(x,x')\leq \delta}\abs{\mathbb G_{n_2}(x)-\mathbb G_{n_2}(x')}$ is monotonically increasing in $\delta$ a.s. and $\delta(\epsilon_{n_2})\to 0$, it must hold that for any fixed $\delta>0$
\begin{equation*}
\sup_{x,x'\in\mathcal X_S\text{ s.t. }\rho(x,x')\leq \delta(\epsilon_{n_2})}\abs{\mathbb G_{n_2}(x)-\mathbb G_{n_2}(x')}\leq\sup_{x,x'\in\mathcal X_S\text{ s.t. }\rho(x,x')\leq \delta}\abs{\mathbb G_{n_2}(x)-\mathbb G_{n_2}(x')}\text{\ \ a.s.}
\end{equation*}
when $n_2$ is sufficiently large, therefore for any $\epsilon>0$ the first term in \eqref{error1 bound} can be controlled as
\begin{eqnarray}
\notag&&\limsup_{n_2\to\infty} P_{\bm{\xi}_{1:n_2}}\Big(\sup_{x,x'\in\mathcal X_S\text{ s.t. }\rho(x,x')\leq \delta(\epsilon_{n_2})}\abs{\mathbb G_{n_2}(x)-\mathbb G_{n_2}(x')}>\epsilon\Big)\\
&\leq&\limsup_{n_2\to\infty} P_{\bm{\xi}_{1:n_2}}\Big(\sup_{x,x'\in\mathcal X_S\text{ s.t. }\rho(x,x')\leq \delta}\abs{\mathbb G_{n_2}(x)-\mathbb G_{n_2}(x')}>\epsilon\Big).\label{continuity of empirical process}
\end{eqnarray}
Due to \eqref{uniform equicontinuity in probability} the right hand side of \eqref{continuity of empirical process} can be made arbitrarily small by sending $\delta\to 0$, hence the left hand side of \eqref{continuity of empirical process} must be identical to zero. Since $\epsilon$ is arbitrary, by definition $\sup_{x,x'\in\mathcal X_S\text{ s.t. }\rho(x,x')\leq \delta(\epsilon_{n_2})}\abs{\mathbb G_{n_2}(x)-\mathbb G_{n_2}(x')}=o_p(1)$, which in turn leads to $E_1 = o_p\big(\frac{1}{\sqrt{n_2}}\big)$.

We now go back to the representation \eqref{representation of constraint value} of $H(x^*(\hat s^*))$ to conclude the coverage guarantee. From \eqref{representation of constraint value} we see that
\begin{eqnarray*}
&&\liminf_{n_2\to\infty,\epsilon_S\to 0}P_{\bm{\xi}_{1:n_2}}(H(x^*(\hat s^*))\geq \gamma)\\
&\geq &\liminf_{n_2\to\infty,\epsilon_S\to 0}P_{\bm{\xi}_{1:n_2}}\Big(\gamma-\frac{1}{\sqrt{n_2}}\mathbb G_{n_2}(x_S^*)+z_{1-\beta}\frac{\hat \sigma (x_S^*)}{\sqrt{n_2}}+E_1+E_2\geq \gamma\Big)\\
&= &\liminf_{n_2\to\infty,\epsilon_S\to 0}P_{\bm{\xi}_{1:n_2}}\Big(-\frac{1}{\sqrt{n_2}}\mathbb G_{n_2}(x_S^*)+z_{1-\beta}\frac{\hat \sigma (x_S^*)}{\sqrt{n_2}}+o_p\big(\frac{1}{\sqrt{n_2}}\big)\geq 0\Big)\\
&= &\liminf_{n_2\to\infty,\epsilon_S\to 0}P_{\bm{\xi}_{1:n_2}}\Big(\frac{\mathbb G_{n_2}(x_S^*)}{\hat \sigma (x_S^*)}+o_p(1)\leq z_{1-\beta}\Big)\text{\ \ since }\hat \sigma (x_S^*)\to \sigma (x_S^*)\text{ a.s. and }\sigma^2 (x_S^*)>0\text{ (Assumption \ref{non-degenerate variance})}\\
&=&1-\beta
\end{eqnarray*}
where in the last equality we use Slutsky's theorem to justify that $\frac{\mathbb G_{n_2}(x_S^*)}{\hat \sigma (x_S^*)}+o_p(1)$ weakly converges to the standard normal.\hfill\Halmos

\proof{Proof of Theorem \ref{asymptotic tight coverage:gaussian}.}Following the proof of Theorem \ref{asymptotic joint:gaussian}, we see that in order to conclude the tight feasibility confidence level it suffices to show that the inequality gap of \eqref{towards tight coverage} is of order $o_p\big(\frac{1}{\sqrt{n_2}}\big)$, i.e.,
\begin{equation}\label{asymptotic binding}
\gamma\leq\hat H(x^*(\hat s^*))-z_{1-\beta}\frac{\hat \sigma(x^*(\hat s^*))}{\sqrt{n_2}}\leq \gamma+o_p\big(\frac{1}{\sqrt{n_2}}\big).
\end{equation}
Indeed, once the second inequality in \eqref{asymptotic binding} is shown, we can use the representation \eqref{representation of constraint value} and apply Slutsky's theorem, like in the proof of Theorem \ref{asymptotic joint:gaussian}, to get
\begin{eqnarray*}
&&\lim P_{\bm{\xi}_{1:n_2}}(H(x^*(\hat s^*))\geq \gamma)\\
&= &\lim P_{\bm{\xi}_{1:n_2}}\Big(\gamma-\frac{1}{\sqrt{n_2}}\mathbb G_{n_2}(x_S^*)+z_{1-\beta}\frac{\hat \sigma (x_S^*)}{\sqrt{n_2}}+E_1+E_2+o_p\big(\frac{1}{\sqrt{n_2}}\big)\geq \gamma\Big)\\
&= &\lim P_{\bm{\xi}_{1:n_2}}\Big(-\frac{1}{\sqrt{n_2}}\mathbb G_{n_2}(x_S^*)+z_{1-\beta}\frac{\hat \sigma (x_S^*)}{\sqrt{n_2}}+o_p\big(\frac{1}{\sqrt{n_2}}\big)\geq 0\Big)\text{\ \ since }E_1,E_2=o_p\big(\frac{1}{\sqrt{n_2}}\big)\\
&= &\lim P_{\bm{\xi}_{1:n_2}}\Big(\frac{\mathbb G_{n_2}(x_S^*)}{\hat \sigma (x_S^*)}+o_p(1)\leq z_{1-\beta}\Big)\\
&=&1-\beta.
\end{eqnarray*}

Now we prove the second inequality in \eqref{asymptotic binding}. By Proposition \ref{unique optimum} the optimal parameter set is a singleton $S^*=\{s^*\}$. Moreover, in the case $H(x_S^*)=\gamma$ Assumption \ref{non-binding} forces that $s^*\notin \{\tilde s_1,\ldots,\tilde s_{M-1}\}\cup \{s_l,s_u\}$. Suppose $s^*\in(\tilde s_{i^*},\tilde s_{i^*+1})$ for some $0\leq i^*\leq M-1$ (note that $\tilde s_0=s_l,\tilde s_M=s_u$). Assumption \ref{piecewise uniform continuous solution curve} then ensures that the parameter-to-solution mapping $x^*(\cdot)$ is uniformly continuous in some neighborhood $\mathcal N(s^*)\subseteq (\tilde s_{i^*},\tilde s_{i^*+1})$ of $s^*$. Since $\mathcal N(s^*)$ is contained in a compact set, the standard deviation function $\sigma(x^*(\cdot))$ is uniformly continuous in $\mathcal N(s^*)$. Moreover, the semimetric $\rho(x^*(\cdot),x^*(\cdot))$ between two solutions is also uniformly continuous in $\mathcal N(s^*)\times \mathcal N(s^*)$. Therefore as $\epsilon_S\to 0$ the following holds
\begin{align*}
\omega_{\sigma}(2\epsilon_S):&=\sup_{s,s'\in\mathcal N(s^*)\text{ s.t. }\norm{s-s'}<2\epsilon_S}\abs{\sigma(x^*(s))-\sigma(x^*(s'))}=o(1)\\
\omega_{\rho}(2\epsilon_S):&=\sup_{s,s'\in\mathcal N(s^*)\text{ s.t. }\norm{s-s'}<2\epsilon_S}\rho(x^*(s),x^*(s'))=o(1).
\end{align*}
According to the criterion of choosing $\hat s^*$ we must have for every parameter value $s_j$ either $f(x^*(s_j))\geq f(x^*(\hat s^*))$ or $\hat H(x^*(s_j))-z_{1-\beta}\frac{\hat \sigma(x^*(s_j))}{\sqrt{n_2}}<\gamma$.
Therefore if $\hat s^*\in\mathcal N(s^*)$, say $\hat s^*=s^{i^*}_{j^*}$, and $s^{i^*}_{j^*-1}\in\mathcal N(s^*)$ as well, then because $s^{i^*}_{j^*-1}<s^{i^*}_{j^*}$ and the parameter-to-objective mapping $f(x^*(\cdot))$ is strictly increasing in $s$ it must hold that
\begin{equation}\label{theorem of the alternative}
\hat H(x^*(s^{i^*}_{j^*-1}))-z_{1-\beta}\frac{\hat \sigma(x^*(s^{i^*}_{j^*-1}))}{\sqrt{n_2}}<\gamma.
\end{equation}
We shall use this fact to derive \eqref{asymptotic binding}. For convenience, we denote by $B(s,\epsilon):=\{s'\in\mathcal S\vert \norm{s'-s}<\epsilon\}$ the ball of radius $\epsilon>0$ centered at $s$. Because $\hat s^*\to s^*$ a.s. and $\epsilon_S\to 0$, it is implied that $P_{\bm{\xi}_{1:n_2}}(B(\hat s^*,2\epsilon_S)\subseteq \mathcal N(s^*)) \to 1$ as $n_2\to\infty$. Thus we can write
\begin{eqnarray*}
&&\hat H(x^*(\hat s^*))-z_{1-\beta}\frac{\hat \sigma(x^*(\hat s^*))}{\sqrt{n_2}}\\
&\leq& \infty\cdot\mathbf{1}\{B(\hat s^*,2\epsilon_S)\not\subseteq \mathcal N(s^*)\}+\big(\hat H(x^*(s^{i^*}_{j^*}))-z_{1-\beta}\frac{\hat \sigma(x^*(s^{i^*}_{j^*}))}{\sqrt{n_2}}\big)\mathbf{1}\{B(\hat s^*,2\epsilon_S)\subseteq \mathcal N(s^*)\}\text{\ \ where }s^{i^*}_{j^*}=\hat s^*\\
&\leq& o_p\big(\frac{1}{\sqrt{n_2}}\big)+\big(\hat H(x^*(s^{i^*}_{j^*-1}))-z_{1-\beta}\frac{\hat \sigma(x^*(s^{i^*}_{j^*-1}))}{\sqrt{n_2}}\big)\mathbf{1}\{B(\hat s^*,2\epsilon_S)\subseteq \mathcal N(s^*)\}+\\
&& \hspace{5ex}\Big(\abs{\hat H(x^*(s^{i^*}_{j^*-1}))-\hat H(x^*(s^{i^*}_{j^*}))}+\frac{z_{1-\beta}}{\sqrt{n_2}}\abs{\hat \sigma(x^*(s^{i^*}_{j^*-1}))-\hat \sigma(x^*(s^{i^*}_{j^*}))}\Big)\mathbf{1}\{B(\hat s^*,2\epsilon_S)\subseteq \mathcal N(s^*)\}\\
&\leq& o_p\big(\frac{1}{\sqrt{n_2}}\big)+\gamma+\Big(\abs{H(x^*(s^{i^*}_{j^*-1}))-H(x^*(s^{i^*}_{j^*}))}+\frac{1}{\sqrt{n_2}}\abs{\mathbb G_{n_2}(x^*(s^{i^*}_{j^*-1}))-\mathbb G_{n_2}(x^*(s^{i^*}_{j^*}))}+\\
&&\hspace{5ex}\frac{z_{1-\beta}}{\sqrt{n_2}}\big(2\norm{\hat\sigma-\sigma}_{\mathcal F_{\mathcal X_S}}+\abs{\sigma(x^*(s^{i^*}_{j^*-1}))-\sigma(x^*(s^{i^*}_{j^*}))}\big)\Big)\mathbf{1}\{B(\hat s^*,2\epsilon_S)\subseteq \mathcal N(s^*)\}\\
&&\hspace{5ex}\text{where the $\gamma$ term comes from \eqref{theorem of the alternative}}\\
&\leq& o_p\big(\frac{1}{\sqrt{n_2}}\big)+\gamma+\Big(o\big(\frac{1}{\sqrt{n_2}}\big)+\frac{1}{\sqrt{n_2}}\sup_{x,x'\in\mathcal X_S\text{ s.t. }\rho(x,x')\leq \omega_{\rho}(2\epsilon_S)}\abs{\mathbb G_{n_2}(x)-\mathbb G_{n_2}(x')}+\\
&&\hspace{5ex}\frac{z_{1-\beta}}{\sqrt{n_2}}\big(2\norm{\hat\sigma-\sigma}_{\mathcal F_{\mathcal X_S}}+\omega_{\sigma}(2\epsilon_S)\big)\Big)\mathbf{1}\{B(\hat s^*,2\epsilon_S)\subset \mathcal N(s^*)\}\\
&&\text{\ \ where the $o\big(\frac{1}{\sqrt{n_2}}\big)$ terms is due to condition \eqref{high resolution}}\\
&\leq& o_p\big(\frac{1}{\sqrt{n_2}}\big)+\gamma+o\big(\frac{1}{\sqrt{n_2}}\big)+\frac{1}{\sqrt{n_2}}\sup_{x,x'\in\mathcal X_S\text{ s.t. }\rho(x,x')\leq \omega_{\rho}(2\epsilon_S)}\abs{\mathbb G_{n_2}(x)-\mathbb G_{n_2}(x')}+\frac{z_{1-\beta}}{\sqrt{n_2}}\big(o_p\big(1)+o(1)\big)\\
&=&\gamma+o_p\big(\frac{1}{\sqrt{n_2}}\big)+\frac{1}{\sqrt{n_2}}\sup_{x,x'\in\mathcal X_S\text{ s.t. }\rho(x,x')\leq \omega_{\rho}(2\epsilon_S)}\abs{\mathbb G_{n_2}(x)-\mathbb G_{n_2}(x')}.
\end{eqnarray*}
Since $\omega_{\rho}(2\epsilon_S)\to 0$, through an argument similar to \eqref{continuity of empirical process} the asymptotically uniform equicontinuity of $\mathbb G_{n_2}$ results in $\sup_{x,x'\in\mathcal X_S\text{ s.t. }\rho(x,x')\leq \omega_{\rho}(2\epsilon_S)}\abs{\mathbb G_{n_2}(x)-\mathbb G_{n_2}(x')}=o_p(1)$. This finally leads to the upper bound in \eqref{asymptotic binding}, hence concludes the theorem.\hfill\Halmos

\proof{Proof of Theorem \ref{asymptotic joint:gaussian supremum}.}We first treat the unnormalized validator (Algorithm \ref{algo:supremum_unnormalized}). As the first step, we introduce some notations and a few auxiliary Gaussian processes. Let $S^o:=S\backslash\{\tilde s_1,\ldots,\tilde s_{M-1}\}$ be the parameter space after excluding the pathological points $\{\tilde s_1,\ldots,\tilde s_{M-1}\}$. We denote by
$$\{\overline{\mathbb G}(s):s\in S^o\}$$
the Gaussian process with mean zero and covariance structure $\mathrm{Cov}(s,s')=\mathrm{Cov}_F(h(x^*(s),\xi),h(x^*(s'),\xi))$, and by
$$\{\overline{\mathbb G}'(s):s\in S^o\}$$
the Gaussian process with mean zero and covariance structure $\mathrm{Cov}(s,s')=\frac{1}{n_2}\sum_{i=1}^{n_2}(h(x^*(s),\xi_i)-\hat H(x^*(s)))(h(x^*(s'),\xi_i)-\hat H(x^*(s')))$ where $\hat H(x^*(s))=(1/n_2)\sum_{i=1}^{n_2}h(x^*(s),\xi_i)$ is the sample mean at $x^*(s)$ and $\hat H(x^*(s'))$ is the sample mean at $x^*(s')$. For a generic stochastic process $\{Y(\theta):\theta\in\Theta\}$ over some set $\Theta$, we denote by
$$\psi_{1-\beta}(\{Y(\theta):\theta\in\Theta\})$$
the $1-\beta$ quantile of $\sup_{\theta\in \Theta}Y(\theta)$. We can formally express the critical value calibrated in Algorithm \ref{algo:supremum_unnormalized} as $q_{1-\beta}=\psi_{1-\beta}(\{\overline{\mathbb G}'(s):s\in \{s_1,\ldots,s_p\}\})$, and $\bar{q}_{1-\beta}=\psi_{1-\beta}(\{\overline{\mathbb G}(s):s\in S^o\})$. Under Assumption \ref{Donsker h}, the Gaussian process $\overline{\mathbb G}$ as the weak limit of the empirical process $\{\sqrt{n_2}((1/n_2)\sum_{i=1}^{n_2}h(x^*(s),\xi_i)-H(x^*(s))):s\in S^o\}$ is a tight Borel measurable element in $l^{\infty}(S^o):=\{f:f\text{ is a function }S^o\to \R\text{ such that }\sup_{s\in S^o}\lvert f(s)\rvert< \infty\}$, therefore the sample path of $\overline{\mathbb G}$ is uniformly continuous with respect to the semimetric $\rho(s,s'):=\sqrt{\mathrm{Var}_F(h(x^*(s),\xi)-h(x^*(s'),\xi))}$ almost surely (Example 1.5.10 in \cite{van1996weak}). Note that, under Assumptions \ref{continuous constraint} and \ref{piecewise uniform continuous solution curve}, on each continuous piece of the solution path this semimetric is continuous in the pair $s,s'$ with respect to the Euclidean metric on $S$. In other words, almost surely the sample path of $\overline{\mathbb G}$ is continuous with respect to the Euclidean metric $d(s,s'):=\lvert s-s'\rvert$ on each piece $(\tilde s_i,\tilde s_{i+1})$. Therefore, by continuity, every countable dense (w.r.t. the Euclidean metric) subset $S^o_c\subset S^o$, e.g., the set of all rational $s$, renders
$$\sup_{s\in S^o}\overline{\mathbb G}(s)=\sup_{s\in S^o_c}\overline{\mathbb G}(s)\text{ almost surely}.$$
Suppose $S^o_c=\{s^o_j\}_{j=1}^{\infty}$, then $\max_{1\leq j\leq k}\overline{\mathbb G}(s^o_j)$ monotonically increases in $k$ towards the limit $\sup_{s\in S^o}\overline{\mathbb G}(s)$ almost surely, and almost sure convergence implies convergence in distribution therefore $\psi_{1-\beta}(\{\overline{\mathbb G}(s):s\in\{s^o_1,\ldots,s^o_k\}\})$ monotonically increases in $k$ towards the limit $\bar{q}_{1-\beta}$, i.e.,
\begin{equation}\label{separability of G}
    \psi_{1-\beta}(\{\overline{\mathbb G}(s):s\in\{s^o_1,\ldots,s^o_k\}\})\leq \bar{q}_{1-\beta},\text{ and }\lim_{k\to\infty}\psi_{1-\beta}(\{\overline{\mathbb G}(s):s\in\{s^o_1,\ldots,s^o_k\}\})=\bar{q}_{1-\beta}.
\end{equation}

As the second step, we want to show that $q_{1-\beta}$ converges to $\bar{q}_{1-\beta}$ almost surely. Under Assumptions \ref{Donsker h} and \ref{moment2}, Theorem 10.6 in \cite{kosorok2008introduction} states that, for almost every realization of the data sequence $\{\xi_i\}_{i=1}^{\infty}$, the Gaussian process $\overline{\mathbb G}'$ weakly converges to $\overline{\mathbb G}$ as $n_2\to\infty$. By the continuous mapping theorem, $\sup_{s\in S^o}\overline{\mathbb G}'(s)$ also weakly converges to $\sup_{s\in S^o}\overline{\mathbb G}(s)$ almost surely, therefore $\lim_{n_2\to\infty}\psi_{1-\beta}(\{\overline{\mathbb G}'(s):s\in S^o\})= \bar{q}_{1-\beta}$ almost surely. It is obvious that $q_{1-\beta}\leq \psi_{1-\beta}(\{\overline{\mathbb G}'(s):s\in S^o\})$, hence we have established that $\limsup_{n_2\to\infty,\epsilon_S\to 0}q_{1-\beta}\leq \bar{q}_{1-\beta}$. To show the other direction, we exploit the separability \eqref{separability of G} of $\overline{\mathbb G}$. For each $i$, let $s_{j_i}\in \{s_1,\ldots,s_p\}$ be such that $s_{j_i}\to s^o_i$ as $\epsilon_S\to 0$. Consider two more Gaussian processes $\{\overline{\mathbb G}'(s):s\in\{s_{j_1},\ldots,s_{j_k}\}\}$ and $\{\overline{\mathbb G}(s):s\in\{s_{j_1},\ldots,s_{j_k}\}\}$. For a fixed $k$, let $\hat \Sigma',\hat \Sigma$ be the covariance matrices of $\{\overline{\mathbb G}'(s):s\in\{s_{j_1},\ldots,s_{j_k}\}\}$ and $\{\overline{\mathbb G}(s):s\in\{s_{j_1},\ldots,s_{j_k}\}\}$ respectively, and let $\Sigma$ be the covariance matrix of $\{\overline{\mathbb G}(s):s\in\{s^o_1,\ldots,s^o_k\}\}$. Assumption \ref{continuous constraint} and the convergence of each $s_{j_i}$ to $s^o_i$ ensure that $\hat\Sigma\to \Sigma$ as $\epsilon_S\to 0$. To argue that $\hat \Sigma'-\hat \Sigma\to \mathbf{0}\in \R^{k\times k}$, where $\mathbf{0}$ denotes the $k\times k$ matrix with zero entries, we need the F-Glivenko-Contelli property of the product class $\mathcal F\cdot \mathcal F:=\{f(\cdot)=h(x,\cdot)h(x',\cdot):x,x'\in\mathcal X\}$. F-Donskerness implies F-Glivenko-Contelli, therefore $\mathcal F=\{h(x,\cdot):x\in\mathcal X\}$ is F-Glivenko-Contelli under Assumption \ref{Donsker h}, which together with Assumption \ref{moment2} forces the product class $\mathcal F\cdot \mathcal F$ to be F-Glivenko-Contelli by statement (ii) of Corollary 9.27 from \cite{kosorok2008introduction}. As a result, we have
\begin{equation*}
   \sup_{x,x'\in\mathcal X}\big\vert \frac{1}{n_2}\sum_{i=1}^{n_2}(h(x,\xi_i)-\hat H(x))(h(x',\xi_i)-\hat H(x'))-\mathrm{Cov}_F(h(x,\xi),h(x',\xi))\big\rvert\to 0\text{ as }n_2\to\infty
\end{equation*}
almost surely, where $\hat H(x)=(1/n_2)\sum_{i=1}^{n_2}h(x,\xi_i)$ and $\hat H(x')$ is similar. In particular $\hat \Sigma'-\hat \Sigma\to \mathbf{0}\in \R^{k\times k}$ as desired. Combining this with the convergence of $\hat \Sigma$ to $\Sigma$, we conclude that $\lim_{n_2\to\infty,\epsilon_S\to 0}\hat \Sigma'= \Sigma$ almost surely. Since the distribution of a zero mean multivariate Gaussian is uniquely determined by its covariance matrix, we must have $\{\overline{\mathbb G}'(s):s\in\{s_{j_1},\ldots,s_{j_k}\}\}$ weakly converges to $\{\overline{\mathbb G}(s):s\in\{s^o_1,\ldots,s^o_k\}\}$, and hence $\psi_{1-\beta}(\{\overline{\mathbb G}'(s):s\in\{s_{j_1},\ldots,s_{j_k}\}\})\to \psi_{1-\beta}(\{\overline{\mathbb G}(s):s\in\{s^o_1,\ldots,s^o_k\}\})$ almost surely. Note that $\psi_{1-\beta}(\{\overline{\mathbb G}'(s):s\in\{s_{j_1},\ldots,s_{j_k}\}\})\leq \psi_{1-\beta}(\{\overline{\mathbb G}'(s):s\in\{s_1,\ldots,s_p\}\})=q_{1-\beta}$, hence
\begin{equation*}
    \liminf_{n_2\to\infty,\epsilon_S\to 0}q_{1-\beta}\geq \psi_{1-\beta}(\{\overline{\mathbb G}(s):s\in\{s^o_1,\ldots,s^o_k\}\})\text{ for each }k.
\end{equation*}
This together with \eqref{separability of G} gives $\liminf_{n_2\to\infty,\epsilon_S\to 0}q_{1-\beta}\geq \bar{q}_{1-\beta}$. Altogether we have $\lim_{n_2\to\infty,\epsilon_S\to 0}q_{1-\beta}=\bar{q}_{1-\beta}$ almost surely.

The rest of the proof closely follows that of Theorem \ref{asymptotic joint:gaussian}. We only highlight some modifications. First, each occurrence of $\frac{z_{1-\beta}\hat\sigma_j}{\sqrt{n_2}}$, $\frac{z_{1-\beta}\hat\sigma(x^*(\hat s^*))}{\sqrt{n_2}}$ and $\frac{z_{1-\beta}\hat\sigma(x^*_S)}{\sqrt{n_2}}$ shall be replaced by $\frac{q_{1-\beta}}{\sqrt{n_2}}$. Second, the second error $E_2$ in \eqref{representation of constraint value} is no longer present, and the series of inequalities in the last paragraph become
\begin{eqnarray*}
&&\liminf_{n_2\to\infty,\epsilon_S\to 0}P_{\bm{\xi}_{1:n_2}}(H(x^*(\hat s^*))\geq \gamma)\\
&\geq &\liminf_{n_2\to\infty,\epsilon_S\to 0}P_{\bm{\xi}_{1:n_2}}\Big(-\frac{1}{\sqrt{n_2}}\mathbb G_{n_2}(x_S^*)+\frac{q_{1-\beta}}{\sqrt{n_2}}+o_p\big(\frac{1}{\sqrt{n_2}}\big)\geq 0\Big)\\
&=&\liminf_{n_2\to\infty,\epsilon_S\to 0}P_{\bm{\xi}_{1:n_2}}\Big(-\frac{1}{\sqrt{n_2}}\mathbb G_{n_2}(x_S^*)+\frac{\bar{q}_{1-\beta}}{\sqrt{n_2}}+o_p\big(\frac{1}{\sqrt{n_2}}\big)\geq 0\Big)\text{\ \ since }q_{1-\beta}\to \bar{q}_{1-\beta}\text{ a.s.}\\
&= &\liminf_{n_2\to\infty,\epsilon_S\to 0}P_{\bm{\xi}_{1:n_2}}\Big(\frac{\mathbb G_{n_2}(x_S^*)}{\sigma (x_S^*)}+o_p(1)\leq \frac{\bar{q}_{1-\beta}}{\sigma (x_S^*)}\Big)\\
&=&\Phi\big(\frac{\bar{q}_{1-\beta}}{\sigma (x_S^*)}\big)\text{\ \ by Slutsky's theorem}.
\end{eqnarray*}
This completes the proof for Algorithm \ref{algo:supremum_unnormalized}.

Now we prove the results for Algorithm \ref{algo:supremum_normalized} by a similar argument. Consider the Gaussian process
$$\{\tilde{\mathbb G}(s):s\in S^o\}$$
with mean zero and covariance structure $\mathrm{Cov}(s,s')=\mathrm{Cov}_F(h(x^*(s),\xi),h(x^*(s'),\xi))/(\sigma(x^*(s))\sigma(x^*(s')))$, and the Gaussian process
$$\{\tilde{\mathbb G}'(s):s\in S^o\}$$
with mean zero and covariance structure $\mathrm{Cov}(s,s')=\big[\frac{1}{n_2}\sum_{i=1}^{n_2}(h(x^*(s),\xi_i)-\hat H(x^*(s)))(h(x^*(s'),\xi_i)-\hat H(x^*(s')))\big]/(\hat\sigma(x^*(s))\hat\sigma(x^*(s')))$ where $\hat H(x^*(s))=(1/n_2)\sum_{i=1}^{n_2}h(x^*(s),\xi_i),\hat\sigma^2(x^*(s))=(1/n_2)\sum_{i=1}^{n_2}(h(x^*(s),\xi_i)-\hat H(x^*(s)))^2$ and $\hat H(x^*(s')),\hat\sigma^2(x^*(s'))$ are similarly defined. We have $\tilde{q}_{1-\beta}=\psi_{1-\beta}(\{\tilde{\mathbb G}(s):s\in S^o\})$, and $q_{1-\beta}=\psi_{1-\beta}(\{\tilde{\mathbb G}'(s):s\in \{s_1,\ldots,s_p\}\})\leq \psi_{1-\beta}(\{\tilde{\mathbb G}'(s):s\in S^o\})$. Under the depicted conditions, Lemma 3 from \cite{lam2016recovering} states that, for almost every realization of the data sequence $\{\xi_i\}_{i=1}^{\infty}$, the Gaussian process $\tilde{\mathbb G}'$ weakly converges to $\tilde{\mathbb G}$, so we have $\limsup_{n_2\to\infty,\epsilon_S\to 0}q_{1-\beta}\leq \limsup_{n_2\to\infty,\epsilon_S\to 0}\psi_{1-\beta}(\{\tilde{\mathbb G}'(s):s\in S^o\})= \tilde{q}_{1-\beta}$ almost surely. By a similar argument based on the separability of $\tilde{\mathbb G}$ and the uniform convergence of covariance as in the case of Algorithm \ref{algo:supremum_unnormalized}, we can show the other direction $\liminf_{n_2\to\infty,\epsilon_S\to 0}q_{1-\beta}\geq \tilde{q}_{1-\beta}$ and thereby conclude that $\lim_{n_2\to\infty,\epsilon_S\to 0}q_{1-\beta}= \tilde{q}_{1-\beta}$ almost surely. The rest of the proof for Algorithm \ref{algo:supremum_normalized} also follows that of Theorem \ref{asymptotic joint:gaussian}, but with each occurrence of $z_{1-\beta}$ replaced by $q_{1-\beta}$. The display in the last paragraph should be modified to be
\begin{eqnarray*}
&&\liminf_{n_2\to\infty,\epsilon_S\to 0}P_{\bm{\xi}_{1:n_2}}(H(x^*(\hat s^*))\geq \gamma)\\
&\geq &\liminf_{n_2\to\infty,\epsilon_S\to 0}P_{\bm{\xi}_{1:n_2}}\Big(-\frac{1}{\sqrt{n_2}}\mathbb G_{n_2}(x_S^*)+\frac{q_{1-\beta}\hat\sigma(x^*_S)}{\sqrt{n_2}}+o_p\big(\frac{1}{\sqrt{n_2}}\big)\geq 0\Big)\\
&=&\liminf_{n_2\to\infty,\epsilon_S\to 0}P_{\bm{\xi}_{1:n_2}}\Big(-\frac{1}{\sqrt{n_2}}\mathbb G_{n_2}(x_S^*)+\frac{\tilde{q}_{1-\beta}\sigma(x^*_S)}{\sqrt{n_2}}+o_p\big(\frac{1}{\sqrt{n_2}}\big)\geq 0\Big)\\
&&\text{\ \ since }q_{1-\beta}\to \tilde{q}_{1-\beta}\text{ and }\hat\sigma(x^*_S)\to\sigma(x^*_S)>0\text{ a.s.}\\
&= &\liminf_{n_2\to\infty,\epsilon_S\to 0}P_{\bm{\xi}_{1:n_2}}\Big(\frac{\mathbb G_{n_2}(x_S^*)}{\sigma (x_S^*)}+o_p(1)\leq \tilde{q}_{1-\beta}\Big)\\
&=&\Phi(\tilde{q}_{1-\beta})\text{\ \ by Slutsky's theorem}.
\end{eqnarray*}

Lastly, by stochastic dominance of the supremum of the Gaussian process $\overline{\mathbb G}$ or $\tilde{\mathbb G}$ over each of its marginal Gaussian component, it is straightforward that $\bar{q}_{1-\beta}\geq z_{1-\beta}\sigma(x^*_S)$ and that $\tilde{q}_{1-\beta}\geq z_{1-\beta}$, therefore both $\Phi\big(\frac{\bar{q}_{1-\beta}}{\sigma(x^*_S)}\big)$ and $\Phi(\tilde{q}_{1-\beta})$ are at least $1-\beta$.\hfill\Halmos

\section{Proofs of Results in Section \ref{sec:formulation showcase}}\label{sec:application proofs}
We first provide a lemma on the continuity of the solution path $x^*(s)$:
\begin{lemma}\label{continuity of solution path}
Suppose the formulation $OPT(s)$ satisfies Assumptions \ref{monotonicity} and \ref{convex decision space}-\ref{bounded level set}. If $\hat{\mathcal F}(s)=\{x:g_t(x,s)\leq 0,t=1,\ldots,T\}$ for some finite $T$ where each $g_t$ is jointly continuous in $x,s$ and convex in $x$ for every fixed $s$, the objective $f(x)$ is continuous, and $OPT(s):=\min\{f(x):x\in\mathcal X\cap \hat{\mathcal F}(s)\}$ has a unique solution $x^*(s)$ for all $s\in[\underline{s},\overline{s}]$, then the solution path $x^*(s)$ is continuous on $[\underline{s},\overline{s}]$.
\end{lemma}
\proof{Proof of Lemma \ref{continuity of solution path}.}The lemma is an application of Proposition 4.4 from \cite{bonnans2013perturbation}. Based on the discussion following Proposition 4.4 in \cite{bonnans2013perturbation}, we argue one by one that $OPT(s)$ satisfies assumptions (i)-(iv) of Proposition 4.4. Assumption (i): The objective $f(x)$ is continuous and independent of $s$ hence it's jointly continuous in $x,s$. Assumption (ii): The constraints of $OPT(s)$ can be formulated as $(f_1(x,s),\ldots,f_R(x,s),w'_1x-z_1,\ldots,w'_Lx-z_L,g_1(x,s),\ldots,g_T(x,s))\in [0,+\infty)^{R+L+T}$, where the left hand side is a vector of continuous functions and the right hand size is a closed convex cone. Assumption (iii) is implied by our Assumptions \ref{monotonicity} and \ref{bounded level set}. Assumption (iv): Since $OPT(s)$ is convex, our Assumptions \ref{monotonicity} and \ref{Slater's condition} ensure Slater's condition for $OPT(s)$ for all $s\in[s_l,s_u]$, and Slater's condition implies Robinson's constraint qualification, a sufficient condition for assumption (iv). Therefore the set-valued mapping $\mathcal X^*(s):=\{x\in \mathcal X\cup \hat{\mathcal F}(s):f(x)=\min_{x\in \mathcal X\cup \hat{\mathcal F}(s)}f(x)\}$ is upper semicontinuous at every $s\in [\underline{s},\overline{s}]$. When the optimal solution $x^*(s)$ for $OPT(s)$ is unique, upper semicontinuity implies continuity, hence $x^*(s)$ is continuous on $ [\underline{s},\overline{s}]$.\hfill\Halmos

The second lemma we present concerns the uniqueness of $x^*(s)$ for linear objectives:
\begin{lemma}\label{SCI uniqueness}
Consider an optimization problem in the form of $\min c'x$ subject to $f_r(x)\leq 0$ for $r=1,\ldots,R$ and $Ax\leq b$ for $A=[a_1,\ldots,a_L]'\in\R^{L\times d}$ and $b\in\R^L$ where each $f_r$ is continuous and convex and $c$ is a non-zero vector. For each $f_r$ and any two solutions $x_1\neq x_2$ such that $f_r(x_1)=f_r(x_2)=0$, assume $f_r(\theta x_1+(1-\theta)x_2)<0$ for any $\theta\in(0,1)$. If any $k\leq d-1$ rows of $A$ does not satisfy the SCI condition, then the optimal solution must be unique whenever one exists.
\end{lemma}
\proof{Proof of Lemma \ref{SCI uniqueness}.}Suppose there are two optimal solutions $x_1,x_2$. By convexity any solution in the form $\theta x_1+(1-\theta)x_2$ for $\theta\in[0,1]$ is also optimal, and because of the condition on $f_r$ we can assume that $f_r(x_1)<0,f_r(x_2)<0$ for all $r=1,\ldots,R$. Therefore, only the linear constraints can be binding on the line segment $\theta x_1+(1-\theta)x_2,\theta\in[0,1]$. Let $A_ox\leq b_o$ be the binding linear constraints on the segment where $A_o$ consists of rows of $A$ and $b_o$ contains the corresponding components of $b$, then it is clear that solution of the form $\theta x_1+(1-\theta)x_2$ is optimal for the linear program $\min c'x$ subject to $A_ox\leq b_o$. Since $A_ox_1=A_ox_2=b_o$, we have $A_o(x_2-x_1)=0$ hence the rank of $A_o$ is at most $d-1$. Now consider the dual $\min b'_oy$ subject to $A'_oy=-c,y\geq 0$. Since the rank of $A_o$ is at most $d-1$, by removing linearly dependent rows, the constraint $A'_oy=-c$ can be simplified to $\tilde A'_oy=-\tilde c$ where $\tilde A'_o$ has at most $d-1$ linearly independent rows. Let $y^*$ be an optimal basic feasible solution of the dual with the simplified constraint $\tilde A'_oy=-\tilde c$, then $y^*$ has at most $d-1$ non-zero (positive) components. However as a feasible solution $y^*$ has to satisfy $A'_oy^*=-c$ therefore the SCI condition holds for the rows of $A$ corresponding to the positive components of $y^*$, leading to a contradiction.\hfill\Halmos

\proof{Proof of Theorem \ref{continuity of SAA}.}We only need to verify the conditions of Lemma \ref{continuity of solution path}. In both cases (i) and (ii), $\hat{\mathcal F}(s)=\{x:\gamma+s-\frac{1}{n}\sum_{i=1}^nh(x,\xi_i)\leq 0\}$ and $\gamma+s-\frac{1}{n}\sum_{i=1}^nh(x,\xi_i)$ is obviously jointly continuous in $x,s$ and convex in $x$, and also $f(x)$ is continuous. Therefore, it only remains to check uniqueness of $x^*(s)$ in order to apply Lemma \ref{continuity of solution path}.

In case (i), the strict convexity of $f(x)$ forces the solution $x^*(s)$ to be unique. In case (ii), we first treat the case when $h$ is linear in $x$. We first note that for such $h$ the SAA takes the form $-\big(\frac{1}{n}\sum_{i=1}^nA(\xi_i)\big)'x\leq \frac{1}{n}\sum_{i=1}^nb(\xi_i)-\gamma-s$. Therefore each constraint of $OPT(s)$ is either linear or strictly convex, and thanks to Lemma \ref{SCI uniqueness} it remains to show that the SCI condition is not satisfied for each $s\in S$. For any $k\leq d-2$ rows $\{w_{l(1)},\ldots,w_{l(k)}\}$ of $W$ and the coefficient vector $-\frac{1}{n}\sum_{i=1}^nA(\xi_i)$, we want to show that the SCI condition does not hold for $\{w_{l(1)},\ldots,w_{l(k)},-\frac{1}{n}\sum_{i=1}^nA(\xi_i)\}$. Suppose SCI does hold, then we have the representation $-\frac{1}{n}\sum_{i=1}^nA(\xi_i) = \sum_{j=1}^k\lambda_j w_{l(j)}+\lambda_cc$, i.e., $-\frac{1}{n}\sum_{i=1}^nA(\xi_i)$ lies in the subspace of dimension spanned by $\{w_{l(1)},\ldots,w_{l(k)},c\}$. However, $-\frac{1}{n}\sum_{i=1}^nA(\xi_i)$ has a density hence lies in any given subspace of dimension $\leq d-1$ with probability zero. Therefore almost surely SCI does not hold for $\{w_{l(1)},\ldots,w_{l(k)},-\frac{1}{n}\sum_{i=1}^nA(\xi_i)\}$. If only linear coefficients from $Wx\leq z$ are considered, SCI condition is again not satisfied by the condition imposed. Therefore almost surely SCI is not satisfied for $OPT(s)$. By noting that the SCI condition is independent of $s$ since $s$ is on the right hand side, we conclude that almost surely SCI is not satisfied for all $s\in S$. When $h(x,\xi)$ is strictly concave in $x$, Lemma \ref{SCI uniqueness} can be directly applied to show the uniqueness of $x^*(s)$.\hfill\Halmos

\proof{Proof of Theorem \ref{continuity of phi DRO}.}We first argue that the constraint function $g(x,s)=\inf\Big\{\sum_{i=1}^nw_ih(x,\xi_i):\sum_{i=1}^n\frac{1}{n}\phi(nw_i)\leq s,\sum_{i=1}^nw_i=1,w_i\geq 0\text{ for all } i\Big\}$ is jointly continuous in $x,s$. Viewing both $x,s$ as parameters of the optimization problem defining $g(x,s)$, one can easily check that the assumptions of Proposition 4.4 from \cite{bonnans2013perturbation} are satisfied, hence $g(x,s)$ as the optimal value of the optimization problem is continuous in the parameters $x,s$. It is also obvious that $g(x,s)$ is concave in $x$ for every $s$ because of its representation as the minimum of a family of concave functions. By Lemma \ref{continuity of solution path} it remains to show the uniqueness of $x^*(s)$.

Case (i) follows from the strict convexity of $f$ as in Theorem \ref{continuity of SAA}. In case (ii), we would like to show that the constraint function $g(x,s)$ is strictly concave in $x$. Indeed, due to compactness an optimal weight vector $\mathbf w^*$ must exist for the minimization problem defining $g(x,s)$. Consider $x_1\neq x_2$ and $\theta\in(0,1)$, and let $\mathbf w^*$ be the minimizing weight vector that gives the worst-case value $g(\theta x_1+(1-\theta)x_2,s)$ at the solution $\theta x_1+(1-\theta)x_2$. Then because of the strict concavity of $h(x,\xi)$ in $x$, we have
\begin{eqnarray*}
g(\theta x_1+(1-\theta)x_2,s)&=&\sum_{i=1}^nw^*_ih(\theta x_1+(1-\theta)x_2,\xi_i)\\
&>&\sum_{i=1}^nw^*_i(\theta h(x_1,\xi_i)+(1-\theta)h(x_2,\xi_i))\\
&\geq&\theta g(x_1,s)+(1-\theta)g(x_2,s).
\end{eqnarray*}
Therefore $g(x,s)$ is strictly concave in $x$, and uniqueness of $x^*(s)$ follows from the SCI condition not being satisfied and applying Lemma \ref{SCI uniqueness}. In case (iii), the strict concavity of $g(x,s)$ can be shown as follows. Due to the strict convexity of $\phi$, for each decision $x$ the minimizing weight vector $\mathbf w^*$ not only exists but also is unique. Let $x_1\neq x_2$, then by the condition there must be some $\theta'\in [0,1)$ such that $\widehat{\mathrm{Corr}}(x_1,\theta' x_1+(1-\theta')x_2)\neq 1$, and let $\mathbf w^{1*}, \mathbf w^{\theta'*}$ be the respective minimizing weight vectors. We argue that $\mathbf w^{1*}\neq\mathbf w^{\theta'*}$. The condition $\lim_{x\to 0+}\phi(x)=+\infty$ ensures positive components of the minimizing weight vector, hence the optimality condition entails $h(x_1,\xi_i)=\lambda_1\phi'(nw_i^{1*})-\beta_1$ and $h(\theta' x_1+(1-\theta')x_2,\xi_i)=\lambda_{\theta'}\phi'(nw_i^{\theta'*})-\beta_{\theta'}$ for all $i=1,\ldots,n$ and some constants $\lambda_1,\beta_1,\lambda_{\theta'},\beta_{\theta'}$. Therefore the empirical correlation between $h(x_1,\xi)$ and $h(\theta'x_1+(1-\theta')x_2,\xi)$ takes the form
\begin{equation*}
\widehat{\mathrm{Corr}}(x_1,\theta' x_1+(1-\theta')x_2)=\frac{(1/n)\sum_{i=1}^n(\phi'(nw^{1*}_i)-\bar{\phi'_1})(\phi'(nw^{\theta'*}_i)-\bar{\phi'_{\theta'}})}{\sqrt{(1/n)\sum_{i=1}^n(\phi'(nw^{1*}_i)-\bar{\phi'_1})^2(1/n)\sum_{i=1}^n(\phi'(nw^{\theta'*}_i)-\bar{\phi'_{\theta'}})^2}}
\end{equation*}
where $\bar{\phi'_1}=(1/n)\sum_{i=1}^n\phi'(nw^{1*}_i),\bar{\phi'_{\theta'}}=(1/n)\sum_{i=1}^n\phi'(nw^{\theta' *}_i)$. If $\mathbf w^{1*}=\mathbf w^{\theta'*}$, we have $\widehat{\mathrm{Corr}}(x_1,\theta' x_1+(1-\theta')x_2)=1$, a contradiction. Therefore, if $\theta'>0$, we have $g(\theta' x_1+(1-\theta')x_2,s)=\sum_{i=1}^nw_1^{\theta'*}h(\theta' x_1+(1-\theta')x_2,\xi_i)\geq \theta'\sum_{i=1}^nw_1^{\theta'*}h(x_1,\xi_i)+(1-\theta')\sum_{i=1}^nw_1^{\theta'*}h(x_2,\xi_i)>\theta' g(x_1,s)+(1-\theta')g(x_2,s)$, hence $g(\theta x_1+(1-\theta)x_2,s)>\theta g(x_1,s)+(1-\theta)g(x_2,s)$ for all $\theta \in(0,1)$ by the (non-strict) concavity of $g(x,s)$ in $x$. Otherwise, if $\theta'=0$, i.e., $\widehat{\mathrm{Corr}}(x_1,x_2)\neq 1$, then by continuity there exists a small enough $\theta>0$ such that $\widehat{\mathrm{Corr}}(x_1,\theta x_1+(1-\theta)x_2)\neq 1$, hence things reduce to the previous case.\hfill\Halmos

\proof{Proof of Theorem \ref{continuity of Wasserstein DRO}.}Strong duality results from \cite{gao2016distributionally} or \cite{blanchet2016robust} show that the constraint function takes the form
\begin{equation*}
g(x,s):=\inf\Big\{\mathbb E_{G}[h(x,\xi)]:d_p(G,F_n )\leq s\Big\}=\sup_{\lambda\geq 0}\Big\{\frac{1}{n}\sum_{i=1}^n\inf_{\xi\in\Xi}(h(x,\xi)+\lambda \norm{\xi-\xi_i}^p)-\lambda s^p\Big\}.
\end{equation*}
We would like to show that $g(x,s)$ is jointly continuous in $x,s$. Let $a(x,\lambda)=\frac{1}{n}\sum_{i=1}^n\inf_{\xi\in\Xi}(h(x,\xi)+\lambda \norm{\xi-\xi_i}^p)$. It is clear that $a(x,\lambda)-\lambda s^p\leq \frac{1}{n}\sum_{i=1}^nh(x,\xi_i)-\lambda s^p\leq \frac{1}{n}\sum_{i=1}^nh(x,\xi_i)-\lambda s_l^p$ by taking $\xi=\xi_i$ in each infimum, and that $a(x,0)=\inf_{\xi\in\Xi}h(x,\xi)$. For each $x$ choose $\overline{\lambda}(x)$ so that $\frac{1}{n}\sum_{i=1}^nh(x,\xi_i)-\overline{\lambda}(x) s_l^p=\inf_{\xi\in\Xi}h(x,\xi)$. Since $h(x,\xi)$ is jointly continuous in $x,\xi$ and $\Xi$ is compact, $h(x,\xi)$ is uniformly continuous in $x,\xi$ on $[x_o-\delta,x_o+\delta]\times \Xi$ for given $x_o$ and $0<\delta<\infty$. This uniform continuity implies that the infimum $\inf_{\xi\in\Xi}h(x,\xi)$ is continuous in $x$ because $\lvert \inf_{\xi\in\Xi}h(x,\xi)-\inf_{\xi\in\Xi}h(x_o,\xi) \rvert\leq \sup_{\xi\in\Xi}\lvert h(x,\xi)-h(x_o,\xi)\rvert\to 0$ as $x\to x_o$. Therefore $\overline{\lambda}(x)$ is continuous in $x$ and
\begin{equation*}
g(x,s)=\sup_{0\leq\lambda\leq \overline{\lambda}(x)}(a(x,\lambda)-\lambda s^p)=\sup_{0\leq\lambda\leq \overline{\lambda}(x)}\Big\{\frac{1}{n}\sum_{i=1}^n\inf_{\xi\in\Xi}(h(x,\xi)+\lambda \norm{\xi-\xi_i}^p)-\lambda s^p\Big\}.
\end{equation*}
Since $\Xi$ is compact, by an argument similar to the one used to prove the continuity of $\inf_{\xi\in\Xi}h(x,\xi)$ we see that $a(x,\lambda)$ is jointly continuous in $x,\lambda$. Continuity of $a(x,\lambda)$ and $\overline{\lambda}(x)$ leads to the joint continuity of $g(x,s)$ in $x,s$. To explain, for a fixed $x$ and some $\delta>0$, define $\overline{\lambda}_{\delta}:=\sup_{x'\text{ s.t. }\norm{x'-x}_2\leq \delta}\overline{\lambda}(x')$, so for all $x',s'$ such that $\norm{x'-x}_2\leq \delta$ and $\lvert s'-s\rvert\leq \delta$ we have
\begin{eqnarray*}
\lvert g(x',s')-g(x,s)\rvert &=&\lvert \sup_{0\leq\lambda\leq \overline{\lambda}_{\delta}}(a(x',\lambda)-\lambda {s'}^p)-\sup_{0\leq\lambda\leq \overline{\lambda}_{\delta}}(a(x,\lambda)-\lambda s^p)\rvert \\
&\leq &\sup_{0\leq\lambda\leq \overline{\lambda}_{\delta}}\lvert a(x',\lambda)-a(x,\lambda)\rvert+\overline{\lambda}_{\delta} \lvert {s'}^p-s^p\rvert\to 0\text{ as }x'\to x,s'\to s
\end{eqnarray*}
where the limit holds because $a(x,\lambda)$ is uniformly continuous on the compact set $\{x':\norm{x'-x}\leq \delta\}\times  [0,\overline{\lambda}_{\delta}]$. Concavity of $g(x,s)$ in $x$ holds because for any probability measure $G$ the expectation $\mathbb E_G[h(x,\xi)]$ is concave in $x$ and the infimum operation perserves concavity.

In order to utilize Lemma \ref{continuity of solution path}, it remains to prove uniqueness of $x^*(s)$ for all $s\in S$. In case (i) uniqueness trivially follows from strict convexity of $f$. In case (ii), we first establish a result concerning the existence of the worst-case distribution:
\begin{lemma}\label{existence of worst-case distribution}
Under the same conditions of Theorem \ref{continuity of Wasserstein DRO}, if $g(x,s)=a(x,\lambda^*)-\lambda^* s^p$ for some $\lambda^*>0$, then there exists a distribution $G^*$ that belongs to the Wasserstein ball and that achieves the worst-case expectation, i.e., $g(x,s)=\mathbb E_{G^*}[h(x,\xi)]$.
\end{lemma}
\proof{Proof of Lemma \ref{existence of worst-case distribution}.}This is a direct consequence of Corollary 1 in \cite{gao2016distributionally}. Note that, since $\Xi$ is compact and $h(x,\xi)$ is continuous, for each decision $x$ the quantity $a(x,\lambda)$ is finite for all $\lambda\geq 0$. Corollary 1 from \cite{gao2016distributionally} then entails the existence of the worst-case distribution if there exists a dual maximizer $\lambda^*>0$.\hfill\Halmos

Consider $x_1\neq x_2$. If there exists some $\theta'\in (0,1)$ such that $g(\theta' x_1+(1-\theta')x_2,s)= a(\theta' x_1+(1-\theta')x_2,\lambda')-\lambda's^p$ for some $\lambda'>0$, then there exists some distribution $G_{\theta'}$ in the Wasserstein ball generating the worst-case expectation $g(\theta' x_1+(1-\theta')x_2,s)$. The strict concavity of $h$ then implies
\begin{eqnarray*}
g(\theta' x_1+(1-\theta')x_2,s)&=&\mathbb E_{G_{\theta'}}[h(\theta' x_1+(1-\theta')x_2,\xi)]\\
&>&\mathbb E_{G_{\theta'}}[\theta' h(x_1,\xi)+(1-\theta')h(x_2,\xi)]\\
&\geq &\theta'g(x_1,s)+(1-\theta')g(x_2,s).
\end{eqnarray*}
Since $g(x,s)$ is (non-strictly) concave in $x$, the above strict inequality at a certain $\theta'$ extends to all other $\theta$, i.e., $g(\theta x_1+(1-\theta)x_2,s)>\theta g(x_1,s)+(1-\theta) g(x_2,s)$ for all $\theta \in(0,1)$. Otherwise if $g(\theta x_1+(1-\theta)x_2,s)=a(\theta x_1+(1-\theta)x_2,0)>a(\theta x_1+(1-\theta)x_2,\lambda)-\lambda s^p$ for all $\theta\in(0,1)$ and $\lambda >0$, since $a(x,0)=\inf_{\xi\in\Xi}h(x,\xi)$ we still have the strict concavity of $g(\theta x_1+(1-\theta)x_2,s)$ in $\theta$. Therefore, according to Lemma \ref{SCI uniqueness}, the solution $x^*(s)$ is unique.\hfill\Halmos

\proof{Proof of Theorem \ref{continuity of moment DRO}.}We first show that each $g_i(x,s):=\sup_{(\mu,\Sigma)\in \mathcal U_i(s)}\mu'x+\sqrt{\frac{1-\alpha_i}{\alpha_i}}\sqrt{x'\Sigma x}$ is jointly continuous in $x$ and $s$. For a fixed pair $x_o,s_o$ and an arbitrary pair $x,s$, we write
\begin{eqnarray}
\notag&&\lvert g_i(x,s)-g_i(x_o,s_o)\rvert\\
\notag&\leq &\lvert g_i(x_o,s)-g_i(x_o,s_o)\rvert+\lvert g_i(x,s)-g_i(x_o,s)\rvert\\
\notag&\leq&\lvert g_i(x_o,s)-g_i(x_o,s_o)\rvert+\sup_{(\mu,\Sigma)\in \mathcal U_i(s)}\lvert \mu'x+\sqrt{\frac{1-\alpha_i}{\alpha_i}}\sqrt{x'\Sigma x}-(\mu'x_o+\sqrt{\frac{1-\alpha_i}{\alpha_i}}\sqrt{x_o'\Sigma x_o})\rvert\\
\notag&\leq&\lvert g_i(x_o,s)-g_i(x_o,s_o)\rvert+\sup_{(\mu,\Sigma)\in \mathcal U_i(s_u)}\lvert \mu'x+\sqrt{\frac{1-\alpha_i}{\alpha_i}}\sqrt{x'\Sigma x}-(\mu'x_o+\sqrt{\frac{1-\alpha_i}{\alpha_i}}\sqrt{x_o'\Sigma x_o})\rvert\\
&&\text{\ \ where $s_u$ is the maximal value for $s$}.\label{moment DRO errors}
\end{eqnarray}
Note that $\mu'x+\sqrt{\frac{1-\alpha_i}{\alpha_i}}\sqrt{x'\Sigma x}$ as a function jointly in $\mu,\Sigma,x$ is continuous, and hence by the compactness of $\mathcal U_i(s_u)$ is uniformly continuous for $(\mu,\Sigma)\in \mathcal U_i(s_u)$ and $x$ in some neighborhood of $x_o$. Uniform continuity implies that the second term in \eqref{moment DRO errors} vanishes as $x\to x_o$. It remains to show that the first term in \eqref{moment DRO errors} also vanishes, i.e., $g_i(x_o,s) \to g_i(x_o,s_o)$, as $s\to s_o$. We first show that as $s\to s_o$
\begin{align}
&\sup_{(\mu_s,\Sigma_s)\in\mathcal U_i(s)}\inf_{(\mu_{s_o},\Sigma_{s_o})\in\mathcal U_i(s_o)}(\norm{\mu_s-\mu_{s_o}}_2+\norm{\Sigma_s-\Sigma_{s_o}}_2)\to 0,\label{limit one}\\
&\sup_{(\mu_{s_o},\Sigma_{s_o})\in\mathcal U_i(s_o)}\inf_{(\mu_s,\Sigma_s)\in\mathcal U_i(s)}(\norm{\mu_s-\mu_{s_o}}_2+\norm{\Sigma_s-\Sigma_{s_o}}_2)\to 0.\label{limit two}
\end{align}
We prove \eqref{limit one} by contradiction. Suppose there exists $(\mu_{s_k},\Sigma_{s_k})\in\mathcal U_i(s_k)$ and $s_k\to s_o$ such that $\inf_{(\mu_{s_o},\Sigma_{s_o})\in\mathcal U_i(s_o)}(\norm{\mu_{s_k}-\mu_{s_o}}_2+\norm{\Sigma_{s_k}-\Sigma_{s_o}}_2)>\epsilon$ for some $\epsilon>0$. Note that all $(\mu_{s_k},\Sigma_{s_k})$ lie in the compact set $\mathcal U_i(s_u)$, hence there is a subsequence converging to some $(\mu_{\infty},\Sigma_{\infty})$ such that $\inf_{(\mu_{s_o},\Sigma_{s_o})\in\mathcal U_i(s_o)}(\norm{\mu_{\infty}-\mu_{s_o}}_2+\norm{\Sigma_{\infty}-\Sigma_{s_o}}_2)\geq\epsilon$, i.e., $(\mu_{\infty},\Sigma_{\infty})\notin \mathcal U_i(s_o)$. Since $\cap_{s>s_o}\mathcal U_i(s)=\mathcal U_i(s_o)$ and $\mathcal U_i(s)$ is non-decreasing in $s$, there exists some $\delta>0$ such that $(\mu_{\infty},\Sigma_{\infty})\notin\mathcal U_i(s)$ for all $s\leq s_o+\delta$, a contradiction with the convergence to $(\mu_{\infty},\Sigma_{\infty})$. To show \eqref{limit two}, suppose there exists $(\mu_k,\Sigma_k)\in\mathcal U_i(s_o)$ and $s_k\to s_o$ such that
\begin{equation}\label{limit three}
\inf_{(\mu_{s_k},\Sigma_{s_k})\in\mathcal U_i(s_k)}(\norm{\mu_k-\mu_{s_k}}_2+\norm{\Sigma_k-\Sigma_{s_k}}_2)>\epsilon
\end{equation}
for some $\epsilon>0$. By compactness, assume without loss of generality that $(\mu_k,\Sigma_k)$ converges to some limit $(\mu_{\infty},\Sigma_{\infty})\in \mathcal U_i(s_o)$. However, the condition $\overline{\cup_{s<s_o}\mathcal U_i(s)}=\mathcal U_i(s_o)$ ensures that $\inf_{(\mu_{s_k},\Sigma_{s_k})\in\mathcal U_i(s_k)}(\norm{\mu_{\infty}-\mu_{s_k}}_2+\norm{\Sigma_{\infty}-\Sigma_{s_k}}_2)\to 0$ as $s_k\to s_o$, which further entails that
\begin{eqnarray*}
&&\inf_{(\mu_{s_k},\Sigma_{s_k})\in\mathcal U_i(s_k)}(\norm{\mu_k-\mu_{s_k}}_2+\norm{\Sigma_k-\Sigma_{s_k}}_2)\\
&\leq &\inf_{(\mu_{s_k},\Sigma_{s_k})\in\mathcal U_i(s_k)}(\norm{\mu_{\infty}-\mu_{s_k}}_2+\norm{\Sigma_{\infty}-\Sigma_{s_k}}_2+\norm{\mu_{\infty}-\mu_k}_2+\norm{\Sigma_{\infty}-\Sigma_k}_2)\\
&\leq &\inf_{(\mu_{s_k},\Sigma_{s_k})\in\mathcal U_i(s_k)}(\norm{\mu_{\infty}-\mu_{s_k}}_2+\norm{\Sigma_{\infty}-\Sigma_{s_k}}_2)+\norm{\mu_{\infty}-\mu_k}_2+\norm{\Sigma_{\infty}-\Sigma_k}_2\\
&\to& 0
\end{eqnarray*}
a contradiction with \eqref{limit three}. This proves \eqref{limit two}. Now we use \eqref{limit one} and \eqref{limit two} to conclude $g_i(x_o,s) \to g_i(x_o,s_o)$ as $s\to s_o$. Since $\mathcal U_i(s_o)$ is compact, there exists an $(\mu^*_{s_0},\Sigma^*_{s_o})\in\mathcal U_i(s_o)$ such that $g_i(x_o,s_o)={\mu_{s_o}^*}'x_o+\sqrt{\frac{1-\alpha_i}{\alpha_i}}\sqrt{x_o'\Sigma^*_{s_o} x_o}$. \eqref{limit two} entails that there exists some $(\mu_s,\Sigma_s)\in \mathcal U_i(s)$ for each $s$ such that $(\mu_s,\Sigma_s)\to (\mu^*_{s_0},\Sigma^*_{s_o})$, therefore $\liminf_{s\to s_o}g_i(x_o,s)\geq \liminf_{s\to s_o}\mu_s'x_o+\sqrt{\frac{1-\alpha_i}{\alpha_i}}\sqrt{x_o'\Sigma_s x_o}=g_i(x_o,s_o)$. On the other hand, for each $s$, compactness of $\mathcal U_i(s)$ implies the existence of some $(\mu^*_s,\Sigma^*_s)\in \mathcal U_i(s)$ such that $g_i(x_o,s)={\mu_s^*}'x_o+\sqrt{\frac{1-\alpha_i}{\alpha_i}}\sqrt{x_o'\Sigma^*_s x_o}$. \eqref{limit one} then implies that there exists corresponding $(\mu^o_s,\Sigma^o_s)\in \mathcal U_i(s_o)$ such that $\norm{\mu^o_s-\mu^*_s}\to 0$ and $\norm{\Sigma^o_s-\Sigma^*_s}\to 0$ as $s\to s_o$. Since $\mu'x_o+\sqrt{\frac{1-\alpha_i}{\alpha_i}}\sqrt{x_o'\Sigma x_o}$ as a function of $(\mu,\Sigma)$ is uniformly continuous on $\mathcal U_i(s_u)$, we have $\limsup_{s\to s_o}g_i(x_o,s)=\limsup_{s\to s_o}{\mu_s^*}'x_o+\sqrt{\frac{1-\alpha_i}{\alpha_i}}\sqrt{x_o'\Sigma^*_s x_o}=\limsup_{s\to s_o}{\mu_s^o}'x_o+\sqrt{\frac{1-\alpha_i}{\alpha_i}}\sqrt{x_o'\Sigma^o_s x_o}\leq g_i(x_o,s_o)$. Altogether we have shown that $g_i(x_o,s) \to g_i(x_o,s_o)$, hence $g_i$ is jointly continuous in $x,s$.

Secondly, we show the uniqueness of $x^*(s)$ so that the desired result follows from applying Lemma \ref{continuity of solution path}. Note that the supremum of a family of convex functions is still convex, therefore each $g_i$ is convex in $x$. In case (i), strictly convexity of $f$ automatically forces uniqueness of $x^*(s)$. In case (ii), we prove uniqueness by either condition (3) or condition (4). Consider $x_1\neq x_2$ such that $g_i(x_1,s)=g_i(x_2,s)=b_i$, and $x_{\theta}:=(1-\theta)x_1+\theta x_2$ for some $\theta \in(0,1)$. Note that it is impossible that $x_1=cx_2$ or $x_2=cx_1$ for some $c\geq 0$ and $c\neq 1$, because otherwise $g_i(x_1,s)= cg_i(x_2,s)=cb_i\neq b_i$ or $g_i(x_2,s)= cg_i(x_1,s)=cb_i\neq b_i$. Let $(\mu_{\theta},\Sigma_{\theta})\in\mathcal U_i(s)$ be such that $g_i(x_{\theta},s)={\mu_{\theta}}'x_{\theta}+\sqrt{\frac{1-\alpha_i}{\alpha_i}}\sqrt{x_{\theta}'\Sigma_{\theta} x_{\theta}}$. Under condition (3), $\Sigma_{\theta}$ is automatically positive definite. Under condition (4), $\Sigma_{\theta}$ can be taken to be $\Sigma_s$ because $x_{\theta}'\Sigma_sx_{\theta}-x_{\theta}'\Sigma_{\theta}x_{\theta}=x_{\theta}'(\Sigma_s-\Sigma_{\theta})x_{\theta}\geq 0$, where the last inequality is due to $\Sigma_{\theta}\preceq \Sigma_s$. That is, in either case, $\Sigma_{\theta}$ can be taken to be positive definite. We then follow the proof of Theorem \ref{continuity of RO ellipsoidal} to show that
\begin{eqnarray*}
g_i(x_{\theta},s)&< &(1-\theta)\big({\mu_{\theta}}'x_1+\sqrt{\frac{1-\alpha_i}{\alpha_i}}\sqrt{x_1'\Sigma_{\theta} x_1}\big)+\theta\big({\mu_{\theta}}'x_2+\sqrt{\frac{1-\alpha_i}{\alpha_i}}\sqrt{x_2'\Sigma_{\theta} x_2}\big)\\
&\leq &(1-\theta)g_i(x_1,s)+\theta g_i(x_2,s)
\end{eqnarray*}
and to conclude uniqueness of $x^*(s)$ for each $s$ using Lemma \ref{SCI uniqueness}.\hfill\Halmos

\proof{Proof of Theorem \ref{continuity of RO polyhedral}.}We first transform the infinitely constrained robust counterpart into finitely many constraints. Note that, since each uncertainty set $\mathcal U_i(s)$ is a bounded polytope, in each robust constraint $\max_{a_i\in\mathcal U_i(s)}a'_ix\leq b_i$ the maximum is attained at a vertex of $\mathcal U_i(s)$. The set of vertices of $\mathcal U_i(s)$ takes the form
\begin{equation*}
    \mathcal V_i(s)=\bigg\{\widetilde{\mathcal W}_i^{-1}z_i+s\widetilde{\mathcal W}_i^{-1}e_i:
    \begin{array}{l}
         \widetilde{\mathcal W}_i\in\R^{d\times d}\text{ is an invertible submatrix of }\mathcal W_i  \\
         \mathcal W_i\widetilde{\mathcal W}_i^{-1}z_i-z_i\leq s(e_i-\mathcal W_i\widetilde{\mathcal W}_i^{-1}e_i)
    \end{array}\bigg\}
\end{equation*}
where the second condition ensures that $\widetilde{\mathcal W}_i^{-1}z_i+s\widetilde{\mathcal W}_i^{-1}e_i\in\mathcal U_i(s)$. The robust counterpart then becomes $v_i'x\leq b_i,v_i\in \mathcal V_i(s)$ for all $i=1,\ldots,K$. We make two important observations for $\mathcal V_i(s)$. First, the number of elements in $\mathcal V_i(s)$ is no more than the number of square submatrices of $\mathcal W_i$ which is finite. Second, the right hand side of $\mathcal W_i\widetilde{\mathcal W}_i^{-1}z_i-z_i\leq s(e_i-\mathcal W_i\widetilde{\mathcal W}_i^{-1}e_i)$ is linear in $s$ hence the system of inequalities are valid for $s$ in some interval of the form $(-\infty,u],[l,-\infty)$ or $[l,u]$, therefore the set of bases corresponding to vertices in $\mathcal V_i(s)$ changes at only finitely many $s$ values. That is, there are $s_l=s_0'\leq s_1'<\cdots<s_{q-1}'<s_q'= s_u$ such that, for each $1\leq j\leq q$, there exist submatrices $\widetilde{\mathcal W}_{i,1},\widetilde{\mathcal W}_{i,2},\ldots,\widetilde{\mathcal W}_{i,k_{i,j}}$ of each $\mathcal W_i$ such that the polyhedral RO $OPT(s)$ for all $s\in [s_{j-1}',s_j']$ can be simply expressed as
\begin{equation*}
\begin{aligned}
    &\min_{x\in\mathcal X}&&c'x\\
    &\text{subject to}&&(\widetilde{\mathcal W}_{i,l}^{-1}z_i+s\widetilde{\mathcal W}_{i,l}^{-1}e_i)'x\leq b_i\text{ for all }i=1,\ldots,K,l=1,\ldots,k_{i,j}.
\end{aligned}
\end{equation*}

The uniqueness of $x^*(s)$ is relatively straightforward to justify. The above representation of the RO and an application of Lemma \ref{SCI uniqueness} suggest that, under the imposed conditions regarding satisfaction of the SCI condition, the solution $x^*(s)$ can be non-unique at only finitely many $s$ values.

We now prove piecewise uniform continuity of the solution path. If at parameter value $\tilde s$ the solution $x^*(\tilde s)$ is not unique, we call it a \textit{non-unique point}. Between every two consecutive non-unique points $\tilde s_j<\tilde s_{j+1}$, $x^*(s)$ is unique hence is continuous in $(\tilde s_j,\tilde s_{j+1})$ due to Lemma \ref{continuity of solution path}. To show that $x^*(s)$ is actually uniformly continuous, it is sufficient and necessary to demonstrate that, as $s$ approaches some non-unique point $\tilde s$, left and right limits $\lim_{s\to \tilde s-}x^*(s),\lim_{s\to \tilde s+}x^*(s)$ exist. Without loss of generality, we focus on left limit. Toward this goal, we first derive a convenient formula of the optimal solution $x^*(s)$ for $s$ in a sufficiently small neighborhood $[\tilde s-\delta,\tilde s)$ of $\tilde s$. As shown in the first step, for sufficiently small $\delta$ the reformulation $OPT(s)$ takes the form
\begin{equation*}
\begin{aligned}
    &\min_x&&c'x\\
    &\text{subject to}&&(\widetilde{\mathcal W}_{i,l}^{-1}z_i+s\widetilde{\mathcal W}_{i,l}^{-1}e_i)'x\leq b_i\text{ for all }i=1,\ldots,K,l=1,\ldots,k_i\\
    &&&Wx\leq z
\end{aligned}
\end{equation*}
for all $s\in [\tilde s-\delta,\tilde s)$. For convenience, we rewrite the above parametric program in a more compact form
\begin{equation}\label{perturbation in matrix}
\begin{aligned}
    &\min_x&&c'x\\
    &\text{subject to}&&(A+s\Delta)x\leq b
\end{aligned}
\end{equation}
where the matrix $A$ contains all $\widetilde{\mathcal W}_{i,l}^{-1}z_i$'s and $W$ as its rows, and the right hand side $b$ has all the corresponding $b_i$'s and $z$ as its components, whereas the perturbation matrix $\Delta$ consists of all the $\widetilde{\mathcal W}_{i,l}^{-1}e_i$'s (and zero entries for the $W$ part of $A$). Note again that $x^*(s)$ is the unique optimal solution of \eqref{perturbation in matrix} for all $s\in[\tilde s-\delta,\tilde s)$. The dual of \eqref{perturbation in matrix} takes the form
\begin{equation}\label{perturbation in matrix:dual}
\begin{aligned}
    &\max_{y}&&b'y\\
    &\text{subject to}&&(A+s\Delta)'y= c\\
    &&&y\leq 0
\end{aligned}
\end{equation}
By the theory of simplex method, for the dual \eqref{perturbation in matrix:dual} there exists some basis $A_{\beta}+s\Delta_{\beta}$, where $\beta$ is a subset of size $d$ of $\{1,2,\ldots,\sum_{i=1}^Kk_i+L\}$ and $A_{\beta},\Delta_{\beta}$ denote the submatrices formed by the corresponding rows of $A,\Delta$, that gives rise to the optimal solution $y_{\beta}^*(s)=(A'_{\beta}+s\Delta'_{\beta})^{-1}c$ to \eqref{perturbation in matrix:dual} (other components of $y^*(s)$ are all zero). Moreover, the corresponding primal optimal solution to \eqref{perturbation in matrix} is $x^*(s)=(A_{\beta}+s\Delta_{\beta})^{-1}b_{\beta}$. By statement (ii) in Lemma 1 from \cite{freund1985postoptimal}, this optimal basis $\beta$ for \eqref{perturbation in matrix:dual} can change for only finitely many times as the parameter $s$ varies, therefore by choosing a small enough $\delta$ this basis $\beta$ remain the same one for all $s\in [\tilde s-\delta,\tilde s)$. That is, the unique optimal solution $x^*(s)=(A_{\beta}+s\Delta_{\beta})^{-1}b_{\beta}$ for all $s\in [\tilde s-\delta,\tilde s)$ and some basis $\beta$. Given this convenient formula, we now establish existence of the left limit. Case (i): $A_{\beta}+\tilde s\Delta_{\beta}$ is invertible. In this case the inverse $(A_{\beta}+s\Delta_{\beta})^{-1}$ must be continuous in $s$ at the non-unique point $\tilde s$, hence the left limit $\lim_{s\to \tilde s-}x^*(s)=(A_{\beta}+\tilde s\Delta_{\beta})^{-1}b_{\beta}$. Case (ii): $b_{\beta}$ is the zero vector. This case is trivial because $x^*(s)$ is also the zero vector hence the left limit exists and is the zero vector. Case (iii): $A_{\beta}+\tilde s\Delta_{\beta}$ is singular and $b_{\beta}$ is a non-zero vector. Note that Assumption \ref{bounded level set} implies that the solution path $\{x^*(s):s\in\mathcal S\}$ is confined within a bounded region, and we shall use this key information to conclude this case. For convenience we reparametrize the solution path as $s_o=(s-\tilde s+\delta)^{-1}$ and $x^*_o(s_o):=x^*(s)$ for $s\in(\tilde s -\delta,\tilde s)$. Letting $A_{\beta}^o=A_{\beta}+(\tilde s-\delta)\Delta_{\beta}$, we can express the reparametrization $x^*_o(s_o)$ as
$$x^*_o(s_o)=s_o\big((A_{\beta}^o)^{-1}\Delta_{\beta}+s_oI_d\big)^{-1}(A_{\beta}^o)^{-1}b_{\beta}\text{ for }s_o\in(\delta^{-1},+\infty)$$
and our goal is to show $\lim_{s_o\to \delta^{-1}+}x^*_o(s_o)$ exists. The matrix $(A_{\beta}^o)^{-1}\Delta_{\beta}$ admits a Jordan decomposition $(A_{\beta}^o)^{-1}\Delta_{\beta}=P^{-1}JP$, where $P$ is an invertible matrix with complex entries and $J$ is the Jordan normal form with the diagonal structure
\begin{equation*}
    J=\prth{\begin{matrix} 
      J_1 & && \\
       &J_2& & \\
       &&\ddots&\\
       &&&J_T
   \end{matrix}}, \text{ with each diagonal block }J_t=\prth{\begin{matrix} 
      \lambda_t &1 && &\\
       &\lambda_t&1 & &\\
       &&\ddots&&\\
       &&&\lambda_t&1\\
       &&&&\lambda_t
   \end{matrix}}
\end{equation*}
where each $\lambda_t$ is an eigenvalue of $(A_{\beta}^o)^{-1}\Delta_{\beta}$. With the Jordan decomposition, the reparametrized solution path takes the form
$$x^*_o(s_o)=s_oP^{-1}\big(J+s_oI_d\big)^{-1}P(A_{\beta}^o)^{-1}b_{\beta}.$$
Moreover, the inverse $\big(J+s_oI_d\big)^{-1}$ can be expressed as
\begin{equation*}
    \big(J+s_oI_d\big)^{-1}=\prth{\begin{matrix} 
      (J_1+s_oI)^{-1} & && \\
       &(J_2+s_oI)^{-1}& & \\
       &&\ddots&\\
       &&&(J_T+s_oI)^{-1}
   \end{matrix}}
\end{equation*}
where each diagonal block, if $J_t\in \R^{d_t\times d_t}$, has the form
\begin{equation}\label{inverse of Jordan form}
    (J_t+s_oI)^{-1}=\prth{\begin{matrix} 
      (\lambda_t+s_o)^{-1} &-(\lambda_t+s_o)^{-2} &\cdots& &(-1)^{d_t-1}(\lambda_t+s_o)^{-d_t}\\
       &(\lambda_t+s_o)^{-1}&-(\lambda_t+s_o)^{-2} &\cdots &\\
       &&\ddots&\ddots&\vdots\\
       &&&(\lambda_t+s_o)^{-1}&-(\lambda_t+s_o)^{-2}\\
       &&&&(\lambda_t+s_o)^{-1}
   \end{matrix}}.
\end{equation}
If we let $(P(A_{\beta}^o)^{-1}b_{\beta})_t$ be the vector of length $d_t$ consisting of the $(1+\sum_{i=1}^{t-1}d_i)$-th to $(\sum_{i=1}^td_i)$-th components of $P(A_{\beta}^o)^{-1}b_{\beta}$, then
$$x^*_o(s_o)=s_oP^{-1}\prth{\begin{matrix} 
      (J_1+s_oI)^{-1}(P(A_{\beta}^o)^{-1}b_{\beta})_1 \\
       (J_2+s_oI)^{-1}(P(A_{\beta}^o)^{-1}b_{\beta})_2\\
       \vdots\\
       (J_T+s_oI)^{-1}(P(A_{\beta}^o)^{-1}b_{\beta})_T
   \end{matrix}}.$$
We argue that $(P(A_{\beta}^o)^{-1}b_{\beta})_t$ must be the zero vector for all $t$ such that $\lambda_t=-\delta^{-1}$. Note that, since $A_{\beta}+\tilde s\Delta_{\beta}$ is singular, some $\lambda_t$ must be $-\delta^{-1}$. Consider a Jordan block $J_t$ with $\lambda_t=-\delta^{-1}$. From the form \eqref{inverse of Jordan form} of the inverse, one can check that $\norm{(J_t+s_oI)^{-1}v}_2\to\infty$ as $s_o\to\delta^{-1}$ for any given non-zero vector $v$. However, the solution $x^*(s)$, hence each $(J_t+s_oI)^{-1}(P(A_{\beta}^o)^{-1}b_{\beta})_t$, is confined to a bounded region, therefore $(P(A_{\beta}^o)^{-1}b_{\beta})_t$ must be zero if $\lambda_t=-\delta^{-1}$. For those blocks with $\lambda_t\neq -\delta^{-1}$, the inverse $(J_t+s_oI)^{-1}$ is continuous in $s_o$ at $s_o=\delta^{-1}$. Altogether, each block $(J_t+s_oI)^{-1}(P(A_{\beta}^o)^{-1}b_{\beta})_t$ is either constantly zero or continuous at $s_o=\delta^{-1}$, therefore $x^*_o(s_o)$ has right limit at $s_o=\delta^{-1}$. This proves the existence of left limit of $x^*(s)$ at $s=\tilde s$.\hfill\Halmos

\proof{Proof of Theorem \ref{continuity of RO ellipsoidal}.}The continuity of the second-order cone constraints in $x,s$ and its convexity in $x$ are straightforward. We only focus on the uniqueness of $x^*(s)$. In case (i) uniqueness trivially follows from strict convexity of $f$. In case (ii), we want to show for each cone constraint that for any $x_1,x_2$ such that $\mu'_ix_1+s\norm{\Sigma_ix_1}_2=\mu'_ix_2+s\norm{\Sigma_ix_2}_2=b_i$ we must have $\mu'_ix_{\theta}+s\norm{\Sigma_ix_{\theta}}_2<b_i$ for any $x_{\theta}=\theta x_1+(1-\theta)x_2$ where $\theta\in (0,1)$. First of all, there exists no $c\geq 0,c\neq 1$ such that $x_1=cx_2$ or $x_2=cx_1$ because otherwise $\mu'_ix_1+s\norm{\Sigma_ix_1}_2=c(\mu'_ix_2+s\norm{\Sigma_ix_2}_2)=cb_i\neq b_i$. Second, if there exists some $c<0$ such that $x_1=cx_2$ or $x_2=cx_1$, then $\mu'_ix_{\theta}+s\norm{\Sigma_ix_{\theta}}_2$ is piecewise linear in $\theta$ and has two pieces with different slopes, therefore $\mu'_ix_{\theta}+s\norm{\Sigma_ix_{\theta}}_2<b_i$ for all $\theta\in (0,1)$. Finally, if $x_1$ and $x_2$ are not parallel, then it is easy to verify that $\norm{\Sigma_ix_{\theta}}_2$ is strictly convex in $\theta$ by examining its second order derivative, therefore we have $\mu'_ix_{\theta}+s\norm{\Sigma_ix_{\theta}}_2<b_i$ again. Together with the SCI condition not being satisfied, we can use Lemma \ref{SCI uniqueness} to conclude the uniqueness of $x^*(s)$ for all $s\in S$. Lemma \ref{continuity of solution path} then implies the desired conclusion. \hfill\Halmos

\section{Finite Sample Performance Guarantees for Univariate Gaussian Validator}\label{sec:finite sample error univariate gaussian}
This section provides finite-sample errors regarding the performance guarantees presented in Theorem \ref{asymptotic joint:gaussian}, focusing on two general classes of constraints: differentiable stochastic constraints (Section \ref{sec:differentiable finite-sample}) and linear chance constraints (Section \ref{sec:CCP finite-sample}).

 \subsection{Differentiable Constraints}\label{sec:differentiable finite-sample}
 In order to derive finite-sample errors, we need stronger versions of Assumptions \ref{Donsker h}-\ref{continuous constraint} and \ref{unique optimal solution}. Assumption \ref{moment2} is replaced by boundedness of the fourth order moment:
 \begin{assumption}\label{moment4}
 $\mathfrak{m}_4:=\big(\mathbb E_F\big[\sup_{x\in\mathcal X}\abs{h(x,\xi)-H(x)}^4\big]\big)^{\frac{1}{4}}<\infty$.
 \end{assumption}

 The $L_2$-continuity condition for the constraint function $h$ in Assumption \ref{continuous constraint} is now strengthened to a differentiability condition:
 \begin{assumption}\label{continuously differentiable constraint}
 The random function $h(\cdot,\xi)$ is continuously differentiable on $\mathcal X$ for almost every $\xi\in\Xi$, and $\nabla H(x)=\mathbb E_F[\nabla h(x,\xi)]$. Assume $\overline{\rho}:=\sup_{x\in\mathcal X}\rho(\mathrm{Cov}_F(\nabla h(x,\xi)))<\infty$ where $\rho(\cdot)$ denotes the largest eigenvalue of a matrix.
 \end{assumption}
 Note that, in the presence of Assumption \ref{moment2}, Assumption \ref{continuously differentiable constraint} implies Assumption \ref{continuous constraint} through the dominated convergence theorem. When the gradient has a square integrable envelope, i.e., $\mathbb E_F[\sup_{x\in\mathcal X}\norm{\nabla h(x,\xi)}^2]<\infty$, and the decision space $\mathcal X$ is compact, Assumption \ref{continuously differentiable constraint} also implies Assumption \ref{Donsker h}.

 We then assume uniqueness of the optimal parameter, and local differentiability of the solution path and the expectation constraint:
 \begin{assumption}[Unique optimal parameter and local differentiability]\label{unique optimal parameter + local differentiability}
 The optimal parameter is unique, i.e., $S^*=\{s^*\}$, and $H(x_S^*)=\gamma$ at the optimal solution $x_S^*=x^*(s^*)$. Moreover, $H(x)$ is continuously differentiable in a neighborhood of $x_S^*$, and the parameter-to-solution mapping $x^*(s)$ is continuously differentiable in a neighborhood of $s^*$. There exists a $\delta>0$ such that for all $s\in[s^*-\delta,s^*+\delta]\subseteq S$ it holds $\frac{1}{2}\leq \nabla H(x^*(s))'\nabla x^*(s)/\nabla H(x^*_S)'\nabla x^*(s^*)\leq 2$ and $\norm{\nabla x^*(s)}_2/\norm{\nabla x^*(s^*)}_2\leq 2$, and that for all $s\leq s^*-\delta$ it holds $H(x^*(s))\leq H(x^*(s^*-\delta))$.
 \end{assumption}
We have the following finite-sample performance bounds for Algorithm \ref{algo:marginal}:
 \begin{theorem}[Finite-sample feasibility guarantee with univariate Gaussian validator]\label{finite sample error:gaussian+differentiable constraint}
 Suppose Assumptions \ref{continuous obj}, \ref{non-degenerate variance}-\ref{monotonicity} and \ref{moment4}-\ref{unique optimal parameter + local differentiability} hold, and $\{s_1,\ldots,s_p\}\cap \{\tilde s_1,\ldots,\tilde s_{M-1}\}=\emptyset$. Recall the mesh size $\epsilon_S = \sup_{s\in S}\inf_{1\leq j\leq p}\abs{s-s_j}$. Denote by $c^*:=\nabla H(x^*_S)'\nabla x^*(s^*)/\norm{\nabla x^*(s^*)}_2$, and by $C$ some universal constant. For any $t>0$ such that
\begin{equation*}
   2\epsilon_S<\mathrm{err}(p,n_2,t):=\frac{4(1+z_{1-\beta})\mathfrak{m}_4}{c^*\norm{\nabla x^*(s^*)}_2}\sqrt{\frac{t\log p}{n_2}}<\frac{\delta}{2}
\end{equation*}
it holds for the parameter $\hat s^*$ output by Algorithm \ref{algo:marginal} that
\begin{equation*}
P_{\bm{\xi}_{1:n_2}}(\abs{\hat s^*-s^*}>2\mathrm{err}(p,n_2,t))\leq \frac{C}{t}.
\end{equation*}
If
 \begin{equation*}
     2\epsilon_S<\frac{4(1+z_{1-\beta})\mathfrak{m}_4}{c^*\norm{\nabla x^*(s^*)}_2}\cdot \frac{(\log p)^{1/4}}{n_2^{3/8}}<\frac{\delta}{2}
 \end{equation*}
 we have
 \begin{equation}\label{feasibility error:gaussian+differentiable constraint}
    P_{\bm{\xi}_{1:n_2}}(x^*(\hat{s}^*)\text{ is feasible for \eqref{stoc_opt}})\geq 1-\beta-C(1+z_{1-\beta})^2\big(\frac{\mathfrak{m}_4}{\sigma(x^*_S)}\big)^3\big(1+\frac{\sqrt{\bar{\rho}}}{c^*}\big)^{\frac{2}{3}}\big(\frac{(\log p)^2}{n_2}\big)^{\frac{1}{4}}.
 \end{equation}
 \end{theorem}
\proof{Proof of Theorem \ref{finite sample error:gaussian+differentiable constraint}.}First we present a lemma concerning moment inequalities for the maximal deviation of sample means:
\begin{lemma}\label{maximal deviation:empirical process of finite cardinality}
Let $\mathcal G$ be function class of finite cardinality, and $G(\xi):=\max_{g\in\mathcal G}\lvert g(\xi)\rvert$ be the envelope function. Suppose $\xi_1,\ldots,\xi_n$ are i.i.d. observations from a common distribution $F$, then for any $k\geq 1$ we have
\begin{equation*}
    \sqrt{n}\prth{\mathbb E\big[\max_{g\in\mathcal G}\big\lvert\frac{1}{n}\sum_{i=1}^ng(\xi_i)-\mathbb E_F[g(\xi)]\big\rvert^k\big]}^{1/k}\leq C\sqrt{1+\log\lvert \mathcal G\rvert}\prth{\mathbb E_F[(G(\xi))^{\tilde k}]}^{1/\tilde k}
\end{equation*}
where $\tilde k=\max(2,k)$, the constant $C$ only depends on $k$, and $\lvert \mathcal G\rvert$ denotes the cardinality of $\mathcal G$.
\end{lemma}
\proof{Proof of Lemma \ref{maximal deviation:empirical process of finite cardinality}.}This is a direct consequence of Theorem 2.14.1 from \cite{van1996weak}. To apply that theorem, note that the covering number of the function class $\mathcal G$ is at most $\lvert \mathcal G\rvert$, hence its entropy integral is at most $\sqrt{1+\log\lvert \mathcal G\rvert}$.\hfill\Halmos

We use Lemma \ref{maximal deviation:empirical process of finite cardinality} to derive tail bounds for various maximal deviations. Denote by $H_j = H(x^*(s_j))$, and $\sigma^2_j=\sigma^2(x^*(s_j))$ for convenience. Applying Lemma \ref{maximal deviation:empirical process of finite cardinality} to $\{h(x^*(s_j),\cdot)-H_j:j=1,\ldots,p\}$ with $k=4$ gives
\begin{eqnarray*}
    n_2^2\mathbb E_{\bm{\xi}_{1:n_2}}\big[\big(\max_{j}\big\lvert\hat H_j-H_j\big\rvert\big)^4\big]&\leq& C(\log p)^2\mathbb E_F[(\max_{j}\lvert h(x^*(s_j),\xi)-H_j\rvert)^4]\\
    &\leq &C(\log p)^2\mathfrak{m}_4^4
\end{eqnarray*}
where $C$ is a universal constant (because $k$ is fixed at $4$) and $\mathbb E_{\bm{\xi}_{1:n_2}}$ denotes the expectation conditioned on Phase one data and with respect to Phase two data. Similarly applying the lemma to the squared class $\{(h(x^*(s_j),\cdot)-H_j)^2-\sigma^2_j:j=1,\ldots,p\}$ with $k=2$ gives
\begin{eqnarray*}
    &&n_2\mathbb E_{\bm{\xi}_{1:n_2}}\big[\big(\max_{j}\big\lvert\frac{1}{n_2}\sum_{i=1}^{n_2}(h(x^*(s_j),\xi_i)-H_j)^2-\sigma^2_j\big\rvert\big)^2\big]\\
    &\leq& C\log p\mathbb E_F[(\max_{j}\lvert (h(x^*(s_j),\xi)-H_j)^2-\sigma^2_j\rvert)^2]\\
    &\leq &C\log p\mathfrak{m}_4^4.
\end{eqnarray*}
By Markov's inequality, for any $t_1>0$ we have
\begin{equation*}
    \max_{j=1,\ldots,p}\abs{\hat H_j-H_j}\leq\frac{\mathfrak{m}_4t_1}{\sqrt{n_2}}
\end{equation*}
with probability at least $1-C(\log p)^2/t_1^4$ and
\begin{equation}\label{maximal deviation:sample variance}
    \max_{j=1,\ldots,p}\abs{\hat \sigma^2_j-\sigma^2_j}\leq \max_{j}\big\lvert\frac{1}{n_2}\sum_{i=1}^{n_2}(h(x^*(s_j),\xi_i)-H_j)^2-\sigma^2_j\big\rvert+\max_{j=1,\ldots,p}(\hat H_j-H_j)^2\leq\frac{\mathfrak{m}_4^2t_1}{\sqrt{n_2}}+\frac{\mathfrak{m}_4^2t_1^2}{n_2}
\end{equation}
with probability at least
\begin{equation*}
1-\frac{C(\log p)^2}{t_1^4}-\frac{C(\log p)}{t_1^2}.
\end{equation*}
Note that, when the upper bound \eqref{maximal deviation:sample variance} holds, $\max_j\hat \sigma^2_j\leq \max_{j}\sigma^2_j+\frac{\mathfrak{m}_4^2t_1}{\sqrt{n_2}}+\frac{\mathfrak{m}_4^2t_1^2}{n_2}\leq \mathfrak{m}_4^2+\frac{\mathfrak{m}_4^2t_1}{\sqrt{n_2}}+\frac{\mathfrak{m}_4^2t_1^2}{n_2}$. Therefore for any $t_1>0$
\begin{eqnarray}
    \notag\max_{j=1,\ldots,p}\abs{\hat H_j-z_{1-\beta}\frac{\hat\sigma_j}{\sqrt{n_2}}-H_j}&\leq& \frac{\mathfrak{m}_4t_1}{\sqrt{n_2}}+z_{1-\beta}\sqrt{\frac{\mathfrak{m}_4^2}{n_2}+\frac{\mathfrak{m}_4^2t_1}{n_2^{3/2}}+\frac{\mathfrak{m}_4^2t_1^2}{n_2^2}}\\
    &\leq& (1+z_{1-\beta})\mathfrak{m}_4\frac{1+t_1}{\sqrt{n_2}}\leq 2(1+z_{1-\beta})\mathfrak{m}_4\frac{t_1}{\sqrt{n_2}}\label{maximal deviation:sample mean adjusted by std}
\end{eqnarray}
for all $t_1\geq 1$ with probability at least
\begin{equation*}
1-\frac{C(\log p)^2}{t_1^4}-\frac{C(\log p)}{t_1^2}.
\end{equation*}
For every constant $\epsilon<\delta$, the solution path $x^*(s)$ is differentiable for $s\in [s^*-\epsilon,s^*+\epsilon]$. Therefore for any $s_j$ such that $\lvert s_j-s^*\rvert\leq \epsilon$, by differentiability we have
\begin{eqnarray*}
    \abs{[h(x^*(s_j),\xi)-H_j]-[h(x^*_S,\xi)-H(x^*_S))]}&=&\abs{\int_{s^*}^{s_j}(\nabla h(x^*(s),\xi)-\nabla H(x^*(s)))'\nabla x^*(s)ds}\\
    &\leq&\int_{s^*-\epsilon}^{s^*+\epsilon}\abs{(\nabla h(x^*(s),\xi)-\nabla H(x^*(s)))'\nabla x^*(s)}ds.
\end{eqnarray*}
The right hand side of the above inequality serves as an envelope function of the function class $\{[h(x^*(s_j),\xi)-H_j]-[h(x^*_S,\xi)-H(x^*_S))]:\lvert s_j-s^*\rvert\leq \epsilon\}$. Assumption \ref{continuously differentiable constraint} entails that $\mathbb E_F[\abs{(\nabla h(x^*(s),\xi)-\nabla H(x^*(s)))'\nabla x^*(s)}^2]\leq \bar{\rho}\norm{\nabla x^*(s)}_2^2$ for all $s$, therefore by Jensen's inequality (or Minkowski's integral inequality)
\begin{eqnarray*}
    &&\mathbb E_F\big[\big(\int_{s^*-\epsilon}^{s^*+\epsilon}\abs{(\nabla h(x^*(s),\xi)-\nabla H(x^*(s)))'\nabla x^*(s)}ds\big)^2\big]\\
    &\leq& \big(\int_{s^*-\epsilon}^{s^*+\epsilon}\sqrt{\mathbb E_F[\abs{(\nabla h(x^*(s),\xi)-\nabla H(x^*(s)))'\nabla x^*(s)}^2]}ds\big)^2\\
    &\leq&\bar{\rho}\big(\int_{s^*-\epsilon}^{s^*+\epsilon}\norm{\nabla x^*(s)}_2ds\big)^2\\
    &\leq &16\bar{\rho}\norm{\nabla x^*(s^*)}_2^2\epsilon^2
\end{eqnarray*}
an upper bound for the second moment of the envelope. Now applying Lemma \ref{maximal deviation:empirical process of finite cardinality} with $k=2$ to $\{[h(x^*(s_j),\xi)-H_j]-[h(x^*_S,\xi)-H(x^*_S))]:\lvert s_j-s^*\rvert\leq \epsilon)\}$, and noting that the cardinality does not exceed $p$, we have
\begin{equation*}
    n_2\mathbb E_{\bm{\xi}_{1:n_2}}\big[\max_{j:\lvert s_j-s^*\rvert\leq \epsilon}\big\lvert\hat H_j-H_j-(\hat H(x^*_S)-H(x^*_S))\big\rvert^2\big]\leq C(\log p) \bar{\rho}\norm{\nabla x^*(s^*)}_2^2\epsilon^2
\end{equation*}
which implies through Markov's inequality that for every $t_2>0$
\begin{equation}\label{maximal deviation:equicontinuity of empirical process}
    \max_{j:\lvert s_j-s^*\rvert\leq \epsilon}\big\lvert\hat H_j-H_j-(\hat H(x^*_S)-H(x^*_S))\big\rvert\leq \frac{t_2}{\sqrt{n_2}}
\end{equation}
with probability at least
\begin{equation*}
    1-\frac{C(\log p) \bar{\rho}\norm{\nabla x^*(s^*)}_2^2\epsilon^2}{t_2^2}.
\end{equation*}
Deviation inequalities \eqref{maximal deviation:sample variance}, \eqref{maximal deviation:sample mean adjusted by std} and \eqref{maximal deviation:equicontinuity of empirical process} are the key elements for establishing finite sample error bounds. Lastly, we also need a bound characterizing the modulus of continuity of the variance $\sigma^2(x^*(s))$. For every $s_j$ such that $\lvert s_j-s^*\rvert\leq \epsilon$
\begin{eqnarray*}
\lvert \sigma^2_j-\sigma^2(x^*_S)\rvert&=&\lvert\mathbb E_F[(h(x^*(s_j),\xi)-H_j)^2]-\mathbb E_F[(h(x^*_S,\xi)-H(x^*_S))^2] \rvert\\
&=&\abs{\mathbb E_F\big[ \int_{s^*}^{s_j}2(h(x^*(s),\xi)-H(x^*(s)))(\nabla h(x^*(s),\xi)-\nabla H(x^*(s)))'\nabla x^*(s)ds\big]}\\
&\leq &\mathbb E_F\big[ \int_{s^*}^{s_j}2\vert h(x^*(s),\xi)-H(x^*(s))\rvert \lvert(\nabla h(x^*(s),\xi)-\nabla H(x^*(s)))'\nabla x^*(s)\rvert ds\big]\\
&\leq &\mathbb E_F\big[ \int_{s^*-\epsilon}^{s^*+\epsilon}2\lvert h(x^*(s),\xi)-H(x^*(s))\rvert \lvert(\nabla h(x^*(s),\xi)-\nabla H(x^*(s)))'\nabla x^*(s)\rvert ds\big]\\
&=& \int_{s^*-\epsilon}^{s^*+\epsilon}2\mathbb E_F\big[\lvert h(x^*(s),\xi)-H(x^*(s))\rvert \lvert(\nabla h(x^*(s),\xi)-\nabla H(x^*(s)))'\nabla x^*(s)\rvert\big] ds\\
&&\text{\ \ \ by Fubini's theorem}\\
&\leq&\int_{s^*-\epsilon}^{s^*+\epsilon}2\sigma(x^*(s)) \sqrt{\bar{\rho}}\norm{\nabla x^*(s)}_2 ds\text{\ \ by Cauchy Schwartz inequality}\\
&\leq& 8\mathfrak{m}_4 \sqrt{\bar{\rho}}\norm{\nabla x^*(s^*)}_2\epsilon.
\end{eqnarray*}
That is, for all $\epsilon<\delta$
\begin{equation}\label{deviation:variance}
    \max_{j:\lvert s_j-s^*\rvert\leq \epsilon}\lvert \sigma_j^2-\sigma^2(x^*_S)\rvert\leq 8\mathfrak{m}_4 \sqrt{\bar{\rho}}\norm{\nabla x^*(s^*)}_2\epsilon.
\end{equation}

We first show the deviation inequality for $\hat s^*$. If \eqref{maximal deviation:sample mean adjusted by std} happens, and $t_1$ is such that
\begin{equation}\label{choice of t1,p,n2}
   2\epsilon_S<\epsilon(t_1,n_2):=\frac{4(1+z_{1-\beta})\mathfrak{m}_4}{c^*\norm{\nabla x^*(s^*)}_2}\cdot\frac{t_1}{\sqrt{n_2}}<\frac{\delta}{2} 
\end{equation}
we want to show that $\lvert \hat s^*-s^*\rvert\leq 2\epsilon(t_1,n_2)$. By Assumption \ref{unique optimal parameter + local differentiability}, for any $s\in(s^*,s^*+\delta]$ the constraint value $H(x^*(s))\geq \gamma +\frac{s-s^*}{2}\nabla H(x^*(s^*))'\nabla x^*(s^*)=\gamma +\frac{s-s^*}{2}c^*\norm{\nabla x^*(s^*)}_2$, and similarly $H(x^*(s))\leq \gamma +\frac{s-s^*}{2}c^*\norm{\nabla x^*(s^*)}_2$ for all $s\in[s^*-\delta,s^*)$. Therefore $H_j> \gamma+2(1+z_{1-\beta})\mathfrak{m}_4\frac{t_1}{\sqrt{n_2}}$ for all $s_j\in (s^*+\epsilon(t_1,n_2),s^*+\delta]$ and $H_j< \gamma-2(1+z_{1-\beta})\mathfrak{m}_4\frac{t_1}{\sqrt{n_2}}$ for all $s_j\in [s^*-\delta,s^*-\epsilon(t_1,n_2))$. Under the condition that $2\epsilon_S<\epsilon(t_1,n_2)$ there must be some $\overline{j}$ for which $s_{\overline{j}}\in (s^*+\epsilon(t_1,n_2),s^*+2\epsilon(t_1,n_2))\subset (s^*+\epsilon(t_1,n_2),s^*+\delta]$ and hence $\hat H_{\overline{j}}-z_{1-\beta}\frac{\hat\sigma_{\overline{j}}}{\sqrt{n_2}}>\gamma$ on one hand. On the other hand the solution path has a derivative $\nabla x^*(s)$ that is non-zero in $[s^*-\delta,s^*+\delta]$ hence the parameter-to-objective mapping $v(s)$ strictly increases in $s$ in the same interval. Therefore the picked parameter $\hat s^*\leq s_{\overline{j}}$. Similarly, there exists some $\underline{j}$ such that  $s_{\underline{j}}\in (s^*-2\epsilon(t_1,n_2),s^*-\epsilon(t_1,n_2))\subset (s^*-\delta,s-\epsilon(t_1,n_2)]$ and $H_{\underline{j}}< \gamma-2(1+z_{1-\beta})\mathfrak{m}_4\frac{t_1}{\sqrt{n_2}}$. Since $H(x^*(s^*-\delta))\geq H(x^*(s))$ for all $s\leq s^*-\delta$, we have for all $s_j\leq s_{\underline{j}}$ that $H_j\leq H_{\underline{j}}$ and $\hat H_j-z_{1-\beta}\frac{\hat\sigma_j}{\sqrt{n_2}}<\gamma$, therefore $\hat s^*\leq s_{\underline{j}}$ is impossible. That is, it must be the case that $\hat s^*\in (s_{\underline{j}}, s_{\overline{j}}]\subset [s^*-2\epsilon(t_1,n_2),s^*+2\epsilon(t_1,n_2)]$. This gives the deviation inequality
\begin{equation*}
    P_{\bm{\xi}_{1:n_2}}(\lvert \hat s^*-s^*\rvert> 2\epsilon(t_1,n_2))\leq C\big(\frac{(\log p)^2}{t_1^4}+\frac{\log p}{t_1^2}\big)\leq \frac{C\log p}{t_1^2}
\end{equation*}
provided that \eqref{choice of t1,p,n2} holds. Since the above bound becomes trivial when $(\log p)/t_1^2\geq 1$, hence we can assume $(\log p)/t_1^2< 1$ without loss of generality (and enlarge the universal constant $C$ if necessary) to get
\begin{equation}\label{deviation:parameter}
    P_{\bm{\xi}_{1:n_2}}(\lvert \hat s^*-s^*\rvert> 2\epsilon(t_1,n_2))\leq \frac{C\log p}{t_1^2}.
\end{equation}

Now we derive the finite sample error for the feasibility confidence level. Using the same notation $\epsilon(t_1,n_2)$, we write
\begin{eqnarray*}
&&P_{\bm{\xi}_{1:n_2}}(H(x^*(\hat s^*))\geq \gamma)\\
&\geq&P_{\bm{\xi}_{1:n_2}}(H(x^*(\hat s^*))\geq \gamma,\lvert \hat s^*-s^*\rvert\leq 2\epsilon(t_1,n_2))\\
&\geq&P_{\bm{\xi}_{1:n_2}}(\hat H(x^*(\hat s^*))-H(x^*(\hat s^*))-\frac{z_{1-\beta}\hat \sigma(x^*(\hat s^*))}{\sqrt{n_2}}\leq \hat H(x^*(\hat s^*))-\frac{z_{1-\beta}\hat \sigma(x^*(\hat s^*))}{\sqrt{n_2}}- \gamma,\lvert \hat s^*-s^*\rvert\leq 2\epsilon(t_1,n_2))\\
&\geq&P_{\bm{\xi}_{1:n_2}}(\hat H(x^*(\hat s^*))-H(x^*(\hat s^*))-\frac{z_{1-\beta}\hat \sigma(x^*(\hat s^*))}{\sqrt{n_2}}\leq 0,\lvert \hat s^*-s^*\rvert\leq 2\epsilon(t_1,n_2))\\
&=&P_{\bm{\xi}_{1:n_2}}(\hat H(x^*_S)-H(x^*_S)-\frac{z_{1-\beta}\sigma(x^*_S)}{\sqrt{n_2}}+\Delta_H+\Delta_{\sigma}\leq 0,\lvert \hat s^*-s^*\rvert\leq 2\epsilon(t_1,n_2))\\
&&\text{\ where }\Delta_H = (\hat H(x^*(\hat s^*))-H(x^*(\hat s^*)))-(\hat H(x^*_S)-H(x^*_S)),\ \Delta_{\sigma}=(z_{1-\beta}/\sqrt{n_2})(\sigma(x^*_S)-\hat \sigma(x^*(\hat s^*)))\\
&\geq&P_{\bm{\xi}_{1:n_2}}\big(\hat H(x^*_S)-H(x^*_S)-\frac{z_{1-\beta}\sigma(x^*_S)}{\sqrt{n_2}}+\max_{j:\lvert s_j-s^*\rvert\leq 2\epsilon(t_1,n_2)}\big\lvert\hat H_j-H_j-(\hat H(x^*_S)-H(x^*_S))\big\rvert+\\
&&\hspace{9ex}\frac{z_{1-\beta}}{\sqrt{n_2}}\max_{j:\lvert s_j-s^*\rvert\leq 2\epsilon(t_1,n_2)}\lvert \sigma(x^*_S)-\hat \sigma_j\rvert\leq 0,\ \lvert \hat s^*-s^*\rvert\leq 2\epsilon(t_1,n_2)\big)\\
&\geq&P_{\bm{\xi}_{1:n_2}}\big(\frac{\sqrt{n_2}(\hat H(x^*_S)-H(x^*_S))}{\sigma(x^*_S)}+\frac{\sqrt{n_2}}{\sigma(x^*_S)}\max_{j:\lvert s_j-s^*\rvert\leq 2\epsilon(t_1,n_2)}\big\lvert\hat H_j-H_j-(\hat H(x^*_S)-H(x^*_S))\big\rvert+\\
&&\hspace{9ex}\frac{z_{1-\beta}}{\sigma(x^*_S)}\max_{j:\lvert s_j-s^*\rvert\leq 2\epsilon(t_1,n_2)}\lvert \sigma(x^*_S)-\hat \sigma_j\rvert\leq z_{1-\beta},\ \lvert \hat s^*-s^*\rvert\leq 2\epsilon(t_1,n_2)\big).
\end{eqnarray*}
It follows from \eqref{maximal deviation:sample variance} and \eqref{deviation:variance} that
\begin{equation*}
    \max_{j:\lvert s_j-s^*\rvert\leq 2\epsilon(t_1,n_2)}\lvert \sigma^2(x^*_S)-\hat \sigma^2_j\rvert\leq \frac{\mathfrak{m}_4^2t_1}{\sqrt{n_2}}+\frac{\mathfrak{m}_4^2t_1^2}{n_2}+16\mathfrak{m}_4 \sqrt{\bar{\rho}}\norm{\nabla x^*(s^*)}_2\epsilon(t_1,n_2)
\end{equation*}
with probability at least $1-C(\log p)/t_1^2-C(\log p)^2/t_1^4$. If $\frac{\mathfrak{m}_4^2t_1}{\sqrt{n_2}}+\frac{\mathfrak{m}_4^2t_1^2}{n_2}+16\mathfrak{m}_4 \sqrt{\bar{\rho}}\norm{\nabla x^*(s^*)}_2\epsilon(t_1,n_2)\leq \sigma^2(x^*_S)/4$, it follows from mean value theorem that with at least the same probability
\begin{equation}\label{deviation:std local}
    \max_{j:\lvert s_j-s^*\rvert\leq 2\epsilon(t_1,n_2)}\lvert \sigma(x^*_S)- \hat \sigma_j\rvert\leq \frac{1}{\sigma(x^*_S)}\big(\frac{\mathfrak{m}_4^2t_1}{\sqrt{n_2}}+\frac{\mathfrak{m}_4^2t_1^2}{n_2}+16\mathfrak{m}_4 \sqrt{\bar{\rho}}\norm{\nabla x^*(s^*)}_2\epsilon(t_1,n_2)\big).
\end{equation}
For now we assume $\frac{\mathfrak{m}_4^2t_1}{\sqrt{n_2}}+\frac{\mathfrak{m}_4^2t_1^2}{n_2}+16\mathfrak{m}_4 \sqrt{\bar{\rho}}\norm{\nabla x^*(s^*)}_2\epsilon(t_1,n_2)\leq \sigma^2(x^*_S)/4$ holds so that the bound \eqref{deviation:std local} is valid. Later on we shall show that this is without loss of generality. We proceed as
\begin{eqnarray*}
&&P_{\bm{\xi}_{1:n_2}}(H(x^*(\hat s^*))\geq \gamma)\\
&\geq&P_{\bm{\xi}_{1:n_2}}\big(\frac{\sqrt{n_2}(\hat H(x^*_S)-H(x^*_S))}{\sigma(x^*_S)}+\frac{t_2}{\sigma(x^*_S)}+\frac{z_{1-\beta}}{\sigma^2(x^*_S)}\big(\frac{\mathfrak{m}_4^2t_1}{\sqrt{n_2}}+\frac{\mathfrak{m}_4^2t_1^2}{n_2}+16\mathfrak{m}_4 \sqrt{\bar{\rho}}\norm{\nabla x^*(s^*)}_2\epsilon(t_1,n_2)\big)\leq z_{1-\beta}\big)\\
&&\hspace{5ex}-P_{\bm{\xi}_{1:n_2}}\big(\max_{j:\lvert s_j-s^*\rvert\leq 2\epsilon(t_1,n_2)}\big\lvert\hat H_j-H_j-(\hat H(x^*_S)-H(x^*_S))\big\rvert>\frac{t_2}{\sqrt{n_2}}\big)-P_{\bm{\xi}_{1:n_2}}\big(\lvert \hat s^*-s^*\rvert> 2\epsilon(t_1,n_2)\big)\\
&&\hspace{5ex}-P_{\bm{\xi}_{1:n_2}}\big(\max_{j:\lvert s_j-s^*\rvert\leq 2\epsilon(t_1,n_2)}\lvert \sigma(x^*_S)- \hat \sigma_j\rvert> \frac{1}{\sigma(x^*_S)}\big(\frac{\mathfrak{m}_4^2t_1}{\sqrt{n_2}}+\frac{\mathfrak{m}_4^2t_1^2}{n_2}+16\mathfrak{m}_4 \sqrt{\bar{\rho}}\norm{\nabla x^*(s^*)}_2\epsilon(t_1,n_2)\big)\big)\\
&\geq&P_{\bm{\xi}_{1:n_2}}\big(\frac{\sqrt{n_2}(\hat H(x^*_S)-H(x^*_S))}{\sigma(x^*_S)}\leq z_{1-\beta}-\frac{t_2}{\sigma(x^*_S)}-\frac{z_{1-\beta}}{\sigma^2(x^*_S)}\big(\frac{\mathfrak{m}_4^2t_1}{\sqrt{n_2}}+\frac{\mathfrak{m}_4^2t_1^2}{n_2}+16\mathfrak{m}_4 \sqrt{\bar{\rho}}\norm{\nabla x^*(s^*)}_2\epsilon(t_1,n_2)\big)\big)\\
&&\hspace{5ex}-\frac{C(\log p) \bar{\rho}\norm{\nabla x^*(s^*)}_2^2(\epsilon(t_1,n_2))^2}{t_2^2}-\frac{C\log p}{t_1^2}\\
&&\text{\ \ \ by \eqref{maximal deviation:equicontinuity of empirical process}, \eqref{deviation:parameter} and \eqref{deviation:std local}}.
\end{eqnarray*}
To deal with the first probability term, we recall the Berry-Esseen theorem. There exists some universal constant $C$ such that
\begin{equation*}
   \sup_{t\in\R}\lvert P_{\bm{\xi}_{1:n_2}}\big(\frac{\sqrt{n_2}(\hat H(x^*_S)-H(x^*_S))}{\sigma(x^*_S)}\leq t\big)- \Phi(t)\rvert\leq \frac{C\mathbb E_F[\lvert h(x^*_S,\xi)-H(x^*_S)\rvert^3]}{\sigma^3(x^*_S)\sqrt{n_2}}
\end{equation*}
where $\Phi$ is the cumulative distribution function for the standard normal. Noting that $\mathbb E_F[\lvert h(x^*_S,\xi)-H(x^*_S)\rvert^3]\leq \mathfrak{m}_4^3$ and that $\Phi$ has a bounded derivative, we further bound the confidence level as
\begin{eqnarray*}
&&P_{\bm{\xi}_{1:n_2}}(H(x^*(\hat s^*))\geq \gamma)\\
&\geq&\Phi\big( z_{1-\beta}-\frac{t_2}{\sigma(x^*_S)}-\frac{z_{1-\beta}}{\sigma^2(x^*_S)}\big(\frac{\mathfrak{m}_4^2t_1}{\sqrt{n_2}}+\frac{\mathfrak{m}_4^2t_1^2}{n_2}+16\mathfrak{m}_4 \sqrt{\bar{\rho}}\norm{\nabla x^*(s^*)}_2\epsilon(t_1,n_2)\big)\big)-\frac{C\mathfrak{m}_4^3}{\sigma^3(x^*_S)\sqrt{n_2}}\\
&&\hspace{5ex}-\frac{C(\log p) \bar{\rho}\norm{\nabla x^*(s^*)}_2^2(\epsilon(t_1,n_2))^2}{t_2^2}-\frac{C\log p}{t_1^2}\\
&\geq& 1-\beta-C\big(\frac{t_2}{\sigma(x^*_S)}+\frac{z_{1-\beta}}{\sigma^2(x^*_S)}\big(\frac{\mathfrak{m}_4^2t_1}{\sqrt{n_2}}+\frac{\mathfrak{m}_4^2t_1^2}{n_2}+16\mathfrak{m}_4 \sqrt{\bar{\rho}}\norm{\nabla x^*(s^*)}_2\epsilon(t_1,n_2)\big)\big)-\frac{C\mathfrak{m}_4^3}{\sigma^3(x^*_S)\sqrt{n_2}}\\
&&\hspace{5ex}-\frac{C(\log p) \bar{\rho}\norm{\nabla x^*(s^*)}_2^2(\epsilon(t_1,n_2))^2}{t_2^2}-\frac{C\log p}{t_1^2}.
\end{eqnarray*}
Arranging terms gives
\begin{eqnarray}
\notag&&1-\beta-P_{\bm{\xi}_{1:n_2}}(H(x^*(\hat s^*))\geq \gamma)\\
\notag&\leq &C\Big(\frac{t_2}{\sigma(x^*_S)}+\frac{z_{1-\beta}}{\sigma^2(x^*_S)}\big(\frac{\mathfrak{m}_4^2t_1}{\sqrt{n_2}}+\frac{\mathfrak{m}_4^2t_1^2}{n_2}+\mathfrak{m}_4 \sqrt{\bar{\rho}}\norm{\nabla x^*(s^*)}_2\epsilon(t_1,n_2)\big)+\frac{\mathfrak{m}_4^3}{\sigma^3(x^*_S)\sqrt{n_2}}\\
\notag&&\hspace{5ex}+\frac{(\log p) \bar{\rho}\norm{\nabla x^*(s^*)}_2^2(\epsilon(t_1,n_2))^2}{t_2^2}+\frac{\log p}{t_1^2}\Big)\\
&\leq &C\Big(\frac{z_{1-\beta}}{\sigma^2(x^*_S)}\big(\frac{\mathfrak{m}_4^2t_1}{\sqrt{n_2}}+\frac{\mathfrak{m}_4^2t_1^2}{n_2}+\mathfrak{m}_4 \sqrt{\bar{\rho}}\norm{\nabla x^*(s^*)}_2\epsilon(t_1,n_2)\big)+\frac{\mathfrak{m}_4^3}{\sigma^3(x^*_S)\sqrt{n_2}}\label{finite sample error:gsc1}\\
\notag&&\hspace{5ex}+\frac{\big((\log p) \bar{\rho}\big)^{1/3}\norm{\nabla x^*(s^*)}_2^{2/3}(\epsilon(t_1,n_2))^{2/3}}{\sigma^{2/3}(x^*_S)}+\frac{\log p}{t_1^2}\Big)\\
\notag&&\text{\ \ by minimizing the bound over $t_2$}\\
\notag&\leq &C\Big(\frac{z_{1-\beta}}{\sigma^2(x^*_S)}\big(\frac{\mathfrak{m}_4^2t_1}{\sqrt{n_2}}+\mathfrak{m}_4 \sqrt{\bar{\rho}}\norm{\nabla x^*(s^*)}_2\epsilon(t_1,n_2)\big)+\frac{\mathfrak{m}_4^3}{\sigma^3(x^*_S)\sqrt{n_2}}\\
\notag&&\hspace{5ex}+\frac{\big((\log p) \bar{\rho}\big)^{1/3}\norm{\nabla x^*(s^*)}_2^{2/3}(\epsilon(t_1,n_2))^{2/3}}{\sigma^{2/3}(x^*_S)}+\frac{\log p}{t_1^2}\Big)
\end{eqnarray}
where in the last inequality we leave out the terms $\frac{\mathfrak{m}_4^2t_1^2}{n_2}$ because when $\frac{t_1}{\sqrt{n_2}}\leq 1$ it holds that $\frac{t_1^2}{n_2}\leq \frac{t_1}{\sqrt{n_2}}$ hence the former can be absorbed into the latter. Previously we assume that $\frac{\mathfrak{m}_4^2t_1}{\sqrt{n_2}}+\frac{\mathfrak{m}_4^2t_1^2}{n_2}+16\mathfrak{m}_4 \sqrt{\bar{\rho}}\norm{\nabla x^*(s^*)}_2\epsilon(t_1,n_2)\leq \sigma^2(x^*_S)/4$. This is without loss of generality, because otherwise the first error term in \eqref{finite sample error:gsc1} is of constant order which makes the upper bound trivial. Now expanding the $\epsilon(t_1,n_2)$ we further bound the error as follows
\begin{eqnarray}
\notag&&1-\beta-P_{\bm{\xi}_{1:n_2}}(H(x^*(\hat s^*))\geq \gamma)\\
\notag&\leq &C\Big(\frac{z_{1-\beta}\mathfrak{m}_4^2}{\sigma^2(x^*_S)}\big(\frac{t_1}{\sqrt{n_2}}+\Big[\frac{\bar{\rho}(1+z_{1-\beta})^2}{{c^*}^2n_2}\Big]^{1/2}t_1\big)+\frac{\mathfrak{m}_4^3}{\sigma^3(x^*_S)\sqrt{n_2}}+\Big[\frac{(\log p)\bar{\rho}(1+z_{1-\beta})^2\mathfrak{m}_4^2}{\sigma^2(x^*_S){c^*}^2n_2}\Big]^{1/3}t_1^{2/3}+\frac{\log p}{t_1^2}\Big)\\
\notag&\leq &C\Big(\frac{z_{1-\beta}\mathfrak{m}_4^2}{\sigma^2(x^*_S)}\big(1+\frac{\sqrt{\bar{\rho}}(1+z_{1-\beta})}{c^*}\big)\frac{t_1}{\sqrt{n_2}}+(\log p)^{1/3}\Big[\frac{\mathfrak{m}_4}{\sigma(x^*_S)}\big(1+\frac{\sqrt{\bar{\rho}}(1+z_{1-\beta})}{c^*}\big)\frac{t_1}{\sqrt{n_2}}\Big]^{2/3}\\
\notag&&\hspace{5ex}+\frac{\mathfrak{m}_4^3}{\sigma^3(x^*_S)\sqrt{n_2}}+\frac{\log p}{t_1^2}\Big)\\
\notag&\leq &C\Big(\big(\frac{z_{1-\beta}\mathfrak{m}_4}{\sigma(x^*_S)}+(\log p)^{1/3}\big)\Big[\frac{\mathfrak{m}_4}{\sigma(x^*_S)}\big(1+\frac{\sqrt{\bar{\rho}}(1+z_{1-\beta})}{c^*}\big)\frac{t_1}{\sqrt{n_2}}\Big]^{2/3}+\frac{\mathfrak{m}_4^3}{\sigma^3(x^*_S)\sqrt{n_2}}+\frac{\log p}{t_1^2}\Big)\\
\notag&&\text{\ \ since it can be assumed $\frac{\mathfrak{m}_4}{\sigma(x^*_S)}\big(1+\frac{\sqrt{\bar{\rho}}(1+z_{1-\beta})}{c^*}\big)\frac{t_1}{\sqrt{n_2}}\leq 1$}\\
\notag&\leq &C\Big(\frac{(1+z_{1-\beta})^{5/3}(\log p)^{1/3}\mathfrak{m}_4^3}{\sigma^3(x^*_S)n_2^{1/3}}\big(1+\frac{\sqrt{\bar{\rho}}}{c^*}\big)^{2/3}t_1^{2/3}+\frac{\log p}{t_1^2}+\frac{\mathfrak{m}_4^3}{\sigma^3(x^*_S)\sqrt{n_2}}\Big)\text{ \ \ since }\mathfrak{m}_4\geq \sigma(x^*_S)\\
&\leq &C\Big(\frac{(1+z_{1-\beta})^{5/3}(\log p)^{1/3}\mathfrak{m}_4^3}{\sigma^3(x^*_S)n_2^{1/3}}\big(1+\frac{\sqrt{\bar{\rho}}}{c^*}\big)^{2/3}t_1^{2/3}+\frac{\log p}{t_1^2}\Big)\label{finite sample error:gsc final form}
\end{eqnarray}
where in the last inequality we drop the last term since it's dominated by the first when $t_1\geq 1$. Note that \eqref{finite sample error:gsc final form} holds only under the condition \eqref{choice of t1,p,n2}. It is straightforward to see that the bound \eqref{finite sample error:gsc final form} is minimized at
\begin{equation*}
    t_1^*:=\frac{(\log p)^{1/4}\sigma^{9/8}(x^*_S)n_2^{1/8}}{(1+z_{1-\beta})^{5/8}\mathfrak{m}_4^{9/8}(1+\sqrt{\bar{\rho}}/c^*)^{1/4}}
\end{equation*}
by equating the two error terms. Consider $\tilde t_1:=(\log p)^{1/4}n_2^{1/8}$. Since $\tilde t_1=t^*_1(1+z_{1-\beta})^{5/8}\big(\frac{\mathfrak{m}_4}{\sigma(x^*_S)}\big)^{9/8}(1+\sqrt{\bar{\rho}}/c^*)^{1/4}>t^*_1$, the first term dominates at $t_1=\tilde t_1$. Therefore when \eqref{choice of t1,p,n2} is satisfied at $t_1=\tilde t_1$, we have
\begin{equation*}
    1-\beta-P_{\bm{\xi}_{1:n_2}}(H(x^*(\hat s^*))\geq \gamma)\leq C(1+z_{1-\beta})^{\frac{5}{3}}\big(\frac{\mathfrak{m}_4}{\sigma(x^*_S)}\big)^3\big(1+\frac{\sqrt{\bar{\rho}}}{c^*}\big)^{\frac{2}{3}}\big(\frac{(\log p)^2}{n_2}\big)^{\frac{1}{4}}.
\end{equation*}
The desired bound is obtained by replacing $\frac{5}{3}$ with $2$ as the exponent of $1+z_{1-\beta}$.\hfill\Halmos


\subsection{Linear Chance Constraints}\label{sec:CCP finite-sample}
Consider linear chance constraints in the form of $\mathbb P_F(a'_kx\leq b_k\text{ for }k=1,\ldots,K)\geq 1-\alpha$. We assume the following isotropy condition:
 \begin{assumption}[Isotropy]\label{isotropy}
 There exist constants $D_2,D_3$ such that for all unit vector $\nu\in\R^d$ and all $a_k,1\leq k\leq K$, the random variable $a'_k\nu$ has a sub-Gaussian norm at most $D_2$, i.e., $\mathbb E\big[\exp\big(\big(\frac{a'_kx}{D_2}\big)^2\big)\big]\leq 2$, and has a density bounded above by $D_3$. Each $b_k$ is a non-zero constant.
 \end{assumption}
 This assumption stipulates that each $a_k$ has variability of constant order in all directions, and it trivially holds when each $a_k$ is standard Gaussian.

 We have the following finite-sample performance bounds for linear chance constraints:
 \begin{theorem}[Finite-sample chance constraint feasibility guarantee with univariate Gaussian validator]\label{finite sample error:gaussian+linear chance constraint}
 Consider \eqref{chance constraint} with a linear chance constraint $h(x,\xi)=\mathbf{1}(a'_kx\leq b_k\text{ for }k=1,\ldots,K)$ and $0<\alpha<\frac{1}{2}$. Suppose Assumptions \ref{continuous obj}, \ref{piecewise uniform continuous solution curve}-\ref{monotonicity}, and \ref{unique optimal parameter + local differentiability}-\ref{isotropy} hold, and $\{s_1,\ldots,s_p\}\cap \{\tilde s_1,\ldots,\tilde s_{M-1}\}=\emptyset$. Recall the notations $\epsilon_S$ and $c^*$ from Theorem \ref{finite sample error:gaussian+differentiable constraint}. For any $t>0$ such that
 \begin{equation*}
2\epsilon_S<\mathrm{err}(p,n_2,t):=\frac{6(1+z_{1-\beta})}{c^*\norm{\nabla x^*(s^*)}_2}\big(\sqrt{\frac{\alpha\log (4pt)}{n_2}}+\frac{\log (4pt)}{n_2}\big)<\frac{\delta}{2}
 \end{equation*}
 it holds for the parameter $\hat s^*$ output by Algorithm \ref{algo:marginal} that
 $$P_{\bm{\xi}_{1:n_2}}\big(\abs{\hat s^*-s^*}>2\mathrm{err}(p,n_2,t)\big)\leq \frac{1}{t}.
 $$
If
 \begin{equation*}
     2\epsilon_S<\frac{6(1+z_{1-\beta})}{c^*\norm{\nabla x^*(s^*)}_2}\big(\sqrt{\frac{\alpha\log (pn_2)}{n_2}}+\frac{\log (pn_2)}{n_2}\big)<\frac{\delta}{2}
 \end{equation*}
 we have
 \begin{equation}\label{feasibility error:gaussian+chance constraint}
    P_{\bm{\xi}_{1:n_2}}(x^*(\hat{s}^*)\text{ is feasible for \eqref{chance constraint}})\geq 1-\beta-C(1+z_{1-\beta})^2\big(1+\sqrt{\tilde C K}\big(\log(\max\big\{3,\frac{n_2}{\tilde C}\big\})\big)^{\frac{1}{4}}\big)\big(\frac{(\log(pn_2))^3}{\alpha n_2}\big)^{\frac{1}{4}}
 \end{equation}
 where $C$ is a universal constant and
 \begin{equation*}
 \tilde C=\frac{D_2^2D_3\sqrt{\log(2K/\alpha)}}{c^*\min_{1\leq k\leq K}\lvert b_k\rvert}.
 \end{equation*}
\end{theorem}
To get a sense of the effect of the dimension $d$ on the finite-sample error \eqref{feasibility error:gaussian+chance constraint}, suppose  that $D_2,D_3,K,\{b_k,k=1,\ldots,K\}$ are all numbers of constant order and we focus on the number $c^*$. The latter is the derivative of the satisfaction probability $P(x^*(s))$ with respect to the parameter $s$ when the solution path is reparameterized to move at a unit speed. Therefore a proxy for the finite-sample performance of Algorithm \ref{algo:marginal} is the sensitivity of the satisfaction probability along the direction of the solution path. The more sensitive it is, the better is the finite-sample performance. Note that this sensitivity does not explicitly depend on the dimension.

Here we provide the proof Theorem \ref{finite sample error:gaussian+linear chance constraint}:
\proof{Proof of Theorem \ref{finite sample error:gaussian+linear chance constraint}.}The proof follows the same line of argument as that of Theorem \ref{finite sample error:gaussian+differentiable constraint}, but uses a different set of deviation inequalities tailored to bounded random variables. To avoid repetition, we focus on the derivation of these deviation inequalities.

We need the following concentration inequalities for the sample mean and sample variance:
\begin{lemma}[Adapted from \cite{maurer2009empirical}]\label{Bennett's inequality}
Let $X_i,i=1,\ldots,n$ be i.i.d. $[0,1]$-valued random variables, $\sigma^2=\mathrm{Var}(X_1)$, and $\hat \sigma^2$ be the sample variance. Then we have for every $\epsilon\in (0,1)$ that
\begin{equation*}
    P\Big(\big\rvert\frac{1}{n}\sum_{i=1}^nX_i- \mathbb E[X_1]\big\rvert>\sqrt{\frac{2\sigma^2\log(2/\epsilon)}{n}}+\frac{\log(2/\epsilon)}{3n}\Big)\leq \epsilon
\end{equation*}
and
\begin{equation*}
    P\Big(\lvert \hat\sigma -\sigma\rvert>\sqrt{\frac{2\log(2/\epsilon)}{n-1}}\Big)\leq \epsilon.
\end{equation*}
\end{lemma}
\proof{Proof of Lemma \ref{Bennett's inequality}.}Theorem 3 in \cite{maurer2009empirical} gives the following Bennett's inequality
\begin{equation*}
    P\Big(\frac{1}{n}\sum_{i=1}^nX_i< \mathbb E[X_1]-\big(\sqrt{\frac{2\sigma^2\log(1/\epsilon)}{n}}+\frac{\log(1/\epsilon)}{3n}\big)\Big)\leq \epsilon.
\end{equation*}
Applying the above inequality to $1-X_i,i=1,\ldots,n$ and noting that $\mathrm{Var}(1-X_1)=\mathrm{Var}(X_1)$, we have
\begin{equation*}
    P\Big(\frac{1}{n}\sum_{i=1}^nX_i> \mathbb E[X_1]+\big(\sqrt{\frac{2\sigma^2\log(1/\epsilon)}{n}}+\frac{\log(1/\epsilon)}{3n}\big)\Big)\leq \epsilon.
\end{equation*}
The first inequality in the lemma then comes from a union bound. The second inequality in the lemma is a direct consequence of Theorem 10 from \cite{maurer2009empirical}.\hfill\Halmos

Let $H_j:=\mathbb P_F(a'_kx^*(s_j)\leq b_k\text{ for all }k=1,\ldots,K)$ be the satisfaction probability at $x^*(s_j)$, and $\sigma^2_j:=H_j(1-H_j)$ be the variance. Applying Lemma \ref{Bennett's inequality} to each $\mathbf{1}(a'_kx^*(s_j)\leq b_k\text{ for all }k=1,\ldots,K)$ gives
\begin{equation*}
    \lvert \hat H_j-H_j\rvert\leq \sqrt{\frac{2H_j(1-H_j)\log(2/t_1)}{n_2}}+\frac{\log(2/t_1)}{3n_2}
\end{equation*}
with probability at least $1-t_1$, and
\begin{equation*}
    \lvert \hat\sigma_j-\sigma_j\rvert\leq \sqrt{\frac{2\log(2/t_1)}{n_2-1}}
\end{equation*}
with probability at least $1-t_1$. Using a union bound, we have
\begin{equation}\label{maximal deviation CCP:sample mean}
    \lvert \hat H_j-H_j\rvert\leq \sqrt{\frac{2H_j(1-H_j)\log(2p/t_1)}{n_2}}+\frac{\log(2p/t_1)}{3n_2}\text{ for all }j=1,\ldots,p
\end{equation}
with probability at least $1-t_1$, and that
\begin{equation}\label{maximal deviation CCP:std}
    \lvert \hat\sigma_j-\sigma_j\rvert\leq \sqrt{\frac{2\log(2p/t_1)}{n_2-1}}\text{ for all }j=1,\ldots,p
\end{equation}
with probability at least $1-t_1$. When \eqref{maximal deviation CCP:std} happen, we also have
\begin{equation*}
    \hat\sigma_j\leq \sigma_j+\lvert \hat\sigma_j-\sigma_j\rvert\leq \sqrt{H_j(1-H_j)}+\sqrt{\frac{2\log(2p/t_1)}{n_2-1}}\text{ for all }j=1,\ldots,p.
\end{equation*}
Together with \eqref{maximal deviation CCP:sample mean}, we can conclude that, with probability at least $1-2t_1$, for all $j=1,\ldots,p$
\begin{eqnarray}
    \notag\lvert \hat H_j-\frac{z_{1-\beta}\hat\sigma_j}{\sqrt{n_2}}-H_j\rvert&\leq& \frac{z_{1-\beta}}{\sqrt{n_2}} \big(\sqrt{H_j(1-H_j)}+\sqrt{\frac{2\log(2p/t_1)}{n_2-1}}\big)+\sqrt{\frac{2H_j(1-H_j)\log(2p/t_1)}{n_2}}+\frac{\log(2p/t_1)}{3n_2}\\
    \notag&\leq& 2\big(z_{1-\beta}+\sqrt{\log (2p/t_1)}\big)\big(\sqrt{\frac{H_j(1-H_j)}{n_2}}+\frac{\sqrt{\log (2p/t_1)}}{n_2}\big)\\
    &\leq& 2(1+z_{1-\beta})\sqrt{\log (2p/t_1)}\big(\sqrt{\frac{H_j(1-H_j)}{n_2}}+\frac{\sqrt{\log (2p/t_1)}}{n_2}\big)\label{maximal deviation CCP:sample mean adjusted by std}\\
    \notag&&\text{\ \ if we assume that }p\geq 2\text{ so that }\log (2p/t_1)>1.
\end{eqnarray}
Deviation bounds \eqref{maximal deviation CCP:std} and \eqref{maximal deviation CCP:sample mean adjusted by std} are CCP counterparts of \eqref{maximal deviation:sample variance} and \eqref{maximal deviation:sample mean adjusted by std}. Now we try to derive the CCP counterpart of \eqref{maximal deviation:equicontinuity of empirical process}. For any $\epsilon<\delta$ and every parameter value $s_j\in [s^*-\epsilon,s^*+\epsilon]$ we have by differentiability
\begin{equation*}
    \lvert a'_kx^*(s_j)-a'_kx^*_S\rvert =\lvert \int_{s^*}^{s_j}a'_k\nabla x^*(s)ds\rvert\leq \int_{s^*}^{s_j}\lvert a'_k\nabla x^*(s)\rvert ds\leq\eta_k(\epsilon):=\int_{s^*-\epsilon}^{s^*+\epsilon}\lvert a'_k\nabla x^*(s)\rvert ds.
\end{equation*}
Note that the sub-Gaussian norm $\norm{\cdot}_{\psi_2}:\{X\text{ is a random variable}:\norm{X}_{\psi_2}<\infty\}\to \R$ is a convex mapping, therefore by Jensen's inequality
\begin{equation*}
    \norm{\eta_k(\epsilon)}_{\psi_2}\leq \int_{s^*-\epsilon}^{s^*+\epsilon}\norm{a'_k\nabla x^*(s)}_{\psi_2} ds\leq \int_{s^*-\epsilon}^{s^*+\epsilon}D_2\norm{\nabla x^*(s)}_2 ds\leq 4D_2\norm{\nabla x^*(s^*)}_2\epsilon.
\end{equation*}
With the above bound of $\eta_k(\epsilon)$, we want to quantify the closeness of the linear chance constraint at the solutions $x^*(s_j)$ and $x^*_S$. We apply a union bound to obtain
\begin{eqnarray*}
    &&\lvert \mathbf{1}(a'_kx^*(s_j)\leq b_k\text{ for all }k=1,\ldots,K)-\mathbf{1}(a'_kx^*_S\leq b_k\text{ for all }k=1,\ldots,K) \rvert\\
    &\leq&\sum_{k=1}^K\lvert \mathbf{1}(a'_kx^*(s_j)\leq b_k)-\mathbf{1}(a'_kx^*_S\leq b_k) \rvert\\
    &\leq &\sum_{k=1}^K\mathbf{1}(a'_kx^*(s_j)\leq b_k<a'_kx^*_S\text{ or }a'_kx^*(s_j)> b_k\geq a'_kx^*_S)\\
    &\leq &\sum_{k=1}^K\mathbf{1}(a'_kx^*_S-\eta_k(\epsilon)\leq b_k<a'_kx^*_S\text{ or }a'_kx^*_S+\eta_k(\epsilon)> b_k\geq a'_kx^*_S)\\
    &\leq &\sum_{k=1}^K\mathbf{1}(\lvert a'_kx^*_S- b_k\rvert\leq \eta_k(\epsilon)).
\end{eqnarray*}
Noting that difference of two indicator functions takes values in $\{-1,0,1\}$, we have
\begin{eqnarray}
    \notag&&\mathbb E_{F}\big[\big(\mathbf{1}(a'_kx^*(s_j)\leq b_k\text{ for all }k=1,\ldots,K)-\mathbf{1}(a'_kx^*_S\leq b_k\text{ for all }k=1,\ldots,K)\big)^2\big]\\
    \notag&=&\mathbb E_{F}\big[\lvert\mathbf{1}(a'_kx^*(s_j)\leq b_k\text{ for all }k=1,\ldots,K)-\mathbf{1}(a'_kx^*_S\leq b_k\text{ for all }k=1,\ldots,K)\rvert\big]\\
    &\leq &\sum_{k=1}^K\mathbb P_F(\lvert a'_kx^*_S- b_k\rvert\leq \eta_k(\epsilon)).\label{local upper bound}
\end{eqnarray}
In order to derive an upper bound for each of the $K$ probabilities above, we first need a lower bound for $\norm{x^*_S}_2$. If there are some $\tilde k\in\{1,2,\ldots,K\}$ such that $b_{\tilde k}<0$, then
\begin{equation*}
    1-\alpha = \mathbb P_F(a'_kx^*_S\leq b_k\text{ for all }k=1,\ldots,K)\leq\mathbb P_F(a'_{\tilde k}x^*_S\leq b_{\tilde k})\leq 2\exp\big(-\frac{\min_k\lvert b_k\rvert^2}{D_2^2\norm{x^*_S}_2^2}\big)
\end{equation*}
where in the last inequality Assumption \ref{isotropy} is used. This forces $\norm{x^*_S}_2\geq \frac{\min_k\lvert b_k\rvert}{D_2\sqrt{\log(2/1-\alpha)}}$.
Otherwise if all $b_k>0$ then
\begin{equation*}
    \alpha = \mathbb P_F(\max_{k=1,\ldots,K}a'_kx^*_S- b_k>0)\leq \sum_{k=1}^K\mathbb P_F(a'_kx^*_S>b_k)\leq 2K\exp\big(-\frac{\min_k\lvert b_k\rvert^2}{D_2^2\norm{x^*_S}_2^2}\big)
\end{equation*}
which forces $\norm{x^*_S}_2\geq \frac{\min_k\lvert b_k\rvert}{D_2\sqrt{\log(2K/\alpha)}}$. When $\alpha >1/2$, the second lower bound dominates hence $\norm{x^*_S}_2\geq \frac{\min_k\lvert b_k\rvert}{D_2\sqrt{\log(2K/\alpha)}}$ always holds. Now we go back to \eqref{local upper bound} and notice that for each $k$ and every $c>0$
\begin{eqnarray*}
\mathbb P_F(\lvert a'_kx^*_S- b_k\rvert\leq \eta_k(\epsilon))&\leq&\mathbb P_F(\eta_k(\epsilon)>\varepsilon)+\mathbb P_F(\lvert a'_kx^*_S- b_k\rvert\leq \varepsilon)\\
&\leq& \mathbb P_F(\eta_k(\epsilon)>c)+\mathbb P_F(\lvert a'_kx^*_S- b_k\rvert\leq c)\\
&\leq& 2\exp\big( -\frac{c^2}{16D_2^2\norm{\nabla x^*(s^*)}_2^2\epsilon^2}\big)+\frac{2D_3c}{\norm{x^*_S}_2}\\
&\leq& 2\exp\big( -\frac{c^2}{16D_2^2\norm{\nabla x^*(s^*)}_2^2\epsilon^2}\big)+\frac{2D_2D_3\sqrt{\log(2K/\alpha)}c}{\min_k\lvert b_k\rvert}.
\end{eqnarray*}
With $c=4D_2\norm{\nabla x^*(s^*)}_2\epsilon\cdot \sqrt{\log\Big(\max\Big\{e,\frac{\min_k\lvert b_k\rvert}{D_2^2D_3\norm{\nabla x^*(s^*)}_2\sqrt{\log(2K/\alpha)}\epsilon}\Big\}\Big)}$, the above bound gives
\begin{equation*}
\mathbb P_F(\lvert a'_kx^*_S- b_k\rvert\leq \eta_k(\epsilon))\leq  10\tilde\epsilon\sqrt{\log\big(\max\{e,\frac{1}{\tilde \epsilon}\}\big)}
\end{equation*}
where $\tilde \epsilon:=\frac{D_2^2D_3\norm{\nabla x^*(s^*)}_2\sqrt{\log(2K/\alpha)}\epsilon}{\min_k\lvert b_k\rvert}$. From the union bound \eqref{local upper bound} it follows that
\begin{eqnarray}
    \notag&&\mathbb E_{F}\big[\big(\mathbf{1}(a'_kx^*(s_j)\leq b_k\text{ for all }k=1,\ldots,K)-\mathbf{1}(a'_kx^*_S\leq b_k\text{ for all }k=1,\ldots,K)\big)^2\big]\\
   &\leq& \sigma^2_{\epsilon}:=10K\tilde\epsilon\sqrt{\log\big(\max\{e,\frac{1}{\tilde \epsilon}\}\big)}\label{bound for local envelope CCP}
\end{eqnarray}
for all $s_j\in [s^*-\epsilon,s^*+\epsilon]$. In particular, $\sigma^2_{\epsilon}$ is a valid upper bound for the variance of each $\tilde h(x^*(s_j),\xi):=\mathbf{1}(a'_kx^*(s_j)\leq b_k\text{ for all }k=1,\ldots,K)-\mathbf{1}(a'_kx^*_S\leq b_k\text{ for all }k=1,\ldots,K)$ since the second moment always upper bounds the variance. Note that $(\tilde h(x^*(s_j),\xi)+1)/2$ is $[0,1]$-valued, hence applying Lemma \ref{Bennett's inequality} to $(\tilde h(x^*(s_j),\xi)+1)/2$ reveals that for all $s_j\in [s^*-\epsilon,s^*+\epsilon]$
\begin{equation*}
    \lvert\hat H_j-H_j-(\hat H(x^*_S)-H(x^*_S))\rvert\leq \sqrt{\frac{2\sigma^2_{\epsilon}\log(2/t_2)}{n_2}}+\frac{2\log(2/t_2)}{3n_2}
\end{equation*}
with probability at least $1-t_2$, therefore
\begin{equation}\label{equicontinuity of local empirical process CCP}
    \max_{j:\lvert s_j-s^*\rvert\leq \epsilon}\lvert\hat H_j-H_j-(\hat H(x^*_S)-H(x^*_S))\rvert\leq \sqrt{\frac{2\sigma^2_{\epsilon}\log(2p/t_2)}{n_2}}+\frac{2\log(2p/t_2)}{3n_2}
\end{equation}
with probability at least $1-t_2$, a counterpart of \eqref{maximal deviation:equicontinuity of empirical process}.

As in the proof of Theorem \ref{finite sample error:gaussian+differentiable constraint}, we first derive the deviation inequality for $\hat s^*$. Based on \eqref{maximal deviation CCP:sample mean adjusted by std}, we would like to find $\underline{H}$ such that, with high probability, for all $H_j\leq \underline{H}$ we have $\hat H_j-\frac{z_{1-\beta}\hat\sigma_j}{\sqrt{n_2}}<1-\alpha$, as well as $\overline{H}$ such that all $H_j\geq \overline{H}$ satisfies $\hat H_j-\frac{z_{1-\beta}\hat\sigma_j}{\sqrt{n_2}}>1-\alpha$. Given the bound \eqref{maximal deviation CCP:sample mean adjusted by std}, it suffices to $\underline{H}$ and $\overline{H}$ such that
\begin{eqnarray}
    H-2(1+z_{1-\beta})\sqrt{\log (2p/t_1)}\big(\sqrt{\frac{H(1-H)}{n_2}}+\frac{\sqrt{\log (2p/t_1)}}{n_2}\big)&>&1-\alpha,\text{ for all }H\geq \overline{H}\label{large H}\\
    H+2(1+z_{1-\beta})\sqrt{\log (2p/t_1)}\big(\sqrt{\frac{H(1-H)}{n_2}}+\frac{\sqrt{\log (2p/t_1)}}{n_2}\big)&<&1-\alpha,\text{ for all }H\leq \underline{H}.\label{small H}
\end{eqnarray}
For \eqref{large H}, since we must have $\overline{H}>1-\alpha$, it holds that $H(1-H)< \alpha(1-\alpha)<\alpha$ for all $H\geq \overline{H}$. Therefore $\overline{H}:=1-\alpha+2(1+z_{1-\beta})\sqrt{\log (2p/t_1)}\big(\sqrt{\frac{\alpha}{n_2}}+\frac{\sqrt{\log (2p/t_1)}}{n_2}\big)$ satisfies \eqref{large H}. For \eqref{small H}, since the left hand side is monotonic in $H$, we only need to find a $\underline{H}$ for which the inequality in \eqref{small H} holds true. If
\begin{equation}\label{condition for lower bound of H}
(1+z_{1-\beta})\sqrt{\log (2p/t_1)}\big(\sqrt{\frac{\alpha}{n_2}}+\frac{\sqrt{\log (2p/t_1)}}{n_2}\big)<\frac{\sqrt{2}}{4}\alpha
\end{equation}
then one can verify that $\underline{H}:=1-\alpha-2\sqrt{2}(1+z_{1-\beta})\sqrt{\log (2p/t_1)}\big(\sqrt{\frac{\alpha}{n_2}}+\frac{\sqrt{\log (2p/t_1)}}{n_2}\big)$ satisfies \eqref{small H} by noting that $\underline{H}>1-2\alpha$ and hence $\underline{H}(1-\underline{H})< 2\alpha$. In order for \eqref{condition for lower bound of H} to hold, we consider $p,t_1,n_2$ satisfying the following counterpart of \eqref{choice of t1,p,n2}
\begin{equation}\label{choice of p,t1,n2 CCP}
    2\epsilon_s<\epsilon(p,t_1,n_2):=\frac{6(1+z_{1-\beta})\sqrt{\log (2p/t_1)}}{c^*\norm{\nabla x^*(s^*)}_2}\big(\sqrt{\frac{\alpha}{n_2}}+\frac{\sqrt{\log (2p/t_1)}}{n_2}\big)< \frac{\delta}{2}.
\end{equation}
We explain why \eqref{choice of p,t1,n2 CCP} implies \eqref{condition for lower bound of H}. Assumption \ref{unique optimal parameter + local differentiability} stipulates that $1\geq H(x^*(s^*+\delta))\geq H(x^*_S)+\frac{1}{2}\nabla H(x^*_S)'\nabla x^*(s^*)\delta=1-\alpha+\frac{1}{2}c^*\norm{\nabla x^*(s^*)}_2\delta$, which leads to $c^*\norm{\nabla x^*(s^*)}_2\delta\leq 2\alpha$. The second inequality in \eqref{choice of p,t1,n2 CCP} then gives $(1+z_{1-\beta})\sqrt{\log (2p/t_1)}\big(\sqrt{\frac{\alpha}{n_2}}+\frac{\sqrt{\log (2p/t_1)}}{n_2}\big)<\alpha/6<\sqrt{2}\alpha/4$. Similar to the proof of Theorem \ref{finite sample error:gaussian+differentiable constraint}, when \eqref{choice of p,t1,n2 CCP} holds and \eqref{maximal deviation CCP:sample mean adjusted by std} happens, we must have $\hat s^*\in[s^*-2\epsilon(p,t_1,n_2),s^*+2\epsilon(p,t_1,n_2)]$. Therefore under the condition \eqref{choice of p,t1,n2 CCP}
\begin{equation}\label{maximal deviation CCP:parameter}
   P_{\bm{\xi}_{1:n_2}}\big( \lvert \hat s^*-s^*\rvert>2\epsilon(p,t_1,n_2)\big)\leq 2t_1.
\end{equation}

Now we proceed to deal with the finite sample confidence error. Following the same steps of bounding the feasibility confidence level, we have
\begin{eqnarray*}
&&P_{\bm{\xi}_{1:n_2}}(H(x^*(\hat s^*))\geq \gamma)\\
&\geq&P_{\bm{\xi}_{1:n_2}}\big(\frac{\sqrt{n_2}(\hat H(x^*_S)-H(x^*_S))}{\sigma(x^*_S)}+\frac{\sqrt{n_2}}{\sigma(x^*_S)}\max_{j:\lvert s_j-s^*\rvert\leq 2\epsilon(p,t_1,n_2)}\big\lvert\hat H_j-H_j-(\hat H(x^*_S)-H(x^*_S))\big\rvert+\\
&&\hspace{9ex}\frac{z_{1-\beta}}{\sigma(x^*_S)}\max_{j:\lvert s_j-s^*\rvert\leq 2\epsilon(p,t_1,n_2)}\lvert \sigma(x^*_S)-\hat \sigma_j\rvert\leq z_{1-\beta},\ \lvert \hat s^*-s^*\rvert\leq 2\epsilon(p,t_1,n_2)\big).
\end{eqnarray*}
We bound the deviation of sample standard deviation as follows
\begin{eqnarray*}
\max_{j:\lvert s_j-s^*\rvert\leq 2\epsilon(p,t_1,n_2)}\lvert \sigma(x^*_S)-\hat \sigma_j\rvert&\leq&\max_{j:\lvert s_j-s^*\rvert\leq 2\epsilon(p,t_1,n_2)}\lvert \sigma(x^*_S)-\sigma_j\rvert+\max_{j=1,\ldots,p}\lvert \sigma_j-\hat \sigma_j\rvert\\
&\leq&\max_{j:\lvert s_j-s^*\rvert\leq 2\epsilon(p,t_1,n_2)}\lvert \sqrt{\alpha(1-\alpha)}-\sqrt{H_j(1-H_j)}\rvert+\max_{j=1,\ldots,p}\lvert \sigma_j-\hat \sigma_j\rvert.
\end{eqnarray*}
The second error is taken care of by \eqref{maximal deviation CCP:std}. To bound the first error, we note that by Assumption \ref{unique optimal parameter + local differentiability} we have $\max_{j:\lvert s_j-s^*\rvert\leq 2\epsilon(p,t_1,n_2)}\lvert H_j-(1-\alpha)\rvert\leq 2c^*\norm{\nabla x^*(s^*)}_2\cdot 2\epsilon(p,t_1,n_2)=24(1+z_{1-\beta})\sqrt{\log (2p/t_1)}\big(\sqrt{\frac{\alpha}{n_2}}+\frac{\sqrt{\log (2p/t_1)}}{n_2}\big)$. Therefore if $24(1+z_{1-\beta})\sqrt{\log (2p/t_1)}\big(\sqrt{\frac{\alpha}{n_2}}+\frac{\sqrt{\log (2p/t_1)}}{n_2}\big)\leq \alpha/2$, it follows from applying mean value theorem that
\begin{equation*}
    \max_{j:\lvert s_j-s^*\rvert\leq 2\epsilon(p,t_1,n_2)}\lvert \sqrt{\alpha(1-\alpha)}-\sqrt{H_j(1-H_j)}\rvert\leq \frac{1}{\sqrt{\alpha}}24(1+z_{1-\beta})\sqrt{\log (2p/t_1)}\big(\sqrt{\frac{\alpha}{n_2}}+\frac{\sqrt{\log (2p/t_1)}}{n_2}\big).
\end{equation*}
Similar to \eqref{deviation:std local}, we can argue that $24(1+z_{1-\beta})\sqrt{\log (2p/t_1)}\big(\sqrt{\frac{\alpha}{n_2}}+\frac{\sqrt{\log (2p/t_1)}}{n_2}\big)\leq \alpha/2$ can be assumed without loss of generality so that the above bound can be assumed to hold. Together with \eqref{maximal deviation CCP:std}, we have
\begin{eqnarray}
\notag&&\max_{j:\lvert s_j-s^*\rvert\leq 2\epsilon(p,t_1,n_2)}\lvert \sigma(x^*_S)-\hat \sigma_j\rvert\\
&\leq&\frac{24}{\sqrt{\alpha}}(1+z_{1-\beta})\sqrt{\log (2p/t_1)}\big(\sqrt{\frac{\alpha}{n_2}}+\frac{\sqrt{\log (2p/t_1)}}{n_2}\big)+\sqrt{\frac{2\log(2p/t_1)}{n_2-1}}\text{ for all }j=1,\ldots,p\label{deviation CCP:local std}
\end{eqnarray}
with probability at least $1-t_1$. Now we can further bound the confidence level
\begin{eqnarray*}
&&P_{\bm{\xi}_{1:n_2}}(H(x^*(\hat s^*))\geq \gamma)\\
&\geq&P_{\bm{\xi}_{1:n_2}}\Big(\frac{\sqrt{n_2}(\hat H(x^*_S)-H(x^*_S))}{\sigma(x^*_S)}+\frac{1}{\sigma(x^*_S)}\big(\sqrt{2\sigma^2_{2\epsilon(p,t_1,n_2)}\log(2p/t_2)}+\frac{2\log(2p/t_2)}{3\sqrt{n_2}}\big)\\
&&\hspace{5ex}+\frac{24z_{1-\beta}}{\sigma(x^*_S)\sqrt{\alpha}}(1+z_{1-\beta})\sqrt{\log (2p/t_1)}\big(\sqrt{\frac{\alpha}{n_2}}+\frac{\sqrt{\log (2p/t_1)}}{n_2}\big)+\frac{z_{1-\beta}}{\sigma(x^*_S)}\sqrt{\frac{2\log(2p/t_1)}{n_2-1}}\leq z_{1-\beta}\Big)\\
&&\hspace{5ex}-t_2-3t_1\text{\ \ \ by \eqref{equicontinuity of local empirical process CCP}, \eqref{maximal deviation CCP:parameter} and \eqref{deviation CCP:local std}}.
\end{eqnarray*}
Like in the proof of Theorem \ref{finite sample error:gaussian+differentiable constraint}, applying Berry-Esseen theorem to the first probability on the right hand side and absorbing various constants into the universal constant $C$ give
\begin{eqnarray*}
&&1-\beta-P_{\bm{\xi}_{1:n_2}}(H(x^*(\hat s^*))\geq \gamma)\\
&\leq&C\Big(\frac{1}{\sqrt{\alpha n_2}}+\frac{1}{\sqrt{\alpha}}\big(\sqrt{2\sigma^2_{2\epsilon(p,t_1,n_2)}\log(2p/t_2)}+\frac{2\log(2p/t_2)}{3\sqrt{n_2}}\big)\\
&&\hspace{5ex}+\frac{z_{1-\beta}}{\alpha}(1+z_{1-\beta})\sqrt{\log (2p/t_1)}\big(\sqrt{\frac{\alpha}{n_2}}+\frac{\sqrt{\log (2p/t_1)}}{n_2}\big)+\frac{z_{1-\beta}}{\sqrt{\alpha}}\sqrt{\frac{2\log(2p/t_1)}{n_2}}+t_1+t_2\Big)\\
&\leq&C\Big(\frac{1}{\sqrt{\alpha}}\big(\sqrt{\sigma^2_{2\epsilon(p,t_1,n_2)}\log(2p/t_2)}+\frac{\log(2p/t_2)}{\sqrt{n_2}}\big)+(1+z_{1-\beta})^2\sqrt{\frac{\log(2p/t_1)}{\alpha n_2}}+t_1+t_2\Big)\\
&&\text{\ \  by keeping dominant terms only}\\
&\leq&C\Big(\frac{1}{\sqrt{\alpha}}\big(\sqrt{\sigma^2_{2\epsilon(p,1/n_2,n_2)}\log(pn_2)}+\frac{\log(pn_2)}{\sqrt{n_2}}\big)+(1+z_{1-\beta})^2\sqrt{\frac{\log(pn_2)}{\alpha n_2}}+\frac{1}{n_2}\Big)\\
&&\text{\ \  by taking }t_1=t_2=\frac{2}{n_2}\\
&\leq&C\Big(\sqrt{\frac{\log(pn_2)}{\alpha}}\sigma_{2\epsilon(p,2/n_2,n_2)}+(1+z_{1-\beta})^2\frac{\log(pn_2)}{\sqrt{\alpha n_2}}\Big)\\
&\leq&C\Big(\sqrt{\frac{\log(pn_2)}{\alpha}}\sigma_{2\epsilon(p,2/n_2,n_2)}+(1+z_{1-\beta})^2\frac{\log(pn_2)}{\sqrt{\alpha n_2}}\Big).
\end{eqnarray*}
It remains to bound the $\sigma_{2\epsilon(p,2/n_2,n_2)}$ term which by the definition \eqref{bound for local envelope CCP} can be expressed as
\begin{equation*}
    \sigma_{2\epsilon(p,2/n_2,n_2)} = \sqrt{20K\tilde \epsilon(p,n_2)}\big(\log(\max\{e,\frac{1}{2\tilde \epsilon(p,n_2)}\})\big)^{\frac{1}{4}}
\end{equation*}
where
\begin{equation*}
\tilde \epsilon(p,n_2)=\frac{6(1+z_{1-\beta})D_2^2D_3\sqrt{\log(2K/\alpha)}}{c^*\min_k\lvert b_k\rvert}\big( \sqrt{\frac{\alpha \log(pn_2)}{n_2}}+\frac{\log(pn_2)}{n_2}\big).
\end{equation*}
Note that $\frac{1}{2\tilde \epsilon(p,n_2)}\leq \frac{c^*\min_k\lvert b_k\rvert  n_2}{D_2^2D_3\sqrt{\log(2K/\alpha)}}=\frac{n_2}{\tilde C}$, hence using this upper bound in the logarithm we have
\begin{eqnarray*}
    \sigma_{2\epsilon(p,2/n_2,n_2)} &\leq&11\sqrt{(1+z_{1-\beta})K\tilde C}\big(\log(\max\big\{e,\frac{n_2}{\tilde C}\big\})\big)^{\frac{1}{4}}\big( \sqrt{\frac{\alpha \log(pn_2)}{n_2}}+\frac{\log(pn_2)}{n_2}\big)^{\frac{1}{2}}\\
    &\leq&11\sqrt{(1+z_{1-\beta})K\tilde C}\big(\log(\max\big\{e,\frac{n_2}{\tilde C}\big\})\big)^{\frac{1}{4}}\big[\big(\frac{\alpha \log(pn_2)}{n_2}\big)^{\frac{1}{4}}+\big(\frac{\log(pn_2)}{n_2}\big)^{\frac{1}{2}}\big]
\end{eqnarray*}
where the second inequality follows because $\sqrt{a+b}\leq \sqrt{a}+\sqrt{b}$ for any $a,b\geq 0$. Substituting $\sigma_{2\epsilon(p,2/n_2,n_2)}$ with its upper bound gives
\begin{eqnarray*}
&&1-\beta-P_{\bm{\xi}_{1:n_2}}(H(x^*(\hat s^*))\geq \gamma)\\
&\leq&C\Big(\big[(1+z_{1-\beta})^2+\sqrt{(1+z_{1-\beta})K\tilde C}\big(\log(\max\big\{e,\frac{n_2}{\tilde C}\big\})\big)^{\frac{1}{4}}\big]\frac{\log(pn_2)}{\sqrt{\alpha n_2}}\\
&&\hspace{5ex}+\sqrt{(1+z_{1-\beta})K\tilde C}\big(\log(\max\big\{e,\frac{n_2}{\tilde C}\big\})\big)^{\frac{1}{4}}\frac{(\log(pn_2))^{3/4}}{(\alpha n_2)^{1/4}}\Big)\\
&\leq&C(1+z_{1-\beta})^2\big(1+\sqrt{\tilde C K}\big(\log(\max\big\{e,\frac{n_2}{\tilde C}\big\})\big)^{\frac{1}{4}}\big)\big(\frac{(\log(pn_2))^{3/4}}{(\alpha n_2)^{1/4}}+\frac{\log(pn_2)}{\sqrt{\alpha n_2}}\big)\\
&\leq&C(1+z_{1-\beta})^2\big(1+\sqrt{\tilde C K}\big(\log(\max\big\{e,\frac{n_2}{\tilde C}\big\})\big)^{\frac{1}{4}}\big)\frac{(\log(pn_2))^{3/4}}{(\alpha n_2)^{1/4}}
\end{eqnarray*}
where the last inequality follows because $\frac{(\log(pn_2))^{3/4}}{(\alpha n_2)^{1/4}}\leq \frac{\log(pn_2)}{\sqrt{\alpha n_2}}$ if $\frac{(\log(pn_2))^{3/4}}{(\alpha n_2)^{1/4}}\leq 1$. Note again that this bound is valid when \eqref{choice of p,t1,n2 CCP} is satisfied at $t_1=\frac{2}{n_2}$. Replacing $e$, the base of the natural logarithm, with $3$ gives the desired bound.\hfill\Halmos

\section{Applying Univariate Gaussian Validator to Formulations with Multidimensional Conservativeness Parameters}\label{sec:multidimensional parameter univariate gaussian}
We consider the case of multidimensional conservativeness parameter, i.e., $S\subset \R^q$ for some $q\geq 2$, and present the asymptotic performance guarantees of the univariate Gaussian validator. We assume the following counterpart of Assumption \ref{piecewise uniform continuous solution curve}:
\begin{assumption}[Piecewise uniformly continuous solution curve]\label{piecewise uniform continuous solution curve:multivariate}
The parameter space $S\subset \R^q$ is compact, and there exist $M$ connected and open subsets $S_1,\ldots,S_M$ of $S$ such that (i) $S_i\cap S_{i'}=\emptyset$ for all $i\neq i'$; (ii) $\mathrm{m}(\cup_{i=1}^MS_i)=\mathrm{m}(S)$ where $\mathrm{m}(\cdot)$ denotes the Lebesgue measure on $\R^q$; and (iii) for each $i=1,\ldots,M$, the optimal solution $x^*(s)$ of $OPT(s)$ exists and is unique for all $s\in S_i$, and $x^*(s)$ as a function of $s$ is uniformly continuous on $S_i$.
\end{assumption}

Similar to the case of scalar parameter, the solution curve $x^*(s)$ on each piece $S_i$ can be continuously extended to the closure $\overline{S_i}:=\cap_{S'\text{ is closed, }S_i\subseteq S'}S'$ under this piecewise uniform continuity assumption. Specifically, for every parameter value $s\in S\backslash \cup_{i=1}^MS_i$, we define the extended parameter-to-solution mapping to be
\begin{equation*}
x^*(s):=\{\lim_{s'\in S_i,s'\to s}x^*(s'):s\in \overline{S_i},i=1,\ldots,M\}.
\end{equation*}
Accordingly, the optimal solution set and optimal parameter set associated with the solution path are defined as
\begin{equation}\label{optimal solution:multidimensional}
    \mathcal X_S^*:=\text{argmin}\{f(x):H(x)\geq \gamma,x=x^*(s)\text{ for }s\in\cup_{i=1}^MS_i\text{ or }x\in x^*(s)\text{ for some }s\in S\backslash \cup_{i=1}^MS_i\}
\end{equation}
and
\begin{equation*}
    S^*:=\{s\in\cup_{i=1}^MS_i:x^*(s)\in \mathcal X_S^*\}\cup\{s\in S\backslash \cup_{i=1}^MS_i:x^*(s) \cap \mathcal X_S^*\neq \emptyset\}.
\end{equation*}


We also assume uniqueness of the optimal solution:
\begin{assumption}[Unique optimal solution]\label{unique optimal solution:multivariate}
The optimal solution set $\mathcal X_S^*$ defined in \eqref{optimal solution:multidimensional} is a singleton $\{x_S^*\}$.
\end{assumption}
Note that in the case of scalar $s$, uniqueness of the optimal solution is a consequence (Proposition \ref{unique optimum}) of several more elementary assumptions among which monotonicity of the robust feasible set with respect to the parameter (Assumption \ref{monotonicity}) plays the key role. However, such notion of monotonicity does not completely carry to the mutidimensional case. For example, one may have a formulation $OPT(s)$ such that the robust feasible set satisfies $\mathrm{Sol}(s)\subseteq \mathrm{Sol}(s')$ whenever $s'\leq s$ component-wise, but $\mathrm{Sol}(s)$ and $\mathrm{Sol}(s')$ are in general not comparable.

We also assume the following counterpart of Assumption \ref{unique optimal solution}:
\begin{assumption}\label{non-degenerate solution path}
For every $\epsilon>0$ there exists an $s\in \cup_{i=1}^MS_i$ such that $H(x^*(s))>\gamma$ and $\Vert x^*(s)-x_S^*\Vert_2<\epsilon$, where $x_S^*$ is the unique optimal solution from Assumption \ref{unique optimal solution:multivariate}.
\end{assumption}
We then have the following asymptotic performance guarantees for Algorithm \ref{algo:marginal}:
\begin{theorem}[Asymptotic joint feasibility$+$optimality guarantee]\label{multivariate asymptotic joint:univariate gaussian}
Suppose Assumptions \ref{continuous obj}-\ref{non-degenerate variance} hold for \eqref{stoc_opt}. Also suppose that Assumptions \ref{piecewise uniform continuous solution curve:multivariate}-\ref{non-degenerate solution path} hold for the formulation $OPT(s)$, and that $\{s_1,\ldots,s_p\}\subseteq \cup_{i=1}^MS_i$. Denote by $\epsilon_S:=\sup_{s\in S}\inf_{j=1,\ldots,p}\norm{s-s_j}_2$ the mesh size, and by $x_S^*$ be the unique optimal solution from Assumption \ref{unique optimal solution:multivariate}. Conditional on Phase one, as Phase two data size $n_2\to\infty$, we have for the output of Algorithm \ref{algo:marginal} that (i) $\lim_{n_2\to\infty,\epsilon_S\to 0}x^*(\hat s^*)=x_S^*$ and $\lim_{n_2\to\infty,\epsilon_S\to 0}d(\hat s^*,\mathcal S^*)=0$ almost surely; and (ii) $\liminf_{n_2\to\infty,\epsilon_S\to 0}P_{\bm{\xi}_{1:n_2}}(H(x^*(\hat s^*))\geq \gamma)\geq 1-\beta$ if $H(x_S^*)=\gamma$, and $\lim_{n_2\to\infty,\epsilon_S\to 0}P_{\bm{\xi}_{1:n_2}}(H(x^*(\hat s^*))\geq \gamma)= 1$ if $H(x_S^*)>\gamma$.
\end{theorem}
\proof{Proof of Theorem \ref{multivariate asymptotic joint:univariate gaussian}.}The proof is the same as that of Theorem \ref{asymptotic joint:gaussian} with straightforward modifications. In particular, $\{\tilde s_1,\ldots,\tilde s_M\}$ shall be replaced by $S\backslash\cup_{i=1}^MS_i$ and the solution set $\mathcal X_S$ is now defined as $\mathcal X_S:=\{x^*(s):s\in\cup_{i=1}^MS_i\}\cup\big(\cup_{s\in S\backslash\cup_{i=1}^MS_i}x^*(s)\big)$.\hfill\Halmos

In order to establish an asymptotically tight feasibility confidence level like in Theorem \ref{asymptotic tight coverage:gaussian}, we further assume uniqueness of the optimal parameter:
\begin{assumption}[Unique optimal parameter]\label{unique optimal parameter:multivariate}
The optimal parameter set $S^*$ is a singleton $\{s^*\}$, and $s^*\in S_{i^*}$ for some $i^*=1,\ldots,M$.
\end{assumption}
We then have the following guarantee:
\begin{theorem}[Asymptotically tight feasibility guarantee]\label{asymptotic tight coverage:gaussian multidimensional}
In addition to the conditions of Theorem \ref{multivariate asymptotic joint:univariate gaussian}, suppose Assumption \ref{unique optimal parameter:multivariate} holds. Suppose also that the parameter-to-objective mapping $v(s)$ satisfies $v(s)<v(s')$ whenever $s<s'$ component-wise and that $H(x_S^*)=\gamma$. For each $j=1,\ldots,p$, let
\begin{equation*}
\tilde j:=\argmin_{j'}\{\norm{s_j-s_{j'}}_2:s_{j'}<s_j\text{ component-wise},s_{j'}\text{ lies on the same piece as }s_j\}
\end{equation*}
and if there is no such feasible $j'$ simply let $\tilde j:=j$. If the mesh is such that
\begin{equation}\label{mesh resolution:multidimension}
\max_{j=1,\ldots,p}\lvert H(x^*(s_j))-H(x^*(s_{\tilde j}))\rvert=o\big(\frac{1}{\sqrt{n_2}}\big)
\end{equation}
then we have for the output of Algorithm \ref{algo:marginal} that $\lim_{n_2\to\infty,\epsilon_S\to 0\text{ s.t. \eqref{mesh resolution:multidimension} holds}}P_{\bm{\xi}_{1:n_2}}(H(x^*(\hat s^*))\geq \gamma)= 1-\beta$.
\end{theorem}
\proof{Proof of Theorem \ref{asymptotic tight coverage:gaussian multidimensional}.}The proof follows exactly that of Theorem \ref{asymptotic tight coverage:gaussian} with straightforward modifications. For example, when bounding $\hat H(x^*(\hat s^*))-z_{1-\beta}\frac{\hat\sigma(\hat s^*)}{\sqrt{n_2}}$ in the proof of Theorem \ref{asymptotic tight coverage:gaussian} we replace the parameter value $s^{i^*}_{j^*}$ output by the algorithm with $s^{i^*}_{j^*-1}$ and use the condition \eqref{high resolution}, whereas now we shall replace the output parameter value $\hat s^*=s_{j^*}$ with $s_{\tilde{j^*}}$ and then use \eqref{mesh resolution:multidimension} to obtain the same bound.\hfill\Halmos

\end{document}